\documentclass[11pt, a4paper]{article}
\usepackage{amsmath,amsthm,amssymb,amsfonts}
\usepackage{a4wide}
\newtheorem{theorem}{Theorem}[section]
\newtheorem{lemma}[theorem]{Lemma}
\newtheorem{proposition}[theorem]{Proposition}
\newtheorem{definition}[theorem]{Definition}

\newtheorem{corollary}[theorem]{Corollary}

\theoremstyle{remark}
\newtheorem{example}[theorem]{Example}
\newtheorem{remark}[theorem]{Remark}

\def\Xint#1{\mathchoice
{\XXint\displaystyle\textstyle{#1}}%
{\XXint\textstyle\scriptstyle{#1}}%
{\XXint\scriptstyle\scriptscriptstyle{#1}}%
{\XXint\scriptscriptstyle\scriptscriptstyle{#1}}%
\!\int}
\def\XXint#1#2#3{{\setbox0=\hbox{$#1{#2#3}{\int}$ }
\vcenter{\hbox{$#2#3$ }}\kern-.6\wd0}}

\def\dashint{\Xint-}

\usepackage{graphicx}
\usepackage[usenames,dvipsnames]{color}
\usepackage[colorlinks=true, pdfstartview=FitV, linkcolor=black, citecolor=blue, urlcolor=blue]{hyperref}

\makeatletter

\@addtoreset{equation}{section}
\makeatother

\newcommand{\E}{{\mathbb E}}

\newcommand{\T}{{\mathbb T}}
\newcommand{\N}{{\mathbb N}}
\newcommand{\Z}{{\mathbb Z}}

\newcommand{\R}{{\mathbb R}}

\newcommand{\mO}{\mathcal{O}}
\newcommand{\mV}{\mathcal{V}}
\newcommand{\mW}{\mathcal{W}}
\newcommand{\mI}{\mathcal{I}}

\newcommand{\pa}{{\partial}}
\newcommand{\na}{{\nabla}}

\newcommand{\eps}{{\varepsilon}}

\def\curl{\hbox{curl \!}}
\def\div{\hbox{div \!}}

\def\ext{{\rm ext}  \,}
\def\Id{{\rm Id}}

\def\sspace{\smallskip \noindent}
\def\mspace{\medskip \noindent}

\title{Analysis of the viscosity of dilute suspensions \\ beyond Einstein's formula} 
\author{David G\'erard-Varet\footnote{Universit\'e Paris Diderot, Sorbonne Paris Cit\'e, Institut de Math\'ematiques de Jussieu-Paris Rive Gauche (UMR 7586), F-75205 Paris, France. Email:  {\textrm  david.gerard-varet@imj-prg.fr}} \hspace{0.5cm} Matthieu Hillairet\footnote{Universit\'e de Montpellier, Institut  Montpelli\'erain Alexandre Grothendieck (UMR5149), 34090 Montpellier, France. Email:  {\textrm  mhillairet@umontpellier.fr} }}

\begin{document}
\maketitle

\begin{abstract}
We provide a mathematical analysis of the effective viscosity of suspensions of spherical particles in a Stokes flow, at low solid volume fraction $\phi$. Our objective is to go beyond the Einstein's approximation $\mu_{eff} = (1+\frac{5}{2}\phi) \mu$. Assuming a lower bound on the minimal distance between the $N$ particles, we are able to identify the $O(\phi^2)$ correction to the effective viscosity, which involves pairwise particle interactions. Applying the methodology developped over the last years on Coulomb gases, we are able to tackle the limit $N \rightarrow +\infty$ of  the $O(\phi^2)$-correction, and provide explicit formula for this limit when the particles centers can be  described by either periodic or stationary ergodic point processes. 
\end{abstract}

\section{Setting of the problem}
Our general concern is the computation of the effective viscosity generated by a suspension of $N$ particles in a fluid flow. We consider  spherical particles of small radius $a$,  centered at $x_{i,N}$, with $N \ge 1$ and $1\le i\le N$.  To lighten notations, we write $x_i$ instead of $x_{i,N}$, and  $B_i = B(x_i, a)$. We assume that the Reynolds number of the fluid flow is small, so that the fluid velocity is governed by the Stokes equation. Moreover, the particles are assumed to be force- and torque-free.  If 
$\mathcal{F} =  \R^3 \setminus  (\cup_i B_i)$ is the fluid domain, governing equations are 
\begin{equation}
\label{Sto}
\left\{
\begin{aligned}
-\mu \Delta u + \na  p & = 0, \quad x \in \mathcal{F} , \\ 
\div u & = 0, \quad x \in \mathcal{F} , \\
u\vert_{B_i} & = u_i  + \omega_i \times (x-x_i), 
\end{aligned}
\right.
\end{equation}
where $\mu$ is the kinematic viscosity, while the constant vectors $u_i$ and $\omega_i$ are Lagrange multipliers associated to the constraints 
\begin{equation}
\label{Sto2}
\begin{aligned}
\int_{\pa B_i} \sigma_\mu(u,p) n\, ds& = 0, \quad & \int_{\pa B_i} \sigma_\mu(u,p) n \times (x-x_i) \, ds  = 0.
\end{aligned}
\end{equation}
Here,  $\sigma_\mu(u,p) := 2\mu D(u) - p I$ is the usual Cauchy stress tensor. The boundary condition at infinity will be specified later on.  

\mspace
We are interested in a situation where the number of particles is large, $N \gg 1$.  We want to understand the additional viscosity created by the particles. Ideally, our goal is to replace the viscosity coefficient $\mu$ in  \eqref{Sto} by an effective viscosity tensor $\mu'$ that would encode the average effect induced by the particles. Note that such replacement can only make sense in the flow region $\mO$ in which the particles are distributed in a dense way. For instance, a finite number of isolated particles will not contribute to the effective viscosity, and should not be taken into account in  $\mO$. The selection of the flow region  is formalized through the following hypothesis on the empirical measure: 
\begin{equation}  \label{H1} \tag{H1}
\begin{aligned}
  & \delta_N = \frac{1}{N} \sum_{i=1}^N \delta_{x_i} \xrightarrow[N \rightarrow +\infty]{} f(x) dx \quad \text{weakly,} \\
& support(f)  = \overline{\mO}, \quad \mO \text{ smooth, bounded and open}, \: f\vert_\mO \in C^1(\overline{\mO}).
\end{aligned}
\end{equation}
Note that we do not ask for regularity of the limit density $f$ over $\R^3$, but only in restriction to $\mO$. Hence, our assumption covers the important case $f = \frac{1}{|\mO|} 1_{\mO}$. 

\mspace
We  investigate the classical regime of dilute suspensions, in which the solid volume fraction 
\begin{equation} \label{def_phi}
\phi = \frac{4}{3}N \pi a^3/|\mO|
\end{equation}
is small, but independent of $N$.  Besides \eqref{H1}, we  make  the separation hypothesis   
 \begin{equation} \label{H2} \tag{H2} 
\min_{i\neq j} |x_i - x_j| \ge c N^{-1/3} \: \text{ for some constant $c > 0$ independent of $N$} .  
\end{equation}
Let us stress that \eqref{H2} is compatible with \eqref{H1} only if the $L^\infty$ norm of $f$ is small enough  (roughly less than $1/c^3$), which in turn forces $\mO$ to be large enough.

\mspace
Our  hope is to replace a model of type \eqref{Sto} by a model of the form 
\begin{equation}
\label{Stoeff}
\left\{
\begin{aligned}
-\mu \Delta u + \na  p & = 0, \quad \div u = 0,  \quad x \in  \R^3 \setminus\mO, \\ 
-2 \div (\mu' D(u')) + \na  p' & = 0,  \quad \div u' = 0, \quad x \in \mO, 
\end{aligned}
\right.
\end{equation}
with the usual continuity conditions on the velocity and the stress: 
\begin{equation} \label{jumpBC}
u = u' \quad \text{at } \: \pa \mO, \quad \sigma_\mu(u,p)n=  \sigma_{\mu'}(u',p')n \quad \text{at } \: \pa \mO. 
\end{equation}
{\it A priori}, $\mu'$ could be inhomogeneous (and should be if the density $f$ seen above is itself non-constant over $\mO$). It could also be anisotropic, if the cloud of particles favours some direction. With this in mind, it is natural to look for   $\mu' = \mu'(x)$ as a general 4-tensor, with $\sigma' = 2 \mu' D(u)$ given in coordinates by $\sigma_{ij} = \mu'_{ij kl} D(u)_{kl}$. By standard classical considerations of mechanics, $\mu'$ should satisfy the relations 
$$  \mu'_{ij kl}  =  \mu'_{jikl} = \mu'_{jilk} = \mu'_{lkji}, $$
namely $\mu'$ should define a  symmetric isomorphism over the space of $3\times3$ symmetric matrices.

\mspace
As we consider a situation in which $\phi$ is small, we may expect $\mu'$ to be a small perturbation of $\mu$, and hopefully admit an expansion in powers of $\phi$: 
\begin{equation} \label{asymptmu}
 \mu' = \mu \Id + \phi \mu_1  + \phi^2 \mu_2 + \dots + \phi^k \mu_k + o(\phi^k).
 \end{equation}
The main mathematical questions are: 
\begin{itemize}
\item Can solutions $u_N$ of \eqref{Sto}-\eqref{Sto2} be approximated by solutions $u_{eff} = 1_{\R^3\setminus\mO} u + 1_{\mO} u'$ of \eqref{Stoeff}-\eqref{jumpBC}, for an appropriate choice of $\mu'$  and an appropriate topology ? 
\item If so, does $\mu'$ admit an expansion of type \eqref{asymptmu}, for some $k$ ? 
\item If so, what are the values of the viscosity coefficients  $\mu_i$, $1 \le i \le k$  ?    
\end{itemize}
Let us stress that,  in most articles about the effective viscosity of suspensions, it is implicitly assumed that the first two questions have a positive answer, at least for $k=1$ or $2$. In other words, the existence of an effective model is taken for granted, and the point is then to answer the third question, or at least to determine the mean values 
\begin{equation} \label{formula_nu_i}
\nu_i := \frac{1}{|\mO|} \int_\mO \mu_i(x) dx
\end{equation}
 of the viscosity coefficients. As we will see in Section \ref{sec1}, these mean values can be determined from the asymptotic behaviour of some integral quantities $\mathcal{I}_N$ as $N \rightarrow +\infty$. These integrals  involve the solutions $u_N$ of \eqref{Sto}-\eqref{Sto2} with condition at infinity
\begin{equation} \label{BC:infinity}
\lim_{|x| \rightarrow +\infty} u(x) - S x = 0 
\end{equation}
where $S$ is an arbitrary symmetric trace-free matrix.

\mspace
The effective viscosity problem for dilute suspensions of spherical particles has a long history, mostly focused on the first order correction created by the suspension, that is $k=1$ in \eqref{asymptmu}. 
The pioneering work on this problem was due to Einstein \cite{Ein}, not mentioning earlier contributions  on the similar conductivity problem  by Maxwell \cite{Max}, Clausius \cite{Cla}, Mossotti \cite{Mos}. 
The celebrated Einstein's formula, 
\begin{equation}  \label{Einstein}
\mu' = \mu + \frac{5}{2} \phi \mu + o(\phi). 
\end{equation}
was derived under the assumption that the particles are homogeneously and isotropically distributed, and neglecting the interactions between particles. In other words, the correction 
$\mu_1 = \frac{5}{2} \mu$ is obtained by summing $N$ times the contribution of one spherical particle to the effective stress. The calculation of Einstein will be seen in Section \ref{sec1}.  It was later extended to the case of an inhomogeneous suspension by Almog and Brenner \cite[page 16]{AlBr}, who found 
\begin{equation}  \label{Almog-Brenner}
\mu_1  =  \frac{5}{2} |\mO| f(x) \mu. 
\end{equation}
The  mathematical justification of  formula \eqref{Einstein} came much later. As far as we know, the first step in this direction was due to Sanchez-Palencia \cite{MR813656} and Levy and Sanchez-Palencia \cite{MR813657}, who recovered Einstein's formula from homogenization techniques, when the suspension is periodically distributed in a bounded domain. Another justification, based on variational principles, is due to Haines and Mazzucato \cite{MR2982744}. They also consider a periodic array of spherical particles in a bounded domain $\Omega$, and define the viscosity coefficient of the suspension in terms of the energy dissipation rate: 
$$ \mu_N  = \frac{\mu}{ |S|^2} \int_{\mathcal{F}} |D(u_N)|^2  $$
where $u_N$ is the solution of \eqref{Sto}-\eqref{Sto2}-\eqref{BC:infinity}, replacing $\R^3$ by $\Omega$. Their main result is that  
$$\mu_N = \mu + \frac{5}{2} \phi \mu + O(\phi^{3/2}).$$
For  preliminary results in the same spirit, see Keller-Rubenfeld \cite{KeRu}. Eventually, a recent work \cite{HiWu} by the second author and Di Wu shows the validity of Einstein's formula under general assumptions of type \eqref{H1}-\eqref{H2}. 
 See also \cite{NiSc} for a similar recent result.

\mspace
Our goal in the present paper is to go beyond this famous formula, and to study the second order correction to the effective viscosity, that is $k=2$ in \eqref{asymptmu}. Results on this problem have split so far into two settings: periodic distributions, and random distributions of spheres. Many different formula have emerged in the literature, after analytical, numerical and experimental studies. In the periodic case, one can refer to the works \cite{Sai,ZuAdBr,NuKe,MR3025042}, or to the more recent work \cite{MR3025042}, dedicated to the case of spherical inclusions of another Stokes fluid with viscosity $\tilde \mu \neq \mu$. Still, in the simple case of a primitive cubic lattice, the expressions for the second order correction differ.  In the random case, the most reknowned analysis is due to Batchelor and Green \cite{BaGr1}, who consider a homogeneous and stationary distribution of spheres, and express the correction $\mu_2$ as an ensemble average that involves the $N$-point correlation function of the process. As pointed out by Batchelor and Green, the natural idea when investigating the effective viscosity up to $O(\phi^2)$  is to replace  the $N$-point correlation function by the $2$-point correlation function, but this leads to a divergent integral. To overcome this difficulty, Batchelor and Green develop what they call a renormalization technique, that was developed earlier by Batchelor to determine the sedimentation speed of a dilute suspension. 
After further analysis of the expression of the two-point correlation function of spheres in a Stokes flow \cite{BaGr2}, completed by numerical computations, they claim that under a pure strain, the particles induce a viscosity of the form 
\begin{equation}  \label{Batchelor-Green}
\mu' =  \mu + \frac{5}{2} \phi \mu + 7.6 \phi^2 \mu + o(\phi^2)
\end{equation}
Although the result of Batchelor and Green is generally accepted by the fluid mechanics community, the lack of clarity about their renormalization technique has led to  debates, see \cite{Hin, Obr, AlBr}. 

\mspace
One main objective in the present paper is to give a rigorous and global  mathematical framework for the computation of 
\begin{equation} \label{def:nu2}
\nu_2 = \frac{1}{|\mO|} \int_\mO \mu_2(x) dx 
\end{equation} 
leading to explicit formula in  periodic and stationary random settings. We will adopt the point of view of the studies mentioned before: we will assume the validity of an effective model of type \eqref{Stoeff}-\eqref{jumpBC}-\eqref{asymptmu} with $k=2$, and will identify the averaged coefficient $\nu_2$.  

\mspace
More precisely, our analysis divides into two parts. The first part, carried in Section \ref{sec1}, has as its main consequence the following 
\begin{theorem} \label{thm1}
Let $(x_{i})_{1 \le i \le N}$  a family of points supported in a fixed compact of  $\R^3$, and satisfying \eqref{H1}-\eqref{H2}. For any trace-free symmetric matrix $S$ and any $\phi > 0$, let $u_N$, resp. $u_{eff}$, the solution of \eqref{Sto}-\eqref{Sto2}-\eqref{BC:infinity} with the radius $a$ of the balls defined through \eqref{def_phi}, resp. the solution  of \eqref{Stoeff}-\eqref{jumpBC}-\eqref{BC:infinity} where $\mu'$ obeys \eqref{asymptmu} with $k=2$, $\mu_1$ being given in \eqref{Almog-Brenner}. 

\sspace
If $u_N - u_{eff} = o(\phi^2)$ in $H^{-\infty}_{loc}(\R^3)$, meaning that for all bounded open set  $U$, there exists $s \in \R$ such that 
 $$\limsup_{N \rightarrow +\infty} \|u_N - u_{eff} \|_{H^s(U)} =o(\phi^2) , \quad \text{as $\phi \rightarrow 0$}, $$
 then, necessarily, the coefficient $\nu_2$ defined in \eqref{def:nu2} satisfies $\: \nu_2 S : S =  \mu  \lim_{N \rightarrow +\infty}  \mV_N$
 where $\nu_2$ was defined in \eqref{def:nu2}, and 
 \begin{equation} \label{def:VN}
 \mV_N := \frac{75 |\mO|}{16 \pi } \left( \frac{1}{N^2}  \sum_{i \neq j} g_S(x_i - x_j)  -  \int_{\R^3 \times \R^3} g_S(x-y) f(x) f(y) dx dy \right) 
 \end{equation}
 with the Calder\'on-Zygmund kernel 
 \begin{equation} \label{def:g_S}
 g_S := -D \left( \frac{S : (x \otimes x)  x}{|x|^5}\right) : S. 
 \end{equation} 
\end{theorem}

\mspace
Roughly, this theorem states that {\em if there is an effective model at order $\phi^2$}, the mean quadratic correction $\nu_2$  is given by the limit of $\mV_N$, defined in \eqref{def:VN}. Note that the integral at the right-hand side of \eqref{def:VN} is well-defined:  $f \in L^2(\R)$ and  $f \rightarrow g_S \star f$ is a Calder\'on-Zygmund operator, therefore continuous on $L^2(\R^3)$. We insist that our result is an {\em if theorem}:  the limit of \eqref{def:VN} does not necessarily exist for any configuration of particles  $x_i = x_{i,N}$ satisfying \eqref{H1}-\eqref{H2}. In particular, it is not clear that an effective model at order $\phi^2$ is available for all such configurations. 

\mspace
Still, the second part of our analysis  shows that for points associated to  stationary random processes (including periodic patterns or Poisson hard core processes), the limit of the functional does exist, and is given by an  explicit formula.  We shall leave for later investigation the problem of approximating $u_N$ by $u_{eff}$ when the limit of $\mathcal{V}_N$ exists. 

\mspace
Our study of functional  \eqref{def:VN} is detailed in Sections \ref{sec2} to \ref{sec4}. It borrows a lot from the mathematical analysis of Coulomb gases, as developped over the last years by Sylvia Serfaty and her coauthors  \cite{MR3353821,MR3455593,MR3046995}. Although our paper is self-contained, we find useful to give a brief account of this analysis here.  As explained in the lecture notes \cite{MR3309890}, one of its main goals  is to understand what configurations of points minimize 
Coulomb energies of the form 
$$H_N = \frac{1}{N^2} \sum_{i\neq j} g(x_i - x_j) + \frac{1}{N} \sum_{i=1}^N V(x_i) $$
where $g(x) = \frac{1}{|x|}$ is a repulsive potential of Coulomb type, and  $V$ is typically a confining potential. It is well-known, see \cite[chapter 2]{MR3309890}, that under suitable assumptions on $V$, the sequence of functionals $H_N$  (seen as a functionals over probability measures by extension  by $+\infty$ outside the set of empirical measures) $\Gamma$-converges to the functional 
$$ H(\lambda) = \int_{\R^3 \times \R^3} g(x-y) d\lambda(x) d\mu(y)  + \int_{\R^3} V(x) d\lambda(x)  $$
Hence,  the empirical measure $\delta_N = \frac{1}{N} \sum_{i=1}^N \delta_{x_i}$ associated to the minimizer $(x_1,\dots,x_N)$ of $H_N$ converges weakly to the minimizer 
$\lambda$ of $H$. 

\mspace
In the series of works  \cite{MR3353821,MR3455593}, see also \cite{MR2945619} on the Ginzburg-Landau model, Serfaty and her coauthors investigate the next order term in the asymptotic expansion of $\min_{x_1,\dots,x_N} H_N$. A keypoint in these works is understanding the behaviour of  (the minimum of)
\begin{equation} \label{calHN}
 \mathcal{H}_N =  \int_{\R^3 \times \R^3\setminus\text{Diag}}  g(x-y) d(\delta_N - \lambda)(x) d(\delta_N -\lambda)(y) 
 \end{equation}
as $N \rightarrow +\infty$. This is done through the notion of renormalized energy. Roughly, the starting point behind this notion is the (abusive) formal identity  
\begin{equation} \label{formal_formula1}
 \text{''} \int_{\R^3 \times \R^3}  g(x-y) d(\delta_N - \lambda)(x) d(\delta_N -\lambda)(y) = \frac{1}{4\pi} \int_{\R^3} |\na h_N|^2 \: \text{''}
 \end{equation} 
where $h_N$ is the solution of $\Delta h_N = 4 \pi (\delta_N - \lambda)$ in $\R^3$. Of course, this identity does not make sense, as both sides are infinite. On one hand, the left-hand side is not well-defined: the potential $g$ is singular at the diagonal, so that the integral with respect to the product of the empirical measures diverges. On the other hand, the right-hand side is not better defined: as the empirical measure does not belong to $H^{-1}(\R^3)$,  $h_N$ is not in $\dot{H}^1(\R^3)$. 

\mspace
Still, as explained in  \cite[chapter 3]{MR3309890}, one can modify this identity, and show a formula  of the form 
\begin{equation} \label{renormalize}
 \mathcal{H}_N = \lim_{\eta \rightarrow 0} \left(  \frac{1}{4\pi}\int_{\R^3} |\na h^\eta_N|^2 - N g(\eta)  \right)
 \end{equation}
where $h_N^\eta$ is an approximation of $h_N$ obtained by regularization of the Dirac masses at the right-hand side of the Laplace equation: $\Delta h_N^\eta = 4 \pi (\delta^\eta_N - \lambda)$ in $\R^3$. Note the removal of  the term $N g(\eta)$ at the right-hand side of \eqref{renormalize}. This term, which goes to infinity as the parameter $\eta \rightarrow 0$, corresponds to the self-interaction of the Dirac masses:  it must be removed, consistently with the fact that the integral defining $\mathcal{H}_N$ excludes the diagonal. This explains the term {renormalized energy}. See \cite[chapter 3]{MR3309890} for more details. 

 \mspace
From there (omitting to discuss the delicate commutation of the limits in $N$ and $\eta$ !),  the asymptotics of $\min_{x_1, \dots, x_N} \mathcal{H}_N$ can be deduced from the one  of $\min_{x_1, \dots, x_N} \int_{\R^3} |\na h_N^\eta|^2$, for fixed $\eta$.  The next step is to show that such  minimum  can be expressed as spatial averages of (minimal) microscopic energies, expressed in terms of solutions of the so-called {\em jellium problems}: see \cite[chapter 4]{MR3309890}. These problems, obtained through rescaling and blow-up of the equation on $h_N^\eta$, are an analogue of cell problems in homogenization. More will be said in Section \ref{sec3}, and we refer to the lecture notes \cite{MR3309890} for all necessary complements. 

\mspace
Thus, the main idea in the second part of our paper is to take  advantage of the analogy between the functionals $\mV_N$ and $\mathcal{H}_N$  to apply the strategy just described. Doing so, we face specific difficulties: our distribution of points is not minimizing an energy, the potential $g_S$ is much more singular than $g$, the reformulation of the functional in terms of an energy  is less obvious, {\it etc}. Still, we are able to reproduce the same kind of scheme. We introduce in Section \ref{sec2} an analogue of the renormalized energy. The analogue of the jellium problem is discussed in Section \ref{sec3}. Finally, in Section \ref{sec4}, we are able to tackle the convergence of $\mV_N$, and give explicit formula for the limit in two cases:  the case of a (properly rescaled) $L \Z^3$-periodic pattern of $M$-spherical particles with centers $a_1$, \dots, $a_M$, and the case of a (properly rescaled) hardcore stationary random process with locally integrable two points correlation function $\rho_2(y,z) = \rho(y-z)$. In the first case, we show that  
\begin{equation} \label{eq:V_N:periodic}
\lim_{N \rightarrow +\infty} \mV_N = \frac{25 L^3}{2 M^2 }  \Bigl( \sum_{i \neq j}  S \na \cdot G_{S,L}(a_i - a_j) \: + \: K  S\na \cdot (G_{S,L} - G_S)(0) \Bigr),  
\end{equation}
where $G_S$ and $G_{S,L}$ are the whole space and  $L\Z^3$-periodic  kernels defined respectively in \eqref{def:G_S_p_S} and \eqref{eq:per_Green}. See Proposition 
\ref{prop: periodic}.  In the special case of a primitive cubic lattice, for which $M=L=1$, we can push further the calculation:  we find that 
$$ \nu_2 S : S = \mu \big( \alpha \sum_{i=1}^3 |S_{ii}|^2 + \beta \sum_{i \neq j} |S_{ij}|^2 \big), $$
with $\alpha \approx 9.48$ and $\beta \approx -2,15$, {\it cf. }Proposition \ref{prop:array} for precise expressions. Our result is  in agreement with \cite{ZuAdBr}. In the random stationary case, if the process has mean intensity one, we show that 
\begin{equation} \label{eq:V_N:random}
\begin{aligned} 
\lim_N \mV_N & = \frac{25}{2}  \lim_{L \rightarrow +\infty}   \frac{1}{L^3} \,  \sum_{\substack{z\neq z' \in \Lambda \cap K_L}} \!\!\! S \na \cdot G_{S,L}(z - z')  \\
&= \frac{25}{2} \lim_{L \rightarrow +\infty} \frac{1}{L^3} \int_{K_L \times K_L} S\na \cdot G_{S,L}(z-z') \rho(z-z') dz dz'. 
\end{aligned}
\end{equation}
These formula open the road to numerical computations of the viscosity coefficients of specific processes, and should in particular allow to check the formula found in the literature \cite{BaGr1, Obr}.  

\mspace
Let us conclude this introduction by pointing out that our analysis falls into the general scope of deriving macroscopic properties of dilute suspensions. In this perspective, it can be related to mathematical studies on the drag or sedimentation speed of suspensions, see \cite{MR2094523,MR2398959,MR3795188,MR3744384,Mec} among many.  See also the recent work \cite{MR3458165} on the conductivity problem.

\section{Expansion of the effective viscosity} \label{sec1}
The aim of this section is to understand the origin of the functional $\mV_N$  introduced in \eqref{def:VN}, and to prove Theorem \ref{thm1}. The outline is the following. We first consider the effective model \eqref{Stoeff}-\eqref{jumpBC}-\eqref{asymptmu}. Given $S$ a symmetric trace-free matrix, and a solution $u_{eff}$ with condition at infinity \eqref{BC:infinity}, we exhibit an integral quantity $\mathcal{I}_{eff} = \mathcal{I}_{eff}(S)$ that involves $u_{eff}$ and allows to recover (partially) the mean viscosity coefficient 
$\nu_2$.  In the next paragraph, we introduce the analogue $\mathcal{I}_{N}$ of $\mathcal{I}_{eff}$, that involves this time the solution $u_N$ of \eqref{Sto}-\eqref{Sto2} and  \eqref{BC:infinity}. In brief, we show that if $u_N$ is $o(\phi^2)$ close to $u_{eff}$, then $\mathcal{I}_{N}$ is $o(\phi^2)$ close to  $\mathcal{I}_{eff}$. Finally, we provide an expansion of $\mathcal{I}_{N}$, allowing to express $\nu_2$ in terms of $\mV_N$. Theorem \ref{thm1}  follows.  

\subsection{Recovering the viscosity coefficients in the effective model}
Let $k \ge 2$, $\mu'$ satisfying \eqref{asymptmu},  with viscosity coefficients $\mu_i$ that may depend on $x$.  Let $S$ symmetric and trace-free. We denote $u_0(x) = Sx$. Let $u_{eff} = 1_{\R^3\setminus\mO} u +  1_{\mO} u'$  the weak solution in $u_0+ \dot{H}^1(\R^3)$ of \eqref{Stoeff}-\eqref{jumpBC}-\eqref{BC:infinity}. 
By a standard energy estimate, one can show the expansion 
$$ u_{eff} - u_0 =  \phi \, u_{eff,1} +  \dots + \phi^k u_{eff,k}  + o(\phi^k) \quad \text{in } \dot{H}^1(\R^3) $$
where the system satisfied by $u_{eff,i} = 1_{\R^3\setminus\mO} u_i  +  1_{\mO} u'_i$  is derived by plugging the expansion in  \eqref{Stoeff}-\eqref{jumpBC} and keeping terms with power $\phi^i$ only.  We notably find 
\begin{equation}
\label{Sto_1}
\left\{
\begin{aligned}
-\mu \Delta u_1 + \na  p_1 & = 0, \quad \div u_1 = 0,  \quad x \in \R^3\setminus\mO, \\ 
-\mu \Delta u'_1 + \na  p'_1 & = 2 \div(\mu_1 D(u_0))   \quad \div u'_1 = 0, \quad x \in \mO, 
\end{aligned}
\right.
\end{equation}
together with the  conditions: $u_1 = 0$ at infinity,
$$ u_1 = u_1' \quad \text{at } \: \pa \mO, \quad \sigma_\mu(u_1,p_1)n   \: = \:  \sigma_\mu(u'_1, p'_1)n \: + \: 2 \mu_1 D(u_0)n \: \text{at } \: \pa \mO.$$  
Similarly, 
\begin{equation}
\label{Sto_2}
\left\{
\begin{aligned}
-\mu \Delta u_2 + \na  p_2 & = 0, \quad  \div u_2 = 0,  \quad x \in  \R^3\setminus\mO, \\ 
-\mu \Delta u'_2 + \na  p'_2 & = 2 \div(\mu_2 D(u_0)) + 2 \div( \mu_1 D(u'_1)),  \quad  \div u'_2 = 0, \quad x \in \mO, 
\end{aligned}
\right.
\end{equation}
together with:  $u_2 = 0$ at infinity,
$$ u_2 = u_2' \quad \text{at } \: \pa \mO, \quad \sigma_\mu(u_2,p_2)n  \: = \:  \sigma_\mu(u'_2, p'_2)n \: + \: 2 \mu_2 D(u_0)n + 2 \mu_1 D(u_1')n \: \text{ at } \: \pa \mO. $$  

\mspace
We now define, inspired by formula (4.11.16) in \cite{Batch_book},  
\begin{equation} \label{Ieff1}
\mathcal{I}_{eff} := \int_{\pa \mO} \sigma_\mu(u-u_0, p_{eff}) n \cdot S x ds  - 2 \mu \int_{\pa \mO} (u-u_0) \cdot Sn ds 
\end{equation}
where $n$ refers to the outward normal. 
We will show that 
\begin{equation} \label{Ieff2}
\mathcal{I}_{eff} = \: 2  |\mO|  \, \left( \phi \nu_1 S :  S  +   \, \phi^2 \nu_2 S : S \right) +  2 \phi^2 \int_{\mO}  \mu_1  D(u_1') : S + o(\phi^2). 
\end{equation}
We first use \eqref{jumpBC} to write 
\begin{align*}
\mathcal{I}_{eff} & = \int_{\pa \mO} \sigma_{\mu'}(u'-u_0, p') n \cdot S x ds + \int_{\pa \mO} \sigma_{\mu'-\mu}(u_0,0) n \cdot S x ds  - 2 \mu \int_{\pa \mO} (u' - u_0) \cdot Sn ds \\ 
& =  \int_{\pa \mO} \sigma_{\mu'}(\phi u'_1 + \phi^2 u'_2, \phi p_1 + \phi^2 p_2) n \cdot S x ds + 2 \int_{\pa \mO} (\phi \mu_1 + \phi^2 \mu_2) S n \cdot S x ds  \\ 
& - 2 \mu \int_{\pa \mO} (\phi u'_1 + \phi^2 u'_2) \cdot Sn ds + o(\phi^2) \\
& = \int_{\pa \mO} \sigma_{\mu}(\phi u'_1 + \phi^2 u'_2, \phi p_1 + \phi^2 p_2) n \cdot S x ds + \phi \int_{\pa \mO} \sigma_{\mu_1}(\phi u'_1, 0) n \cdot S x ds \\
& + 2 \int_{\pa \mO} (\phi \mu_1 + \phi^2 \mu_2) S n \cdot S x ds - 2 \mu \int_{\pa \mO} (\phi u'_1 + \phi^2 u'_2) \cdot Sn ds + o(\phi^2)
\end{align*}
Using the equations satisfied by $u'_1$ and $u'_2$, after integration by parts, we get 
\begin{align*}
& \int_{\pa \mO} \sigma_{\mu}(\phi u'_1 + \phi^2 u'_2, \phi p_1 + \phi^2 p_2) n \cdot S x ds \\ 
= & - \int_{\mO} 2 \div(\phi\mu_1 S + \phi^2 \mu_2 S) \cdot Sx dx - \int_{\mO} 2 \div(\phi^2 \mu_1 D(u'_1)) \cdot S x dx  + 2 \mu \int_{\mO} D(\phi u'_1 + \phi^2 u'_2) : S dx \\
= & 2 |\mO| (\phi \nu_1 S : S + \phi^2 \nu_2 S : S) - 2 \int_{\pa \mO} (\phi\mu_1 + \phi^2 \mu_2) Sn \cdot Sx ds \\ 
+ & 2 \int_{\mO} \phi^2 \mu_1 D(u'_1) : S - 2 \int_{\pa \mO} \phi^2 \mu_1 D(u'_1) n \cdot S x ds  + 2 \mu \int_{\mO} (\phi u'_1 + \phi^2 u'_2) \cdot Sn dx. 
\end{align*}  
 Plugging this last line in the expression for $\mathcal{I}_{eff}$ yields \eqref{Ieff2}. 
 
 \mspace
 We see through formula \eqref{Ieff2} that the expansion of $\mathcal{I}_{eff}$ in powers of $\phi$ gives access to $\nu_1$, and, if $\mu_1$ is known, it further gives access to $\nu_2$. On the basis of the works \cite{AlBr,NiSc} and of the recent paper \cite{HiWu} which considers the same setting as ours,  it is natural to assume  that $\mu_1$ is given by \eqref{Almog-Brenner}. This implies $\nu_1 = \frac{5}{2} \mu$. With such expression of $\mu_1$, and the form of $f$ specified in \eqref{H1}, we can check that 
 $u_S = (5 |\mO|)^{-1} u_{eff,1}$  satisfies 
 \begin{equation} \label{def:uA_mathcalO}
 -\Delta u_S + \na p  = \div(S f) = S\na f, \quad \div u_S = 0 \quad \text{in } \: \R^3, \quad \lim_{|x| \to \infty} u_S(x) = 0. 
 \end{equation}
It follows that 
\begin{equation} \label{Ieff3} 
 \mathcal{I}_{eff} =   5 \phi \mu |\mO| |S|^2 \: +\:  2\phi^2 |\mO|  \nu_2 S : S \: -  \:   50 \mu \phi^2 |\mO|^2 \int_{\R^3} |D(u_S)|^2  + o(\phi^2).
 \end{equation}
 
 \subsection{Recovering the viscosity coefficients in the model with particles}
To determine the possible value of the mean viscosity coefficient $\nu_2$, we must now relate the functional $\mathcal{I}_{eff}$, based on the effective model, to a functional $\mathcal{I}_N$ based on the real model with spherical rigid particles. From now on, we place ourselves under the assumptions of Theorem \ref{thm1}. Note that, thanks to hypothesis \eqref{H2}, the spherical particles do not overlap for $\phi$ small enough, so that a weak solution $u_N \in u_0 + \dot{H}^1(\R^3)$ of \eqref{Sto}-\eqref{Sto2}-\eqref{BC:infinity} exists and is unique.   

\mspace
By integration by parts, for any $R$ such that $\mO \Subset B_R$, we have 
\begin{equation} \label{Ieff4} 
 \mathcal{I}_{eff} = \int_{\pa B_R} \sigma_\mu(u_{eff}-u_0, p_{eff}) n \cdot S x ds  - 2 \mu \int_{\pa B_R} (u_{eff}-u_0) \cdot Sn ds 
 \end{equation}
By analogy with \eqref{Ieff1}, and as all particles remain in a fixed compact $K \supset \mO$ independent of $N$,  we set for any $R$ such that $K \subset B_R$:
\begin{equation} \label{IN} 
 \mathcal{I}_N := \int_{\pa B_R} \sigma_\mu(u_N-u_0, p_N) n \cdot S x ds  - 2 \mu \int_{\pa B_R} (u_N-u_0) \cdot Sn ds 
 \end{equation}
which again does not depend on our choice of $R$ by integration by parts. Now, if $u_{eff}$ and $u_N$ are $o(\phi^2)$-close in the sense of Theorem \ref{thm1}, then 
\begin{equation} \label{IN-Ieff}
\limsup_{N \rightarrow +\infty} |\mathcal{I}_N - \mathcal{I}_{eff}| = o(\phi^2).
\end{equation}
Indeed, $u_N - u_{eff}$ is a solution of a homogenenous Stokes equation outside $K$. By elliptic regularity, we find that $\limsup_{N \rightarrow +\infty} \|u_{eff} - u_N\|_{H^s(K')} = 0$, for any compact $K' \subset \R^3\setminus K$ and any positive $s$. Relation \eqref{IN-Ieff} follows. 

\mspace
We now turn to the most difficult part of this section, that is expanding $\mathcal{I}_N$ in powers of $\phi$. We aim at proving 
\begin{proposition} \label{prop:expansion_IN}
Let $(x_{i})_{1\le i \le N}$,  satisfying \eqref{H1}-\eqref{H2}. For $S$ trace-free and symmetric, for $\phi > 0$, let $u_N$  the solution of \eqref{Sto}-\eqref{Sto2}-\eqref{BC:infinity} with the ball radius $a$  defined through \eqref{def_phi}. Let $\mathcal{I}_N$ as in \eqref{IN}, $\mV_N$ as in \eqref{def:VN}, and $u_S$ the solution of \eqref{def:uA_mathcalO}. One has

\begin{equation}
\mathcal{I}_N =   5 \phi \mu |\mO| |S|^2 \: +\:  2  \phi^2 \mu |\mO| \mV_N -  \:   50 \mu \phi^2 |\mO|^2 \int_{\R^3} |D(u_S)|^2 + o(\phi^2)
 \end{equation}
 
\end{proposition}
\noindent
As before, notation  $A_N = B_N + o(\phi^2)$ means  $\limsup_N |A_N - B_N| = o(\phi^2)$.
Obviously, Theorem \ref{thm1} follows directly from \eqref{Ieff3}, \eqref{IN-Ieff} and from the proposition. 

\mspace
To start the proof, we set $v_N := u_N - u_0$. Note that $v_N \in \dot{H}^1(\mathcal F)$  still satisfies the Stokes equation outside the ball, with $v_N = 0$ at infinity, and 
$v_N = - Sx + u_i + \omega_i \times (x-x_i)$ inside $B_i$.  Moreover, taking into account the identities: 
$$ \int_{\pa B_i} \sigma_{\mu}(u_0,0)n \, ds = 2\mu \int_{\pa B_i} S n = 2 \mu  \int_{B_i} \div S = 0 $$
and 
\begin{equation} \label{calc:momentum}
\begin{aligned}
   \int_{\pa B_i} \sigma_{\mu}(u_0,0)n \times (x-x_i) \, ds & = 2\mu  \int_{\pa B_i} S n \times (x-x_i) \, ds  = 2\mu  \int_{\pa B_i} S (x-x_i) \times n \, ds \\
   & =  2\mu  \int_{B_i} \curl \left( S (x-x_i)\right)\, ds = 0, 
   \end{aligned}
   \end{equation}
one has for all $i$:
$$\int_{\pa B_i} \sigma_\mu(v_N,p_N) n \, ds  = 0, \quad \int_{\pa B_i} \sigma_\mu(v_N,p_N) n \times (x-x_i)  \, ds = 0.$$ 
From the definition \eqref{IN}, we can re-express $\mathcal{I}_N$ as 
\begin{equation}
\mathcal{I}_N =      \sum_{i=1}^N \int_{\pa B_i}  \sigma_\mu(v_N,p_N) n \cdot Sx \, ds  -  2\mu \sum_{i=1}^N \int_{\pa B_i}  \, v_N \cdot Sn  \, ds 
\end{equation}
To obtain an expansion of $\mathcal I_N$ in powers of $\phi$, we will now approximate $(v_N,p_N)$ by some explicit field $(v_{app}, p_{app})$, inspired by the method of reflections. This approximation involves the elementary problem: 
\begin{equation} \label{eq.stresslet}
\left\{
\begin{aligned}
-\mu \Delta v + \na p & = 0 \quad \text{ outside } B(0,a), \\ 
\div v & = 0 \quad \text{ outside } B(0,a), \\
v(x) & = - Sx, \quad x \in B(0,a). 
\end{aligned}   
\right.
\end{equation}
 The solution of \eqref{eq.stresslet} is explicit \cite{GuMo},  and  given by 
 \begin{equation} \label{def.stresslet.us}
 v^s[S] := -\frac{5}{2} S : (x \otimes x) \frac{a^3 x}{|x|^5} - Sx \frac{a^5}{|x|^5}  +\frac{5}{2} (S : x \otimes x) \frac{a^5 x}{|x|^7} = v[S] + O(a^5 |x|^{-4}) 
 \end{equation}
 with 
\begin{equation} \label{def.stresslet.u}
 v[S] :=-\frac{5}{2} S : (x \otimes x) \frac{a^3 x}{|x|^5}.
 \end{equation} 
 The pressure is
 \begin{equation} \label{def.stresslet.p}
 p^s[S]  :=  -5 \mu a^3  \frac{S : (x \otimes x)}{|x|^5}.
  \end{equation}
 We now introduce 
 \begin{equation} \label{def:vapp}
(v_{app},p_{app})(x) := \sum_{i=1}^N (v^s[S],p^s[S])(x-x_i) + \sum_{i=1}^N (v^s[S_i], p^s[S_i])(x-x_i),
\end{equation}  
where 
\begin{equation} \label{def:Sj}
S_i \: := \: \sum_{j \neq i} D(v[S])(x_i-x_j).  
\end{equation} 
In short, the first sum at the right-hand side of \eqref{def:vapp} corresponds to a superposition of $N$ elementary solutions, meaning that the interaction between the balls is neglected.  This sum satisfies the Stokes equation outside the ball, but creates an error at each ball $B_i$, whose leading term is $S_i x$.  This explains the correction by the second sum at the right-hand side of \eqref{def:vapp}. One could of course reiterate the process: as the distance between particles is large compared to their radius, we expect the interactions to be smaller and smaller. This is the principle of the method of reflections that is investigated in \cite{MR3744384}. 
From there, Proposition \ref{prop:expansion_IN} will follow from two facts. Defining   
$$ \mathcal{I}_{app} :=  \sum_{i=1}^N \int_{\pa B_i}  \sigma_\mu(v_{app},p_{app}) n \cdot Sx \, ds  -  2\mu \sum_{i=1}^N \int_{\pa B_i}  \, v_{app} \cdot Sn  \, ds $$
we will show first that 

\begin{equation} \label{Iapp1}
\mathcal{I}_{app} = 5 \phi \mu |S|^2 + 2 \phi^2 \mu  |\mO| \mV_N -  \:  50 \mu \phi^2 |\mO|^2 \int_{\R^3} |D(u_S)|^2
\end{equation}
and then
\begin{equation} \label{Iapp2}
\limsup_{N \rightarrow +\infty} |\mathcal{I}_N - \mathcal{I}_{app}| = o(\phi^2)
\end{equation}
Identity \eqref{Iapp1} follows from a calculation that we now detail.   We define 
  $$ \mI_i(v,p) :=   \int_{\pa B_i}  \bigl( (\sigma(v,p) n \otimes x)  - 2\mu (v \otimes n)  \bigr) \, ds $$
  We have
 \begin{align*}
  \mI_{app} 
 & = \sum_{i} \mI_i(v^s[S](\cdot-x_i), p^s[S](\cdot-x_i)) : S  + \sum_{i} \sum_{j\neq i} \mI_i(v^s[S](\cdot-x_j), p^s[S](\cdot-x_j)) : S  \\
 & + \sum_{i} \mI_i(v^s[S_i](\cdot-x_i), p^s[S_i](\cdot-x_i)) : S   + \sum_{i} \sum_{j\neq i} \mI_i(v^s[S_j](\cdot-x_j), p^s[S_j](\cdot-x_j)) : S \\
 &  =: I_a + I_b + I_c + I_d. 
 \end{align*}
 To treat $I_b$ and $I_d$, we rely on the following property, that is checked easily through integration by parts: {\em for any $(v,p)$ solution of Stokes in $B_i$, and any trace-free symmetric matrix $S$, $\mI_i(v,p) : S = 0$}. As for all $i$ and all $j\neq i$,  $v^s[S](\cdot-x_j)$ or $v^s[S_j](\cdot -x_j)$  is a solution  of Stokes inside $B_i$,  we deduce 
\begin{equation} \label{Ib_Id}
I_b = I_d = 0.
\end{equation}

\mspace
As regards $I_a$, we use the following formula which follows from a tedious calculation \cite{GuMo}: 
for any traceless matrix $S$, 
\begin{equation} \label{eq.Wi}
 \mI_i(v^s[S](\cdot-x_i)) =  \frac{20 \pi}{3} \mu a^3 S. 
\end{equation}
It follows that 
\begin{equation} \label{Ia}
I_a =  N \frac{20 \pi}{3} \mu a^3 |S|^2 =  5 \phi |\mO| \mu |S|^2  
\end{equation}
This term corresponds to the famous Einstein formula for the mean effective viscosity. It is coherent with the expression \eqref{Almog-Brenner} for $\mu_1$, which implies $\nu_1 = \frac{5}{2} \mu$.  

\mspace
Eventually, as regards $I_c$, we can use again \eqref{eq.Wi}, replacing $S$ by $S_i$:

\begin{align}
\nonumber I_c & =   \frac{20 \pi}{3} \mu a^3 \sum_{i}  S_i : S   = \frac{20 \pi}{3}  \mu a^3  \sum_{i}  \sum_{j \neq i} D(v[S])(x_i-x_j) : S \\
\nonumber  & = \frac{75 |\mO|^2}{8\pi} \mu \phi^2 \frac{1}{N^2}  \sum_{i}  \sum_{j \neq i} g_S(x_i-x_j)  & \\ 
\label{Ic}  & = 2 \phi^2\mu |\mO| \mV_N + \phi^2 \frac{75 |\mO|^2}{8\pi} \mu \int_{\R^3 \times \R^3} g_S(x-y) f(x) f(y) dx dy, 
 \end{align}
with $g_S$ defined in \eqref{def:g_S}. In view of \eqref{Ib_Id}-\eqref{Ia}-\eqref{Ic},  to conclude that \eqref{Iapp1} holds, it is enough  to prove
\begin{lemma} \label{lem:sign:gA}
For 	any $f \in L^2(\R^3)$,
\begin{equation} \label{link:gS:energy}   
\int_{\R^3 \times \R^3}   g_S(x-y) f(x) f(y) dx dy  = - \frac{16\pi}{3} \int_{\R^3} |D(u_S)|^2,
\end{equation}
with  $g_S$ defined in \eqref{def:g_S}, and $u_S \in \dot{H}^1(\R^3)$ the solution of \eqref{def:uA_mathcalO}.
\end{lemma}

\noindent
{\em Proof.} Note that both sides of the identity are continuous over $L^2$: the left-hand side is continuous as the Calder\'on-Zygmund operator $f \rightarrow g_S \star f$ is continuous over $L^2$, while the right-hand side is continuous by classical elliptic estimates for the Stokes operator. By density, this is therefore enough to assume that $f \in C^\infty_c(\R^3)$. 
We denote by $U = (U_{ij}), Q = (Q_j)$ the fondamental solution of the Stokes operator. This means that for all $j$, the vector field $U_j = (U_{ij})_{1\le i \le N}$ and the scalar field $Q_j$ satisfy the Stokes equation
\begin{equation} \label{def:Uj}
 -\Delta U_j + \na Q_j = \delta e_j, \quad \div U_j = 0 \quad \text{ in } \R^3.  
\end{equation}
It is well-known, see \cite[page 239]{MR1284205}, that 
$$ U(x) =  \frac{1}{8\pi} \left( \frac{1}{|x|} Id + \frac{x \otimes x}{|x|^3} \right), \quad Q(x) = \frac{1}{4\pi} \frac{x}{|x|^3}. $$
From there, one can deduce the following formula, {\it cf}  \cite[page 290, equation  (IV.8.14)]{MR1284205}:  
$$ \sigma(U_j, Q_j) = - \frac{3}{4\pi} \frac{(x \otimes x) x_j}{|x|^5}.  $$
Using the  Einstein convention for summation, this implies in turn that 
\begin{align} 
\nonumber g_S(x) &  = -  S_{kl} \pa_{x_k} \left( \frac{S : (x \otimes x)  x_l}{|x|^5}\right) =   \frac{4\pi}{3}  S :  S_{kl} \pa_{x_k} \sigma(U_l,Q_l)(x)  \\
\label{formula_g_S} & =  \frac{8\pi}{3}   S :  D S_{kl} \pa_{x_k} U_l  = (S\na) \cdot  (S_{kl} \pa_{x_k} U_l)
\end{align}
where we have used that $S$ is  trace-free to obtain the third equality. Hence, 
\begin{align} 
\int \int g_S(x-y)  f(x) dx f(y) dy 
\nonumber & =  \frac{8\pi}{3}   \int_{\R^3}  \big((S :  D S_{kl} \pa_{x_k} U_l) \star f \big)(y) f(y) dy  \\
\label{eq:gA} & = \frac{8\pi}{3} \int  S : D   S_{kl} \pa_{x_k}   (U_l \star f)(y) \, f(y) dy.
\end{align}
Note that the permutations between the derivatives and the convolution product do not raise any difficulty, as $f \in C_c^\infty(\R^3)$.   
Now,  using  $S_{kl} = S_{lk}$, and denoting by ${\rm St}^{-1}$ the convolution with the fundamental solution (inverse of the Stokes operator), we get 
\begin{equation} \label{def:GA}
 S_{kl} \pa_{x_k}  \int   U_l(y-x)  f(x) dx  = {\rm St}^{-1} (S \na f)(y). 
\end{equation}
Eventually, 
\begin{align*}
 \int \int g_S(x-y)  f(x) f(y) dx  dy  &  = \frac{8\pi}{3} \int  S : \na {\rm St}^{-1} (S \na f)(y) \, f(y) dy \\
& =  - \frac{8\pi}{3} \int  {\rm St}^{-1}(S \na f)(y) \cdot (S \na f)(y) \, dy = - \frac{16\pi}{3} \int_{\R^3} |D(u_S)|^2.
\end{align*}
This concludes the proof of the lemma. 
\begin{remark} 
By polarization of the previous identity, at least  for $f, \tilde f$ smooth and decaying enough, one has
\begin{equation} \label{polarization}
\begin{aligned}
\int \int g_S(x-y) f(y) \tilde f(x) dx & =  - \frac{8\pi}{3} \int  {\rm St}^{-1}(S \na f)(x) \cdot (S \na \tilde f)(x) \, dx \\
& = \frac{8\pi}{3} \int  (S \na) \cdot \big( {\rm St}^{-1}(S \na f)\big)(x)  \, \tilde f(x) \, dx  
\end{aligned}
\end{equation}
 \end{remark}

\mspace
The last step in proving Proposition \ref{prop:expansion_IN}, hence Theorem \ref{thm1},  is to show the bound \eqref{Iapp2}. If $w := v_N - v_{app}$, $q := p_N - p_{app}$, 
$$ \mI_N - \mI_{app} = \sum_{i=1}^N \int_{\pa B_i}  \sigma_\mu(w,q) n \cdot Sx \, ds  -  2\mu \sum_{i=1}^N \int_{\pa B_i}  \, w \cdot Sn  \, ds $$
Direct verifications  show that $v_{app}$, hence $w$,  satisfies the same force- and torque-free conditions as $v$. This means that for any family of constant vectors $u_i$ and $\omega_i$, $1 \le i \le N$, 
$$ \mI_N - \mI_{app}  = \sum_{i=1}^N \int_{\pa B_i}  \sigma_\mu(w,q) n \cdot (Sx - u_i - \omega_i \times (x-x_i)) \, ds  -  2\mu \sum_{i=1}^N \int_{\pa B_i}  \, w \cdot Sn  \, ds $$
By a proper choice of $u_i$ and $\omega_i$, we find 
\begin{align} 
\nonumber \mI_N - \mI_{app}  & = -\sum_{i=1}^N \int_{\pa B_i}  \sigma_\mu(w,q) n \cdot v_N  \, ds  -  2\mu \sum_{i=1}^N \int_{\pa B_i}  \, w \cdot Sn  \, ds \\
\nonumber & = - \int_{\mathcal{F}}  2\mu D(w) :  D(v_N)  \, dx -  2\mu \sum_{i=1}^N \int_{B_i}  \, D(w) : S  \, dx \\
\nonumber & =  -\sum_{i=1}^N \int_{\pa B_i}  \sigma_\mu(v_N,p_N) n \cdot w  \, ds -  2\mu \sum_{i=1}^N \int_{B_i}  \, D(w) : S  \, dx \\
\label{rhs_IN-Iapp} & = -\sum_{i=1}^N \int_{\pa B_i}  \sigma_\mu(v_N,p_N) n \cdot (w +  \tilde u_i + \tilde \omega_i \times (x-x_i)) \, ds -  2\mu \sum_{i=1}^N \int_{B_i}  \, D(w) : S  \, dx
 \end{align}
for any family $(\tilde u_i, \tilde \omega_i)$, using this time that $v_N$ is force- and torque-free. Let $q \ge 2$. By a proper choice of $(\tilde u_i, \tilde \omega_i)$, by Poincar\'e and Korn inequalities, one  can ensure that for all $i$, 
$$ \|w +  \tilde u_i + \tilde \omega_i \times (x-x_i)\|_{W^{1-\frac{1}{q},q}(\pa B_i)} \le C \|D(w)\|_{L^q(B_i)} $$
where 
$$ \|g\|_{W^{1-\frac{1}{q},q}(\pa B_i)} = \inf \Big\{ \frac{1}{a} \|G\|_{L^q(B_i)} + \|\na G\|_{L^q(B_i)}, \quad G\vert_{\pa B_i} = g \Big\}$$
Note that the factor $\frac{1}{a}$ at the right-hand side is consistent with scaling considerations.  Moreover, by standard use of the Bogovskii operator, see \cite{MR1284205},  there exists a constant $C$ (depending only on the constant $c$ in \eqref{H2}) and a field $W \in W^{1,q}(\mathcal{F})$ , zero outside $\cup_{i=1}^N B(x_i,2a)$   satisfying 
\begin{align*}
& \div W = 0 \quad \text{in } \mathcal{F}, \quad W\vert_{B_i} = (w +  \tilde u_i + \tilde \omega_i \times (x-x_i))\vert_{B_i}, \\
& \| D(W) \|^q_{L^q(\mathcal{F})} \le \sum_i   \|w +  \tilde u_i + \tilde \omega_i \times (x-x_i)\|_{W^{1-\frac{1}{q},q}(B_i)}^q. 
\end{align*}
We deduce, with $p \le 2$ the conjugate exponent of $q$:  
\begin{align*}   
&\big|\sum_{i=1}^N \int_{\pa B_i}  \sigma_\mu(v_N,p_N) n \cdot (w +  \tilde u_i + \tilde \omega_i \times (x-x_i)) \, ds \bigr| = 2\mu \big| \int_{\mathcal{F}} D(v_N) : D(W) \bigr| \\ 
& \le 2 \mu \|D(v_N)\|_{L^p(\cup B(x_i,2a))} \|D(W)\|_{L^q(\mathcal{F})} \le C \phi^{1/p-1/2} \|D(v_N)\|_{L^2(\R^3)} \big( \sum_i \|D(w)\|^q_{L^q(B_i)} \big)^{1/q}.
\end{align*} 
By well-known variational properties of the Stokes solution,  $\|D(v_N)\|_{L^2}$ minimizes $\|D(v)\|_{L^2}$ over the set of all $v$ in $\dot{H}^1(\R^3)$ satisfying a boundary condition of the form $v\vert_{B_i} = - Sx + u_i + \omega_i \times (x-x_i)$ for all $i$. By the same considerations as before, based on the Bogovski operator, we infer that 
$$\|D(v_N)\|_{L^2(\R^3)}^2 \le C \sum_{i=1}^N \|D(- Sx)\|_{L^2(B_i)}^2 \le C' \phi  $$
so that 
$$ \big|\sum_{i=1}^N \int_{\pa B_i}  \sigma_\mu(v_N,p_N) n \cdot (w +  \tilde u_i + \tilde \omega_i \times (x-x_i)) \, ds \bigr|  \le C \phi^{1/p} \big( \sum_i \|D(w)\|^q_{L^q(B_i)} \big)^{1/q}. $$
Using this inequality with the first term in \eqref{rhs_IN-Iapp} and applying the H\"older inequality to the second term, we end up with 
\begin{equation} \label{IN-Iapp}
 |\mI_N - \mI_{app}|  \le C \phi^{1/p} \big( \sum_i \|D(w)\|^q_{L^q(B_i)} \big)^{1/q} 
\end{equation} 
To deduce \eqref{Iapp2}, it is now enough to prove that for all $q > 1$, there exists a constant $C$ independent of $N$ or $\phi$ such that  
\begin{equation} \label{bound_w_q}
\sum_i \|D(w)\|^q_{L^q(B_i)} \le  C ( \phi^{1+\frac{2q}{p}} + \phi^{1+\frac{4q}{3}} )
\end{equation}
Indeed, taking $q > 2$, meaning $p < 2$, and combining this inequality with \eqref{IN-Iapp} yields  \eqref{Iapp2}, more precisely 
$$  |\mI_N - \mI_{app}|  \le C  (\phi^{1+\frac{2}{p}} + \phi^{\frac{7}{3}}). $$
In order to show the bound \eqref{bound_w_q}, we must write down the expression for $w\vert_{B_i} = v_N\vert_{B_i} - v_{app}\vert_{B_i}$, where $v_{app}$ was introduced in \eqref{def:vapp}.  A little calculation, using Taylor's formula with integral remainder, shows that    
\begin{equation} \label{expression:w}
w\vert_{B_i}(x)  = w^r_i(x) - D_i (x - x_i)  -  E_i (x - x_i)  -  \mathbf{F}_i\vert_{x}(x-x_i,x-x_i)  
\end{equation}
with $w_i^r$ a rigid vector field (that disappears when taking the symmetric gradient), with  
$$ D_{i} := \sum_{j \neq i}  D(v[S_j])(x_i-x_j) , \quad E_{i} := \sum_{j \neq i} D(v^s[S+S_j]-v[S+S_j])(x_i-x_j) $$ 
and with the bilinear application: 
$$ \mathbf{F}_i\vert_{x} := \sum_{j \neq i}  \int_0^1 (1-t)  \na^2 v^s[S+S_j](t(x-x_i) + x_i -x_j) dt.$$ 
We remind that $v^s[S]$ and $v[S]$ were introduced in \eqref{def.stresslet.us} and \eqref{def.stresslet.u}, while the matrices $S_j$ are defined in \eqref{def:Sj}. Note that the matrices $D_i$ and $S_i$ have the same kind of structure. More precisely, we can define for a collection $(A_1, \dots, A_N)$ of $N$ symmetric matrices, an application 
$$ \mathcal{A} : (A_1, \dots, A_N) \rightarrow (A'_1, \dots, A'_N), \quad A'_i = \sum_{j\neq i} D(v[A_j])(x_i - x_j). $$   
Then, $(S_1, \dots, S_N) =  \mathcal{A}(S, \dots,S)$ and   $(D_1, \dots, D_N) =  \mathcal{A}(S_1, \dots,S_N) =  \mathcal{A}^2(S, \dots,S)$.  
Note that for any matrix $A$, the kernel $D(v[A])$, homogenenous of degree $-3$, is of Calder\'on-Zygmund type. Using this property,  we are able to prove in the appendix the following lemma, which is an adaptation of a result by the second author and Di Wu \cite{HiWu}: 
\begin{lemma} \label{Calderon}
For all $1 < q < +\infty$, there exists a constant $C$, depending on $q$ and on the constant $c$ in \eqref{H2}, such that, if $(A'_1, \dots, A'_N) = \mathcal{A}(A_1, \dots, A_N)$, then 
$$ \sum_{i=1}^N |A'_i|^q \le C \phi^{\frac{q}{p}} \sum_{i=1}^N |A_i|^q $$
\end{lemma}

\noindent
We can now proceed to the proof of \eqref{bound_w_q}. Denoting $w_i^1:= D_i (x-x_i)$, we find by the lemma:
$$  \sum_i \|D(w^1_i)\|^q_{L^q(B_i)} \le  C  a^3  \sum_i  |D_i|^q \le C' a^3 \phi^{\frac{q}{p}} \sum_{i=1}^N |S_i|^q  \le C'' a^3\phi^{\frac{2q}{p}} \sum_{i=1}^N |S|^q \le \mathcal{C} \phi^{1+\frac{2q}{p}}. $$
Then, we notice that for any matrix $A$,  $|D(v^s[A] - v[A])(x)| = O(a^5 |x|^{-5})$. This implies  that $w^2_i :=  E_i (x - x_i)$ satisfies 
\begin{align*} 
 &\sum_i \|D(w^2_i)\|^q_{L^q(B_i)} \le  C  a^3  \sum_i  |E_i|^q \le C'  a^3 a^{5q} \sum_i \Bigl( \sum_{j \neq i} \frac{|S_j| + |S|}{|x_i - x_j|^5} \Bigr)^q 
\end{align*}
By assumption \eqref{H2}, the points $y_i := N^{1/3} x_i$ satisfy  for all $i \neq j$:
$$|y_i - y_j| \ge \frac{1}{2} ( c +  |y_i - y_j|) \ge c.$$
 In particular,  
\begin{align*} 
 \sum_i \|D(w^2_i)\|^q_{L^q(B_i)} & \le  C a^3 \phi^{5q/3} \sum_i \bigl( \sum_{j} \frac{|S| + |S_j|}{(c + |y_i - y_j|)^5} \bigr)^q 
 \end{align*}
 We then make use of the following easy generalization of Young's convolution inequality: 
\begin{equation}  \label{general_convolution} 
 \forall q \ge 1, \quad  \sum_{i}   ( \sum_j |a_{ij} b_j| )^q \le  \max\big(\sup_i \sum_j |a_{ij}|,  \sup_j \sum_i |a_{ij}|\big)^q \sum_i |b_i|^q. 
 \end{equation}
 Applied with $a_{ij} = \frac{1}{(c + |y_i - y_j|)^5}$ and $b_j = |S| + |S_j|$, together with Lemma  \ref{Calderon},  it yields 
 \begin{align*} 
 \sum_i \|D(w^2_i)\|^q_{L^q(B_i)} & \le C a^3  \phi^{5q/3}  \bigl( \sum_j |S|^q + |S_j|^q \bigr)  \le  C' a^3  \phi^{5q/3} (1+\phi^{\frac{q}{p}}) N  \le \mathcal{C} \phi^{1+\frac{5q}{3}}.
\end{align*}

\mspace
It remains to bound the symmetric gradient of $w_i^3 :=  \mathbf{F}_i\vert_{x}(x-x_i,x-x_i)$. By the expression of $v^s$, we get that in $B_i$: 
$$ |D(w_i^3)| \le C \sum_{j \neq i}\left( \frac{a^5}{|x_i - x_j|^5} +  \frac{a^4}{|x_i - x_j|^4}\right) (|S| + |S_j|)  $$
Proceeding as above, we find 
$$  \sum_i \|D(w^3_i)\|^q_{L^q(B_i)}  \le  C a^3  ( \phi^{5q/3} + \phi^{4q/3})  (1+\phi^{\frac{q}{p}}) N  \le C' \phi^{1+\frac{4q}{3}} $$  
As $D(w) = D(w_i^1) + D(w_i^2) + D(w_i^3)$, {\it cf.} \eqref{expression:w},  the previous estimates yield \eqref{bound_w_q}. This concludes the proof of Proposition \ref{prop:expansion_IN}, and therefore the proof of Theorem \ref{thm1}.

\section{The $\phi^2$ correction $\mV_N$ as a renormalized energy} \label{sec2}
We start in this section the asymptotic analysis of the viscosity coefficient 
$$ \mathcal{V}_N =  \frac{75 |\mO|}{16 \pi }  \Big(  \frac{1}{N^2}  \sum_{i \neq j} g_S(x_i - x_j)  -  \int_{\R^3 \times \R^3} g_S(x-y) f(x) f(y) dx dy \Big) $$
As a preliminary step, we will show that there is no loss of generality in assuming  
\begin{equation}  \label{hyp:in_O}
 \forall i \in \{1, \dots, N\}, \quad \text{dist}(x_i, \mO^c) \ge \frac{1}{\ln N}. 
\end{equation} 
 We introduce the set 
 $$ I_{N,ext} = \big\{ 1 \le i \le N, \: \text{dist}(x_i, \mO^c)  \le  \frac{1}{\ln N}  \big\}, \quad \text{and }  N_{ext} = N_{ext}(N) := |I_{N,ext}|.  $$
 By \eqref{H1}-\eqref{H2}, it is easily seen that $N_{ext} = o(N)$ as $N \rightarrow +\infty$. We now show 
 \begin{lemma} \label{lem:in_O}
 $\mV_N$ is uniformly bounded in $N$, and 
 $$  \mathcal{V}_{N,ext} := \mV_N -  \frac{75 |\mO|}{16 \pi}  \Big(  \frac{1}{(N - N_{ext})^2} \sum_{\substack{i\neq j \\ i,j \notin I_{N,ext}}} g_S(x_i - x_j)  -  \int_{\R^3 \times \R^3} g_S(x-y) f(x) f(y) dx dy \Big) $$
goes to zero as $N\rightarrow +\infty$.
 \end{lemma}
 
 \mspace
{\em Proof.}
For any open set $U$, we denote $\dashint_U = \frac{1}{|U|} \int_U$. 

\mspace
Let $d := \frac{c}{4} N^{-1/3} \le  \min_{i \neq j}\frac{|x_i-x_j|}{4}$ by \eqref{H2}.  We write 
 \begin{align*}
\frac{1}{N^2} \sum_{i \neq j} g_S(x_i - x_j)
= &   \frac{1}{N^2} \sum_{i \neq j} \left( g_S(x_i - x_j) - \dashint_{B(x_j,d)} g_S(x_i-y) dy \right) \\
+ &  \frac{1}{N^2} \sum_{i \neq j} \left( \dashint_{B(x_j,d)} g_S(x_i-y) dy - \dashint_{B(x_i,d)} \dashint_{B(x_j,d)}  g_S(x-y) dx dy \right) \\
+ &   \frac{1}{N^2} \sum_{i \neq j}  \dashint_{B(x_i,d)} \dashint_{B(x_j,d)}  g_S(x-y) dx dy \: := I + II + III. 
\end{align*}
 For the first term, with $y_i := N^{1/3} x_i$ and with  \eqref{H2} in mind, that is $|y_i - y_j| \ge c$ for $i \neq j$:
\begin{align*}
\bigl| g_S(x_i - x_j) - \dashint_{B(x_j,d)} g_S(x_i-y) dy \bigr| &  \le  \dashint_{B(x_j,d)} \sup_{z \in [x_j, y]}|\na g_S|(x_i-z)| |x_j - y| dy \\
 & \le C N^{4/3}\frac{d}{(c + |y_i - y_j|)^4}, 
 \end{align*}
 see \eqref{def:g_S}. This yields, by a discrete convolution inequality: 
\begin{align*}
|I|  & \le \frac{CN^{7/3}}{N^2} d  \sup_i \sum_{j}   \frac{1}{(c+|y_i- y_j|)^4} \le  C' N^{1/3} d \le  \mathcal{C} 
\end{align*}
where we have used that $\sum_{j=1}^N  \frac{1}{(c+|y_i - y_j|)^4}$ is uniformly bounded in $N$ and in the index $i$ thanks to the separation assumption.  By similar arguments, $ | II | \le \mathcal{C}$.
As regards the last term, we notice that 
$$ |III| \le \frac{1}{N^2 d^6} \bigl| \int_{\R^3 \times \R^3} g_S(x-y) F_N(x) F_N(y) dy - \sum_{i=1}^N  \int_{\R^3 \times \R^3} g_S(x-y)  1_{B(x_i,d)}(x) 1_{B(x_i,d)}(y) dx dy \bigr|   $$
where $F_N = \sum_{i=1}^N 1_{B(x_i,d)}$. The operator $\mathcal{T} F(x) = \int g_S(x-y) F(y) dy$ is a  Calder\'on-Zygmund operator, and therefore continuous over $L^2$. As $F_N^2 = F_N$ (the balls are disjoint), we find that the $L^2$ norm of $F_N$ is $(N d^3)^{1/2}$ and 
$$    \bigl|   \int_{\R^3 \times \R^3} g_S(x-y) F_N(x) F_N(y) dy \bigr|  \le \|\mathcal{T}\| \|F_N \|_{L^2}^2 \le  \|\mathcal{T}\| N d^3.$$
Similarly, 
$$  \sum_{i=1}^N  \bigl| \int_{\R^3 \times \R^3} g_S(x-y)  1_{B(x_i,\eta)}(x) 1_{B(x_i,\eta)}(y) dx dy \bigr| \le N \|\mathcal{T}\| d^3. $$
It follows that $|III| \le \frac{C}{N d^3}$. With our choice of $d$, the first part of the lemma is proved. 

\mspace
From there, to prove that $\mV_{N,ext}$ goes to zero, as $N_{ext} = o(N)$, it is enough to show that 
$$ \frac{1}{N^2} \big(  \sum_{i \neq j} g_S(x_i - x_j) -  \hspace{-0.3cm} \sum_{\substack{i\neq j, \\ i,j \notin I_{N,ext}}} \hspace{-0.3cm} g_S(x_i - x_j) \big)  \: \rightarrow \: 0. $$
By symmetry, it is enough that 
$$\frac{1}{N^2}  \hspace{-0.3cm}  \sum_{\substack{i  \neq j, \\ i \in I_{N,ext}}}   \hspace{-0.3cm} g_S(x_i - x_j)  \: \rightarrow \: 0.$$ 
This can be shown by a similar decomposition as the previous one. Namely, 
\begin{align*}
\frac{1}{N^2} \sum_{i \neq j} g_S(x_i - x_j)
= &   \frac{1}{N^2}  \sum_{\substack{i   \neq j \\ i \in I_{N,ext}}} \left( g_S(x_i - x_j) - \dashint_{B(x_j,d)} g_S(x_i-y) dy \right) \\
+ &  \frac{1}{N^2} \sum_{\substack{i   \neq j \\ i \in I_{N,ext}}} \left( \dashint_{B(x_j,d)} g_S(x_i-y) dy - \dashint_{B(x_i,d)} \dashint_{B(x_j,d)}  g_S(x-y) dx dy \right) \\
+ &   \frac{1}{N^2} \sum_{\substack{i   \neq j \\ i \in I_{N,ext}}}  \dashint_{B(x_i,d)} \dashint_{B(x_j,d)}  g_S(x-y) dx dy \: := I_{ext} + II_{ext} + III_{ext}. 
\end{align*}
Proceeding as above, we find this time
$$| I_{ext} | + | II_{ext} |  + | III_{ext} | \le \mathcal{C} \frac{N_{ext}}{N} \rightarrow 0 \quad \text{as } N \rightarrow +\infty $$
which concludes the proof. 
\begin{remark} \label{rem:points_in_O}
By Lemma \ref{lem:in_O}, there is no restriction assuming \eqref{hyp:in_O} when studying the asymptotic behaviour of $\mV_N$. {\em Therefore, we make from now on the assumption \eqref{hyp:in_O}.}   
\end{remark}
\mspace
As explained in the introduction, the analysis of $\mV_N$ will rely on the mathematical methods introduced over the last years for Coulomb gases, the core problem being the analysis of a functional of the form \eqref{calHN}. We shall first reexpress $\mathcal{V}_N$ in a similar form. More precisely, we will show 
\begin{proposition} \label{prop:VN_approx_WN}
Denoting 
$$
\mW_N := \frac{75 |\mO|}{16\pi} \int_{\R^3 \times \R^3\setminus Diag} g_S(x-y)  \big(d\delta_N(x) - f(x) dx\big)   \big(d\delta_N(y) - f(y) dy\big) ,
$$
we have $\mV_N = \mW_N + \varepsilon(N)$ where $\varepsilon(N) \to 0$ as $N \to \infty.$ 
\end{proposition}
\begin{remark} \label{rem:meaning_WN}
In the definition of $\mW_N$, the integrals of the form 
$$  \int_{\R^3 \times \R^3\setminus Diag} g_S(x-y) d\delta_N(x) f(y) dy, \: \text{ and } \int_{\R^3 \times \R^3\setminus Diag} g_S(x-y) f(x) dx d\delta_N(y),   $$
that appear when expanding the product, are understood as 
\begin{align*}  
& \int_{\R^3 \times \R^3\setminus Diag} g_S(x-y) d\delta_N(x) f(y) dy  = \frac{8\pi}{3}\frac{1}{N} \sum_{i=1}^N S \na \cdot {\rm St}^{-1} S\na f(x_i), \\
& \int_{\R^3 \times \R^3\setminus Diag} g_S(x-y)  f(x) dx d\delta_N(y) = \frac{8\pi}{3}\frac{1}{N} \sum_{i=1}^N S \na \cdot {\rm St}^{-1} S\na f(x_i), 
\end{align*}
where ${\rm St}$ is the Stokes operator, see \eqref{polarization} and  the proof below for an explanation. 
\end{remark}

\mspace
{\em Proof.} Clearly, 
$$ \mV_N = \frac{75 |\mO|}{16\pi} \int_{\R^3 \times \R^3\setminus Diag}   g_S(x-y)  \Big(d\delta_N(x) d\delta_N(y) -  f(x) f(y) dx dy \Big) $$
so that formally 
\begin{align*}
 \mV_N = \mW_N  &  +  \frac{75 |\mO|}{16\pi}  \int_{\R^3 \times \R^3\setminus Diag} g_S(x-y)  (d\delta_N(x) - f(x)dx) f(y) dy \\
& +  \frac{75 |\mO|}{16\pi}   \int_{\R^3 \times \R^3\setminus Diag} g_S(x-y) f(x) dx (d\delta_N(y) -f(y) dy).
 \end{align*}
Note that it is not obvious that this formal decomposition makes sense, because all three quantities at the right-hand side involve integrals of $g_S(x-y)$ against product measures of the form $d \delta_N(x) f(y) dy$ (or the symmetric one), which may fail to converge because of the singularity of $g_S$. To solve this issue, a rigorous path consists in approximating, at fixed $N$, each Dirac mass $\delta_{x_i}$ by a (compactly suppported) approximation of unity $\rho_\eta(x-x_i)$, where $\eta > 0$ is the approximation parameter and goes to zero. One can then set, for each $\eta$, $\delta_N^\eta(x) :=  \frac{1}{N} \sum_{i=1}^N \rho_\eta(x-x_i)$, leading to the rigorous decomposition
\begin{align*}
 \mV_N^\eta  = \mW_N^\eta &  +  \frac{75 |\mO|}{16\pi}  \int_{\R^3 \times \R^3\setminus Diag} g_S(x-y)  (\delta_N^\eta(x) d(x) - f(x)dx) f(y) dy \\
& +  \frac{75 |\mO|}{16\pi}   \int_{\R^3 \times \R^3\setminus Diag} g_S(x-y) f(x) dx (\delta_N^\eta(y) dy -f(y) dy) 
 \end{align*}  
where  $\mV_N^\eta$, $\mW_N^\eta$ are deduced from $\mV_N$, $\mW_N$  replacing the empirical measure by its regularization. It is easy to show that $\lim_{\eta \rightarrow 0} \mV_N^\eta = \mV_N$. To conclude the proof, we shall establish the following: first,  
\begin{equation} \label{lim_eta_integrale_double}
\lim_{\eta \rightarrow 0} \int_{\R^3 \times \R^3\setminus Diag} g_S(x-y)  \delta_N^\eta(x) dx f(y) dy =  \frac{8\pi}{3} \frac{1}{N} \sum_{i=1}^N S \na {\rm St}^{-1} S\na f(x_i), 
\end{equation}
the same limit holding for the symmetric term. In particular,   \eqref{lim_eta_integrale_double} will show that $\mW_N = \lim_{\eta \rightarrow 0} \mW_N^\eta$ exists, in the sense given in Remark \ref{rem:meaning_WN}.  Then, we will prove 
\begin{equation}  \label{lim_N_integrale_double}
\lim_{N \rightarrow +\infty}  \frac{8\pi}{3} \frac{1}{N} \sum_{i=1}^N S \na {\rm St}^{-1} S\na f(x_i)  =    \int_{\R^3 \times \R^3\setminus Diag} g_S(x-y) f(x) f(y) dx dy. 
\end{equation}
which together with   \eqref{lim_eta_integrale_double} will complete the proof of the proposition. 

\mspace
The limit \eqref{lim_eta_integrale_double} follows from identity \eqref{polarization}. Indeed, for $\eta > 0$, this formula yields 
$$ \int_{\R^3 \times \R^3\setminus Diag} g_S(x-y)  \delta_N^\eta(x) dx f(y) dy = - \frac{8\pi}{3} \int_{\R^3}  {\rm St}^{-1}(S \na f)(x)  \cdot S \na \delta_N^\eta(x) dx.  $$
Now, we remark that due to our assumptions on $f$, by elliptic regularity,  $h = {\rm St}^{-1}(S \na f)(x)$ is $C^1$ inside $\mO$. Moreover, in virtue of  Remark \eqref{rem:points_in_O},  we can assume \eqref{hyp:in_O}.  Hence, as $\eta \rightarrow 0$,   
$$ - \frac{8\pi}{3} \int_{\R^3} h(x)  \cdot S \na \delta_N^\eta(x) dx \: \rightarrow \: - \frac{8\pi}{3}  \langle  S \na \delta_N ,  h \rangle =  
 \frac{8\pi}{3} \frac{1}{N} \sum_{i=1}^N S \na \cdot h(x_i). $$

\mspace
It remains to prove \eqref{lim_N_integrale_double}. 
In the special case where $f \in C^r(\R^3)$ for some $r \in (0,1)$ (implying that it vanishes at $\pa \mO$), classical results on Calder\'on-Zygmund operators yield that the function  $\int_{\R^3}  g_S(x-y) f(x) dx = \frac{8\pi}{3} S \na \cdot h(y)$  is a continuous (even H\"older)  bounded function, so  \eqref{H1} implies straightforwardly 
$$  \int_{(\R^3 \times \R^3)\setminus \text{Diag}}   g_S(x-y) f(x) dx (d\delta_N(y) - f(y)dy) = \int_{\R^3}  \frac{8\pi}{3} S \na \cdot h(y)  (d\delta_N(y) - f(y)dy)   \rightarrow 0. $$
In the general case where $f$ is discontinuous across $\pa \mO$, the proof is a bit more involved. The difficulty lies in the fact that some points $x_i$ get closer to the boundary as $N \rightarrow +\infty$.  

\mspace
Let $\eps > 0$. Under  \eqref{H2},  there exists $c' > 0$ (depending on $c$ only) such that for $N^{-1/3} \le \eps$,
\begin{equation} \label{cardinal}
 \left| \{ i, \: x_i \: \text{belongs to the $c'  \eps$ neighborhood of $\pa \mO$}  \}\right| \le \eps N. 
 \end{equation}
Let $\chi_\eps : \R^3 \rightarrow [0,1]$ be a smooth function such that $\chi_\eps =1$ in a  $c' \eps/4$ neighborhood of $\pa \mO$, $\chi_\eps =0$ outside a  $c' \eps/2$ neighborhood of $\pa \mO$. We write 
\begin{align*}
& \int_{(\R^3 \times \R^3)\setminus \text{Diag}}   g_S(x-y) f(x) dx (d\delta_N(y) - f(y)dy) \\
& =  \int_{(\R^3 \times \R^3)\setminus \text{Diag}}   g_S(x-y) (\chi_\eps f)(x) dx (d\delta_N(y) - f(y)dy) \\ 
&+ \int_{(\R^3 \times \R^3)\setminus \text{Diag}}   g_S(x-y) ((1- \chi_\eps)f)(x) dx (d\delta_N(y) - f(y)dy). 
\end{align*}
By formula \eqref{polarization}, the second term reads 
$$   \int_{(\R^3 \times \R^3)\setminus \text{Diag}}   g_S(x-y) (1- \chi_\eps f)(x) dx (d\delta_N(y) - f(y)dy) = \frac{8\pi}{3}   \int_{\R^3} S \na \cdot u_\eps(y) \,  (d\delta_N(y) - f(y)dy) $$
with $u_\eps = {\rm St}^{-1} S  \na ((1- \chi_\eps) f)$. The source term $(1- \chi_\eps) f$ being $C^1$ and compactly supported, $S \na \cdot  u_\eps$ is H\"older  and bounded, so that, as $N \rightarrow +\infty$,  the integral goes to zero by the weak convergence assumption \eqref{H1}, for any fixed $\eps > 0$. 
As regards the first term, we split it again into 
\begin{align*}
& \int_{(\R^3 \times \R^3)\setminus \text{Diag}}   g_S(x-y) (\chi_\eps f)(x) dx (d\delta_N(y) - f(y)dy) \\ 
& =  \int_{(\R^3 \times \R^3)\setminus \text{Diag}}   g_S(x-y) (\chi_\eps f)(x) dx  \chi_\eps(y)(d\delta_N(y) - f(y)dy) \\ 
& +  \int_{(\R^3 \times \R^3)\setminus \text{Diag}}   g_S(x-y) (\chi_\eps f)(x) dx (1-\chi_\eps)(y) (d\delta_N(y) - f(y)dy) \\
& =  \frac{8\pi}{3} \int_{\R^3} S\na \cdot v_\eps(y) \chi_\eps(y)  (d\delta_N(y) - f(y)dy)  \\
& + \:   \frac{8\pi}{3}  \int_{\R^3} S\na \cdot v_\eps(y) (1-\chi_\eps)(y)  (d\delta_N(y) - f(y)dy) 
\end{align*}
where $v_\eps$ is this time the solution of the Stokes equation with source $S  \na (\chi_\eps f)$. It is H\"older away from $\pa \mO$, so that the last term at the right-hand side goes again to zero as $N \rightarrow +\infty$, by assumption \eqref{H1}. 

\mspace
It remains to handle the first term at the right-hand side. We shall show below that for a proper choice of $\chi_\eps$ one has
\begin{equation} \label{bound_veps}
 \|\na v_\eps\|_{L^\infty} \le C, \quad \text{$C$ independent of $\eps$}. 
 \end{equation}
 Taking  advantage of this fact, we write 
\begin{align*}
& \left|   \frac{8\pi}{3} \int_{\R^3} S\na \cdot v_\eps(y) \chi_\eps(y)  (d\delta_N(y) - f(y)dy) \right|  \\
& \le   \frac{8\pi}{3}  \|S\cdot \na v_\eps\|_{L^\infty(\R^3)} \left(  \frac{1}{N} \left|\{ i, \: \chi_\eps(x_i) >0 \}\right| + \|\chi_\eps f \|_{L^1} \right) \le C \eps 
\end{align*}
where we used property \eqref{cardinal} to obtain the last inequality. With this bound and the convergence to zero of the other terms for fixed $\eps$ and $N \rightarrow +\infty$,  the limit 
\eqref{lim_N_integrale_double} follows.

\mspace
We still have to show that $\na v^\eps$ is uniformly bounded in $L^\infty$ for a good choice of $\chi_\eps$.  We borrow here to the analysis of vortex patches in the Euler equation, initiated by Chemin in 2-d \cite{MR1688875}, extended by Gamblin and Saint-Raymond in 3-d \cite{MR1373741}. First,  as $\mO$ is smooth, one can find a family of five smooth divergence-free vector fields $w_1, \dots, w_5$, tangent at $\pa \mO$ and non-degenerate in the sense that 
$$ \inf_{x \in \R^3} \sum_{i \neq j} | w_i \times w_j | > 0, $$
see \cite[Proposition 3.2]{MR1373741}. We take $\chi_\eps$ in the form  $\chi(t/\eps)$, for a coordinate $t$ transverse to the boundary, meaning that $\pa_t$ is normal at $\pa \mO$. With this choice and the assumptions on $f$, one checks easily that  $\chi_\eps f$ is bounded uniformly in $\eps$ in $L^\infty(\R^3)$ and that  for all $i$,  $w_i \cdot \na (\chi_\eps f)$ is bounded uniformly in $\eps$ in $C^0(\R^3) \subset C^{r-1}(\R^3)$ for all $r \in (0,1)$. Hence, the norm $\|\chi_\eps f\|_{r,W}$ introduced in  \cite[page 395]{MR1373741}, where $W = (w_1, \dots, w_5)$, is bounded uniformly in $\eps$.

\mspace
We then  split the Stokes system 
$$ -\Delta v_\eps + \na p_\eps = S \na (\chi_\eps f), \quad \div v_\eps = 0 $$
into the equations: 
$$   \curl  v_\eps = \Omega_\eps,  \quad \div v_\eps = 0 $$
and 
$$ -\Delta \Omega_\eps  =  \curl  S \na (\chi_\eps f). $$
 Let us  show that $\pa_i\pa_j \Delta^{-1} (\chi_\eps f)$ is bounded uniformly in $\eps$ in $L^\infty$. Let $\chi \in C^\infty_c(\R^3)$, $\chi \ge 0$, $\chi=1$ near zero. Let for all $m \in \R$,  $\Lambda^m(\xi) := (\chi(\xi) + |\xi|^2)^{m/2}$. It is easily seen through Fourier transform that for all $s \in \N$
\begin{equation} \label{basses_freq}
\|  \pa_i\pa_j \chi(D) \Lambda^{-2}(D) \Delta^{-1} (\chi_\eps f) \|_{H^s} \le  C_s  \|\chi_\eps f\|_{L^2}  \le C'_s.
\end{equation}
Moreover, by the calculations in  \cite[page 401]{MR1373741}, replacing $\omega$ with $\chi_\eps f$, we get 
 \begin{equation} \label{hautes_freq}
 \|  \pa_i\pa_j  \Lambda^{-2}(D) (\chi_\eps f) \|_{L^\infty} \le  C \|\chi_\eps f \|_{L^\infty} \ln(2+\frac{\|\chi_\eps f\|_{r,W}}{ \|\chi_\eps f \|_{L^\infty}}) \le C'_r, \quad \forall  0 < r < 1.
 \end{equation}
 Combining \eqref{basses_freq} and  \eqref{hautes_freq}, we find that 
 $$ \pa_i\pa_j \Delta^{-1} (\chi_\eps f) =  \pa_i\pa_j \left( \chi(D) \Lambda^{-2}(D)  \Delta^{-1}  +  \Lambda^{-2}(D)   \right) (\chi_\eps f)$$
 is bounded uniformly in $\eps$ in $L^\infty$  and consequently, 
 $$\|\Omega_\eps \|_{L^\infty} \le  C.$$ 
 Also, by continuity of Riesz transforms over $L^p$, we have 
 $$\forall 1 < p < \infty, \quad \|\Omega_\eps \|_{L^2} \le  C_p \| \chi_\eps f\|_{L^p} \le C'_p.$$ 
Now,  applying  $w_k \cdot \na$ to the equation satisfied by $\Omega_\eps$, we obtain for all $1 \le k \le 5$,
 \begin{align}
\nonumber  -\Delta ( w_k \cdot \na \Omega_\eps) &  =  \curl  S \na ( w_k \cdot \na (\chi_\eps f)) + [ w_k \cdot \na ,   \curl  S \na ]  (\chi_\eps f) + [ w_k \cdot \na , \Delta] \Omega_\eps  \\
\label{FGH} & = \sum_{i,j} \pa_{i} \pa_j F_{i,j,\eps} +  \sum_{i} \pa_{i}  G_{i,\eps} + H_\eps 
 \end{align}
 where $F_{i,j,\eps}$, $G_{i,\eps}$ and $H_\eps$ are combinations of $\Omega_\eps$, $\chi_\eps f$ and $w_k  \cdot \na (\chi_\eps f)$. In particular, they are bounded uniformly in $\eps$ in $L^\infty \cap L^p$, for any $1 < p < \infty$. 
 
 \mspace
 For the first term at the r.h.s., we write with the same cut-off function $\chi$ as before:
  \begin{align*}
  (-\Delta)^{-1}  \sum_{i,j} \pa_{i} \pa_j F_{i,j,\eps}  = \chi(D) (-\Delta)^{-1}  \sum_{i,j} \pa_{i} \pa_j F_{i,j,\eps} + (1-\chi(D))  \sum_{i,j} \pa_{i} \pa_j F_{i,j,\eps}
  \end{align*}
 By continuity of $ (-\Delta)^{-1}   \pa_{i} \pa_j$ over $L^2$, the first term, with low frequencies, belongs to $H^s$ for any $s$, with uniform bound in $\eps$. By continuity of $(1-\chi(D)) (-\Delta)^{-1}   \pa_{i} \pa_j$ over H\"older spaces (Coifman-Meyer theorem), the second term, with high frequencies, is uniformly bounded in $\eps$ in  $C^{r-1}(\R^3)$, for any $0 < r < 1$. 

\mspace
For the second and third terms in \eqref{FGH}, we claim that 
$$ \|(-\Delta)^{-1}  \sum_{i} \pa_{i}  G_{i,\eps}\|_{L^\infty} \le C, \quad   \|(-\Delta)^{-1} H_{\eps}\|_{L^\infty} \le C. $$
 This can be seen easily by expressing these fields as $ \sum_{i} \pa_{i} \Phi \star  G_{i,\eps}$ and $\Phi \star  H_\eps$ with  $\Phi$ the fundamental solution, and by using the uniform  $L^p$ bounds on $G_{i,\eps}$ and $H_\eps$.  Eventually, we find that 
  $$ \|w_k \cdot \na \Omega_\eps\|_{C^{r-1}} \le C_r, \quad \forall 1 \le k \le 5, \quad \forall 0 < r < 1.$$   
  We conclude by \cite[Proposition 3.3]{MR1373741} that $\na v_\eps$ is bounded in $L^\infty(\R^3)$ uniformly in $\eps$.

\subsection{Smoothing}
By Proposition  \ref{prop:VN_approx_WN}, we are left with understanding  the asymptotic behaviour of 
\begin{equation} \label{def:WN}
\mW_N := \frac{75 |\mO|}{16\pi} \int_{\R^3 \times \R^3\setminus Diag} g_S(x-y)  \big(d\delta_N(x) - f(x) dx\big)   \big(d\delta_N(y) - f(y) dy\big) 
\end{equation}
The following field will play a crucial role.  For  $U,Q$ defined in \eqref{def:Uj}, we set
\begin{equation} \label{def:G_S_p_S_0}
G_S(x) \: :=  \: S_{kl} \pa_k U_l(x), \quad p_S(x) =  S_{kl} \pa_k Q_l(x).
\end{equation}
From \eqref{formula_g_S}, we have $ g_S = \frac{8\pi}{3} \, (S \na) \cdot G_S$, and $G_S$ solves in the sense of distributions
\begin{equation} \label{eq_on_G_S}
 -\Delta G_S + \na p_S =    S \na \delta, \quad \div G_S = 0 \: \text{ in } \: \R^3. 
 \end{equation}
Moreover, from the explicit expression 
$$ U_l(x) = \frac{1}{8\pi} \left( \frac{1}{|x|} e_l + \frac{x_l}{|x|^3} x \right), \quad Q_l(x)  = \frac{1}{4\pi} \frac{x_l}{|x|^3}, $$
and taking into account the fact that $S$ is symmetric and trace-free, we get 
\begin{equation} \label{def:G_S_p_S}  
G_S(x)  =  - \frac{3}{8\pi} S_{kl} x_l x_k \frac{x}{|x|^5} =   - \frac{3}{8\pi}  (Sx\cdot x) \frac{x}{|x|^5}, \quad p_S(x) = - \frac{3}{4\pi} \frac{(Sx \cdot x)}{|x|^5
} 
\end{equation}
Let us note that $G_S$ is called a point stresslet in the literature, see \cite{GuMo}. It can be interpreted as the velocity field created in a fluid of viscosity $1$ by a point particle  whose resistance to a strain is given by $-S$.

\mspace
We now come back to the analysis of \eqref{def:WN}. Formal replacement of the function $f$ in Lemma  \ref{lem:sign:gA} by $\delta_N - f$ yields the formula 
\begin{equation} \label{formal_formula2}
  " \int_{\R^3 \times \R^3} g_S(x-y)  \big(d\delta_N(x) - f(x) dx\big)   \big(d\delta_N(y) - f(y) dy\big) = - \frac{16\pi}{3N^2} \int_{\R^3} |D(h_N)|^2 " 
  \end{equation}
  where 
\begin{equation} \label{def:hN}
 h_N(x) :=  \sum_{i=1}^N G_S(x-x_i) - N {\rm St}^{-1} (S \na f) =  \sum_{i=1}^N G_S(x-x_i) - N \int_{\R^3} G_S(x-y) f(y) dy 
 \end{equation}
satisfies 
\begin{equation} \label{eq:hN}
- \Delta h_N + \na q_N = S\na \sum_i \delta_{x_i} \: - \: N S\na f, \quad \div h_N = 0 \quad \text{in } \R^3. 
\end{equation}
The formula \eqref{formal_formula2} is similar to the formula \eqref{formal_formula1}, and is as much abusive, as both sides are infinite. Still, by an appropriate regularization of the source term   $S\na \sum_i \delta_{x_i}$, we shall be able in the end to obtain a rigorous formula, convenient for the study of $\mW_N$. This regularization process is the purpose of the present paragraph.

\mspace

\mspace
For any $\eta > 0$, we denote $B_\eta = B(0,\eta)$, and  define $G_S^\eta$ by:
\begin{equation} \label{def_G_S_eta1}
G_S^\eta = G_S, \: p_S^\eta  = p_S \: \text{ outside }  B_\eta, 
\end{equation}
\begin{equation} \label{def_G_S_eta2}
- \Delta G_S^\eta + \na p_S^\eta  = 0, \quad \div G_S^\eta = 0, \quad G_S^\eta\vert_{\pa B_\eta} = G_S\vert_{\pa B_\eta} \: \text{ in } B_\eta.
\end{equation}
Note that by homogeneity,
\begin{equation} \label{homogeneity_GA}
 G_S^\eta(x) = \frac{1}{\eta^2} G_S^1(x/\eta). 
 \end{equation}  
The field $G_S^\eta$ belongs to $\dot{H}^1(\R^3)$, and solves  
\begin{equation} \label{eq_on_G_S_eta}
- \Delta G_S^\eta + \na p_S^\eta = S^\eta 
\end{equation}
where $S^\eta$ is the measure on the sphere defined by 
\begin{equation} \label{def:Seta}  
S^\eta :=  - \left[2D(G_S^\eta) n - p_S^\eta n \right] \, s^\eta  = - \left[\pa_n G_S^\eta  - p_S^\eta n \right] \, s^\eta
\end{equation}
with $n=\frac{x}{|x|}$ the unit normal vector  pointing {\em outward} $B_\eta$, $\left[F\right] := F\vert_{\pa B_\eta^+} - F\vert_{\pa B_\eta^-}$ the jump at $\pa B_\eta$ (with $\pa B_\eta^+$, resp. $\pa B_\eta^-$, the outer, resp. inner boundary of the ball),  and $s^\eta$  the standard surface measure on $\partial B_\eta$.  We claim the following 
\begin{lemma}  \label{lem:smoothing}
For all $\eta > 0$, 
$\: S^\eta = {\rm div} \, \Psi^\eta \:$ in $\R^3$, where  
\begin{equation} \label{def:Psi_eta}
\begin{aligned}
 \Psi^\eta & :=   \frac{3}{\pi \eta^5} \left( Sx \otimes x + x \otimes Sx - 5\frac{|x|^2}{2} S + \frac{5}{4} \eta^2 S \right) -  2 D(G^\eta_S)(x) + p^\eta_S(x) {\rm Id}, \quad x \in B_\eta,\\ 
 \Psi^\eta & := 0 \text{ outside}. 
 \end{aligned}
 \end{equation}
 Moreover, $\Psi^\eta \rightarrow S\delta$ in the sense of distributions as $\eta \rightarrow 0$, so that $S^\eta \rightarrow S\na \delta$. 
 \end{lemma}
 
 \mspace
{\em Proof of the lemma}. From the explicit formula \eqref{def:G_S_p_S} for $G_S$ and $p_S$, we find 
$$ 2D(G_S) = -\frac{3}{4\pi} \frac{Sx \otimes x + x \otimes Sx}{|x|^5} + \frac{15}{4\pi} \frac{(Sx \cdot x) x \otimes x}{|x|^7} - \frac{3}{4\pi} \frac{Sx\cdot x}{|x|^5} {\rm Id}. $$
so that 
\begin{align} \label{expression_tensor_ext}
(2D(G_S^\eta)n - p^\eta_S n)\vert_{\pa B_\eta^+}  & = (2D(G_S)n  - p_S n)\vert_{\pa B_\eta^+} =  \frac{3}{4\pi |\eta|^3} \left( 4 (Sn \cdot n) n -  Sn \right) 
\end{align} 
Using that $S$ is trace-free, one can check from definition \eqref{def:Psi_eta} that $\div \Psi^\eta = 0$ in the complement of $\pa B_\eta$, while:
\begin{align*}
[\Psi^\eta n] & = -\Psi^\eta n\vert_{\pa B_\eta^-} \\
& = \frac{3}{\pi \eta^3} \Big( (Sn \otimes n)n + (n \otimes Sn)n - \frac{5}{4} Sn \Big)  -  (2 D(G^\eta_S)n + p^\eta_S n)\vert_{\pa B_\eta^-} \\
& = (2D(G_S^\eta)n  - p^\eta_S n)\vert_{\pa B_\eta^+}   -  (2 D(G^\eta_S)n + p^\eta_S n)\vert_{\pa B_\eta^-} 
\end{align*}
where the last equality comes from \eqref{expression_tensor_ext}. Together with \eqref{def:Seta}, it implies the first claim of the lemma.  

\mspace
To compute the limit of $\Psi^\eta$ as $\eta \rightarrow 0$, we write $\Psi^\eta = \Psi^\eta_1 + \Psi^\eta_2$, with 
$$ \Psi_1^\eta = \frac{3}{\pi \eta^5} \left( Sx \otimes x + x \otimes Sx - 5\frac{|x|^2}{2} S + \frac{5}{4} \eta^2 S \right), \quad \Psi_2^\eta =  -  2 D(G^\eta_S)(x) + p^\eta_S(x) {\rm Id}. $$
Let $\varphi \in C^\infty_c(\R^3)$ a test function. We can write 
$ \langle \Psi_1^\eta , \varphi \rangle =  \langle \Psi_1^\eta , \varphi(0) \rangle + \langle \Psi_1^\eta , \varphi -  \varphi(0)  \rangle. $
The second term  is $O(\eta)$, while the first term can be computed using the elementary formula 
$\int_{B_1} x_i x_j dx = \frac{4\pi}{15} \delta_{ij}$. We find 
\begin{equation} \label{lim_Psi1}
 \lim_{\eta \rightarrow 0}    \langle\Psi_1^\eta , \varphi \rangle = \frac{3}{5} S \varphi(0) =  \langle \frac{3}{5}S \delta , \varphi \rangle. 
 \end{equation}
For the second term, using the homogeneity \eqref{homogeneity_GA}, we find again that $\lim_{\eta} \langle \Psi_2^\eta , \varphi \rangle =  \langle \Psi_2^1 , \varphi(0) \rangle$. Note that the pressure $p_S^1$ is defined up to a constant, so that we can always select the one with zero average.  With this choice, we find 
\begin{equation} \label{lim_Psi2}
\begin{aligned}
 \langle \Psi_2^1 , \varphi(0) \rangle & = \int_{B_\eta}  \big( -  2 D(G^\eta_S) + p^\eta_S{\rm Id} \big) \varphi(0) = - 2 \int_{B_1} D(G^1_S) \: \varphi(0) \\
 & = -  \int_{\pa B_1} (n \otimes G^1_S + G^1_S \otimes n) \: \varphi(0) =   - \int_{\pa B_1}  (n \otimes G_S + G_S \otimes n) \: \varphi(0) \\ 
 & =  \frac{3}{4\pi} \int_{\pa B_1} (S n \cdot n) n \otimes n \: \varphi(0) = \frac{2}{5}S \varphi(0) =  \langle \frac{2}{5}S \delta , \varphi \rangle. 
 \end{aligned} 
 \end{equation}
 where the sixth equality comes from the elementary formula $\int_{\pa B_1} n_i n_j n_k n_l  ds^1=  \frac{4\pi}{15} (\delta_{ij} \delta_{kl} + \delta_{ik} \delta_{jl} + \delta_{il} \delta_{jk})$. 
The result follows.

 \mspace
 For later purpose, we also prove here the
 \begin{lemma} \label{lem:int_GA}
  $$ \int_{\pa B_\eta} G_S^\eta d S^\eta = \int_{\pa B_\eta} G_S d S^\eta =  \frac{1}{\eta^3}\left(  \int_{B^1} |\na G_S^1|^2 + \frac{3}{10\pi} |S|^2 \right).$$
  \end{lemma}
 {\em Proof.} 
 \begin{align*}
 \int_{\pa B_\eta} G_S^\eta d S^\eta & = \int_{\pa B_\eta}  G_S^\eta \left( \pa_n G_S^\eta - p_S n \right)\vert_{\pa B_\eta^-} ds^\eta & - \int_{\pa B_\eta}  G_S^\eta \left( \pa_n G_S^\eta - p_S n \right)\vert_{\pa B_\eta^+} ds^\eta \\ 
& =  \int_{B_\eta} |\na G^\eta_S|^2 dx & -  \int_{\pa B_\eta}  G_S \left( \pa_r G_S - p_S e_r \right)\vert_{\pa B_\eta} ds^\eta 
\end{align*}
By \eqref{homogeneity_GA},   $\int_{B_\eta} |\na G^\eta_S|^2 dx = \frac{1}{\eta^3} \int_{B_1} |\na G^1_S|^2 dx$. The second term can be computed with \eqref{def:G_S_p_S}: \begin{align*}  
    & \int_{\pa B_\eta}  G_S \left( \pa_r G_S - p_S e_r \right)\vert_{\pa B_\eta} ds^\eta =   \int_{\pa B_\eta} \left( -\frac{3}{8\pi \eta^2} (S n \cdot n ) n \right) \, \left(   \frac{3}{2\pi \eta^3} (S n \cdot n) n \right) ds^\eta \\
 = & - \frac{9}{16\pi^2 \eta^3} \int_{\pa B_1} (S n \cdot n)^2 ds^1 = -\frac{3}{10\pi} |S|^2.    
 \end{align*}
 
 \subsection{The renormalized energy}
 Thanks to the regularization of $S\na\delta$ introduced in the previous paragraph, {\it cf.} Lemma \ref{lem:smoothing}, we shall be able to set a rigorous alternative to the abusive formula \eqref{formal_formula2}. Specifically, we shall state an identity involving $\mW_N$, defined in \eqref{def:WN},  and the energy of the function 
 \begin{equation} \label{def:hNeta}
 h_N^\eta(x) :=  \sum_{i=1}^N G_S^\eta(x-x_i) + N {\rm St}^{-1} (S \na f) =  \sum_{i=1}^N G_S^\eta(x-x_i) - N \int_{\R^3} G_S(x-y) f(y) dy. 
 \end{equation}
This function solves
 \begin{equation} \label{eq:hNeta}
  - \Delta h_N^\eta + \na p_N^\eta =  \sum_{i=1}^N S^\eta(x-x_i) - N S \na f, \quad \div h_N^\eta = 0
  \end{equation} 
and is a regularization  of $h_N$, {\it cf.} \eqref{def:hN}-\eqref{eq:hN}.

\mspace
The main result of this section is the 
 \begin{proposition}
\begin{equation} \label{asymptote_WN}
\mW_N = - \frac{25|\mO|}{2N^2}  \lim_{\eta \rightarrow 0} \left( \int_{\R^3} | \na h_{N}^\eta |^2  - \frac{N}{\eta^3}( \int_{B^1} |\na G_S^1|^2 + \frac{3}{10\pi} |S|^2 ) \right). 
\end{equation}
\end{proposition}
 {\em Proof.} We assume that $\eta$ is small enough so that $2 \eta < \min_{i \neq j} |x_i - x_j|$. From the explicit expressions \eqref{def:hN}, \eqref{def:hNeta}, we find that $h_N, h_N^\eta = O(|x|^{-2})$,  $\na (h_N, h_N^\eta) = O(|x|^{-3})$ and $p_N, p_N^\eta = O(|x|^{-3})$  at infinity. As these quantities decay fast enough,  we can perform an integration by parts to find 
\begin{align*}
 \int_{\R^3} | \na h_{N}^\eta |^2 
 & =  \langle -\Delta h_N^\eta , h_N^\eta \rangle = \langle -\Delta h_N^\eta + \na p_N^\eta , h_N^\eta \rangle \\ 
 & = \langle \sum_i S^\eta(x-x_i) - N S \na f , h_N \rangle \: +  \: \langle \sum_i S^\eta(x-x_i) - N S \na f , h_N^\eta - h_N \rangle \\
 & =  \sum_i \langle  S^\eta(x-x_i), h_N^i  \rangle +    \sum_i \langle  S^\eta(x-x_i), G_S(x-x_i)  \rangle   \\
 & -  \langle  N S \na f , h_N \rangle \: +  \: \langle \sum_i S^\eta(x-x_i) - N S \na f , h_N^\eta - h_N \rangle  =: a + b + c + d, 
\end{align*}
where we defined $h_N^i := h_N -  G_S(x-x_i)$. 

\mspace
As $h_N^i$ is smooth over the support of $S^\eta(\cdot - x_i)$, we can apply Lemma \ref{lem:smoothing} to obtain
$$ \lim_{\eta \rightarrow 0} a = - \sum_i S \na \cdot h^i_N(x_i).  $$
We can then apply Lemma \ref{lem:int_GA} to obtain
$$  b =  \frac{N}{\eta^3}  ( \int_{B^1} |\na G_S^1|^2 + \frac{3}{10\pi} |S|^2 ). $$
As regards the fourth term, we notice that by our definition \eqref{def_G_S_eta1}-\eqref{def_G_S_eta2} of $G_S^\eta$, and the fact that the balls $B(x_i, \eta)$ are disjoint,  the function 
$h_N - h_N^\eta = \sum_{i} (G_S(x-x_i) - G_S^\eta(x-x_i))$ is zero over $\cup_i \pa B(x_i, \eta)$, which is the support of $\sum_i S^\eta(x-x_i)$. It follows that 
\begin{align*} 
d   =  -N  \langle  S \na f  , h_N^\eta - h_N \rangle & = N \sum_i \int_{B(x_i, \eta)} S\na \cdot G_S^\eta(x-x_i) \left( f(x) - f(x_i) \right) dx \\ 
& - N \sum_i \int_{B(x_i, \eta)} S\na \cdot G_S(x-x_i)  \left( f(x) - f(x_i) \right)dx 
\end{align*}
where we integrated by parts, using that $G_S - G_S^\eta$ is zero outside the balls. Let us  notice that the second integral at the right-hand side converges despite the singularity of  $S\na \cdot G_S$, using the smoothness of $f$ near $x_i$ (by assumption \eqref{hyp:in_O} and  Remark \ref{rem:points_in_O}). Moreover, it goes to zero as $\eta \rightarrow 0$. Using the homogeneity and smoothness properties of $G_S^\eta$ inside $B^\eta$, we also find that the first sum goes to zero with $\eta$, resulting in 
$$  \lim_{\eta \rightarrow 0} d =  0. $$
We end up with 
$$ \lim_{\eta\rightarrow 0 } \left( \int_{\R^3} | \na h_{N}^\eta |^2  - \frac{N}{\eta^3}  ( \int_{B^1} |\na G_S^1|^2 + \frac{3}{10\pi} |S|^2 ) \right) = - \sum_i S \na \cdot h^i_N(x_i)  -  \langle  N S \na f , h_N \rangle $$
 It remains to rewrite properly the right-hand side: we first get 
 \begin{align*}
 - \sum_i S \na \cdot  h^i_N(x_i) & = - \sum_{i \neq j} S \na \cdot G_S(x_i - x_j)  + N  \sum_{i}  \int_{\R^3}  S \na \cdot G_S(x_i-y) f(y) dy \\
  & = - \frac{3N^2}{8\pi}  \int_{\R^3 \times \R^3 \setminus \textrm{Diag}} g_S(x-y) d \delta_N(x) (d\delta_N(y) - f(y) dy) 
  \end{align*}
 and integrating by parts 
 \begin{align*}
 -  \langle  N S \na f , h_N \rangle & = N \int_{\R^3}S \na \cdot h_N(x) f(x) dx \\ 
 & = N  \int_{\R^3} \left(  \sum_i S \na \cdot G_S(x-x_i) - N\int_{\R^3} S \na \cdot G_S(x - y) f(y) dy \right) f(x) dx \\
 & = \frac{3N^2}{8\pi}   \int_{\R^3 \times \R^3}  g_S(x-y) f(x)   dx (d\delta_N(y) - f(y) dy)  dx.
 \end{align*}
The last equality was deduced from the identity $g_S = \frac{8\pi}{3} \, (S \na) \cdot G_S$, see the line after \eqref{def:G_S_p_S_0}. The proposition follows. 
 
 \mspace
 We can refine the previous proposition by the following
 \begin{proposition} \label{prop:h_N_eta_h_N_alpha}
 Let $c > 0$ the constant in \eqref{H2}. There exists $C > 0$ such that: for all  $\alpha < \eta <  \frac{c}{2} N^{-1/3}$, 
\begin{align*}
 \Bigl| \int_{\R^3} | \na h_{N}^\eta |^2 - \int_{\R^3} | \na h_{N}^\alpha |^2 -  N \left( \frac{1}{\eta^3} -  \frac{1}{\alpha^3} \right) ( \int_{B^1} |\na G_S^1|^2 + \frac{3}{10\pi} |S|^2 ) \Bigr| \le C N^2 \eta.
 \end{align*}
 \end{proposition}
 
 \mspace
{\em Proof.}  One has from \eqref{def:hNeta}
$$  h_{N}^\eta  = h_{N}^\alpha + \sum_{i=1}^N (G_S^\eta - G_S^\alpha)(x-x_i). $$
It follows that 
\begin{align*}
\int_{\R^3} | \na h_{N}^\eta |^2 - \int_{\R^3} | \na h_{N}^\alpha |^2 & = \sum_{i,j} \int_{\R^3 }\na (G_S^\eta - G_S^\alpha)(x-x_i) :  \na (G_S^\eta - G_S^\alpha)(x-x_j) \\
& +  2  \sum_{i}  \int_{\R^3} \na h_N^\alpha :  \na(G_S^\eta - G_S^\alpha)(x-x_i) 
\end{align*}
After integration by parts,  
\begin{align*}
\int_{\R^3 }\na (G_S^\eta - G_S^\alpha)(\cdot-x_i) :  \na (G_S^\eta- G_S^\alpha)(\cdot-x_j) =  \langle (S^\eta - S^\alpha)(\cdot -x_i) , (G_S^\eta - G_S^\alpha)(\cdot - x_j) \rangle. 
\end{align*}
while 
\begin{align*}
\int_{\R^3 }  \na h_N^\alpha :  \na(G_S^\eta - G_S^\alpha)(x-x_i)  = \langle \sum_j S^\alpha(\cdot - x_i) - N S \na f , (G_S^\eta - G_S^\alpha)(\cdot-x_i)  \rangle.
\end{align*}
We get 
\begin{align}
& \int_{\R^3} | \na h_{N}^\eta |^2 - \int_{\R^3} | \na h_{N}^\alpha |^2  =  \sum_{i \neq j}   \langle (S^\alpha + S^\eta)(\cdot -x_i) , (G_S^\eta - G_S^\alpha)(\cdot - x_j) \rangle \nonumber \\ 
& - 2 \sum_i N \langle  S \na f , (G_S^\eta - G_S^\alpha)(\cdot-x_i)  \rangle \: + \:  N   \langle (S^\alpha + S^\eta) , (G_S^\eta - G_S^\alpha) \rangle  =: a + b + c  \label{iiiiii}
\end{align}
We note that $G_S^\eta - G_S^\alpha$ is zero outside $B_\eta$,  while $S^\alpha + S^\eta$ is supported in $B_\eta$. Moreover, thanks to \eqref{H2},  for $\alpha < \eta  < \frac{c}{2}$, the balls $B(x_i, \eta)$ are disjoint. We deduce: $a=0$. 

\mspace
After integration by parts, taking into account that $G_S^\eta- G_S^\alpha$ vanishes outside $B_\eta$, we can write  $b =  b_\eta - b_\alpha$ with 
\begin{align*}
 b_\alpha & :=  2\sum_i N \int_{B(x_i, \eta)}  S \na \cdot G_S^\alpha(\cdot-x_i) \, (f - f(x_i))  \\ 
 b_\eta & :=  2\sum_i N \int_{B(x_i,\eta)}  S \na \cdot G_S^\eta(\cdot-x_i) \, (f - f(x_i)).  
 \end{align*}
By assumption \eqref{hyp:in_O}, for $N$ large enough, for all $1\le i \le N$ and all $\eta \le \frac{c}{2} N^{-1/3}$, $B(x_i, \eta)$ is included in $\mO$.  Hence, $f$ is $C^{1}$ in $B(x_i, \eta)$, and 
$$  \bigl| \int_{B(x_i, \eta)}  S \na \cdot G_S^\eta(\cdot-x_i)  \, (f - f(x_i) \bigr|  \le \frac{C}{\eta^3} \|\na f\vert_{\mO}\|_{\infty}  \int_{B(x_i, \eta)} |x-x_i| dx   \le C \eta.   $$
This results in: $b_\eta \le C N^2  \eta$. 

\mspace
Similarly, decomposing $B(x_i, \eta) = B(x_i, \alpha) \cup \Big( B(x_i, \eta) \setminus B(x_i, \alpha) \Big)$, we find 
$$ \bigl| \int_{B(x_i, \eta)}  S \na \cdot G_S^\alpha(\cdot-x_i) \, (f - f(x_i))\bigr| \le C\left( \alpha+  \int_{B(x_i, \eta)} \frac{1}{|x-x_i|^2}  dx\right) \le C' \eta $$
using again that $f$ is Lipschitz over $B(x_i, \eta)$.  We end up with $b_\alpha \le C N^2 \eta$, and finally $b \le C N^2 \eta$. 

\mspace
For the last term $c$ in \eqref{iiiiii}, we first notice that as $G_S^\eta - G_S^\alpha$ is zero outside $B_\eta$:
\begin{align} 
\label{eq:S_alpha_diff_g_S}
& \langle (S^\alpha + S^\eta) , (G_S^\eta - G_S^\alpha) \rangle = \langle S^\alpha , (G_S^\eta - G_S^\alpha) \rangle \\
\nonumber
= &   \langle S^\alpha , G_S^\eta \rangle \: - \:    \langle S^\alpha , G_S \rangle \\
\nonumber
= &  \langle S^\alpha , G_S^\eta \rangle \: - \:  \frac{1}{\alpha^3}\left(  \int_{B^1} |\na G_S^1|^2 + \frac{3}{10\pi} |S|^2 \right)
\end{align}
where we used Lemma \ref{lem:int_GA} in the last line. By the definition of $S^\alpha$, the remaining term splits  into 
\begin{align*}
\langle S^\alpha , G_S^\eta \rangle  = - \int_{\pa B_\alpha^+} \left( \pa_r G_S - p_S e_r \right) \cdot G_S^\eta ds^\alpha  \: + \: \int_{\pa B_\alpha^-} \left( \pa_r G_S^\alpha- p_S^\alpha e_r \right) \cdot G_S^\eta ds^\alpha
\end{align*}
By integration by parts, applied in $B_\eta\setminus B_\alpha$ for the first term and  in $B_\alpha$ for the second term, we get
\begin{align*}
\langle S^\alpha , G_S^\eta \rangle  & = -  \int_{\pa B_\eta^-} \left( \pa_r G_S - p_S e_r \right) \cdot G_S^\eta  ds^\eta + \int_{B_\eta\setminus B_\alpha} \na G_S : \na G_S^\eta  +   \int_{B_\alpha} \na G_S^\alpha: \na G_S^\eta \\ 
& = -  \int_{\pa B_\eta} \left( \pa_r G_S - p_S e_r \right) \cdot G_S^\eta  ds^\eta + \int_{B_\eta} \na G_S^\alpha\cdot \na G_S^\eta \\
& = -  \int_{\pa B_\eta} \left( \pa_r G_S - p_S e_r \right) \cdot G_S^\eta  ds^\eta  +  \int_{\pa B_\eta^-} G_S^\alpha\cdot \left( \pa_r G_S^\eta - p_S^\eta e_r \right) \\
& = \langle  S^\eta , G_S \rangle =  \frac{1}{\eta^3}\left(  \int_{B^1} |\na G_S^1|^2 + \frac{3}{10\pi} |S|^2 \right)
\end{align*}
From there, the conclusion  follows easily. 

\mspace
If we let $\alpha \rightarrow 0$ in Proposition \ref{prop:h_N_eta_h_N_alpha}, combining with Propositions \ref{asymptote_WN}  and \ref{prop:VN_approx_WN}, we find  
\begin{corollary} \label{cor:asymptote:WN}
For all  $ \eta <  \frac{c}{2} N^{-1/3}$, 
\begin{align*}
 \Bigl| \mathcal{V}_N + \frac{25|\mO|}{2N^2} \Bigl(  \int_{\R^3} | \na h_{N}^\eta |^2 - \frac{N}{\eta^3} ( \int_{B^1} |\na G_S^1|^2 + \frac{3}{10\pi} |S|^2 ) \Bigr) \Bigr| \le  \eps(N) 
 \end{align*}
where $\eps(N) \rightarrow 0$ as $N \rightarrow +\infty$.  
\end{corollary}
\noindent
This corollary shows that to understand the limit of $\mV_N$, it is enough to study the limit of 
$$ \frac{25|\mO|}{2N^2} \Bigl(  \int_{\R^3} | \na h_{N}^{\eta_N} |^2 - \frac{N}{\eta_N^3} ( \int_{B^1} |\na G_S^1|^2 + \frac{3}{10\pi} |S|^2 ) \Bigr)$$
for $\eta_N := \eta N^{-1/3}$, $\eta < \frac{c}{2}$ fixed.  For periodic and more general stationary point processes, this will be possible through an homogenization approach. This homogenization approach involves an analogue of a  cell equation, called jellium in the literature on Coulomb gases. We will motivate and introduce this system in the next section.

\section{Blown-up system} \label{sec3}
Formula \eqref{asymptote_WN} suggests to understand at first the behaviour  of  $\int_{\R^3} |\na h_N^\eta|^2$ at fixed $\eta$, when $N \rightarrow +\infty$.  To analyze the system  \eqref{eq:hNeta}, a useful intuition can be taken from classical homogenization problems of the form 
\begin{equation} \label{periodic_homogenization}
- \Delta h_\eps + \na p_\eps =  S \na \left(  \frac{1}{\eps^3} F(x,x/\eps)  - \frac{1}{\eps^3}  \overline{F}(x) \right),  \: \div h_\eps = 0 \: \text{ in a  domain } \Omega, \quad h_\eps\vert_{\pa \Omega} = 0,
\end{equation}
with  $F(x,y)$  periodic in variable $y$, and $\overline{F}(x) := \int_{\T} F(x,y) dy$. Roughly, $\Omega$ would be like $\mO$,  the small scale $\eps$ like  $N^{-1/3}$, the term $ \frac{1}{\eps^3}  F(x,x/\eps)$ would correspond to the  sum of (regularized) Dirac masses, while the term $ \frac{1}{\eps^3}  \overline{F}$ would be an analogue of $N f$. The factor $\frac{1}{\eps^3}$ in front of $F$ is put consistently with the fact that $\sum_i \delta_{x_i}$ has mass $N$. 
The dependence on $x$ of the source term in \eqref{periodic_homogenization} mimics the possible macroscopic inhomogeneity of the point distribution $\{x_i\}$. 

\mspace
In the much simpler model \eqref{periodic_homogenization}, standard arguments show that $h_\eps$  behaves like 
\begin{equation} \label{asympt_h_eps} 
h_\eps(x) \approx  \frac{1}{\eps^2} H(x,x/\eps) 
\end{equation}
where $H(x,y)$ satisfies the cell problem 
$$ - \Delta_y H(x,\cdot) + \na_y P(x,\cdot) = S \na_y F(x,\cdot), \quad {\rm div}_y H(x,\cdot) = 0, \quad y \in \T^3. $$ 
Let us stress that substracting the term  $\frac{1}{\eps^3}  \overline{F}(x)$ in the source term of \eqref{periodic_homogenization} is crucial for the asymptotics \eqref{asympt_h_eps} to hold.  It follows that 
$$ \eps^6 \int_{\Omega} |\na h_\eps|^2 \approx  \int_{\Omega} |\na_y H(x,x/\eps)|^2 dx \xrightarrow[\eps \rightarrow 0]{} \int_{\Omega} \int_{\T^3} |\na_y H(x,y)|^2 dy dx.  $$
Note that the factor $\eps^6$ in front of the left-hand side is coherent with the factor $\frac{1}{N^2}$ at the right-hand side of \eqref{asymptote_WN}.
Note also that  
$$  \int_{\T^3} |\na_y H(x,y)|^2 dy  = \lim_{R \rightarrow +\infty} \frac{1}{R^3} \int_{(-R,R)^3} |\na_y H(x,y)|^2 dy,  $$
Such average over larger and larger boxes may be still meaningful in more general settings, typically in stochastic homogenization. 

\mspace
Inspired by those remarks, and back to  system \eqref{eq:hNeta}, the hope is that  some  homogenization process may take place, at least locally near each $x \in \mO$. More precisely, we hope to recover  $\lim_{N} \mathcal{W}_N$ by summing over $x \in \mO$ some microscopic energy, locally averaged around $x$. This microscopic energy will be deduced from an analogue of the cell problem,  called a {\em jellium}  in the literature on the Ginzburg-Landau model and Coulomb gases. 

\subsection{Setting of the problem} \label{subsec:setting}
We will call {\em point distribution} a locally finite subset of $\R^3$. Given a point distribution  $\Lambda$, we consider the following problem in $\R^3$ 
\begin{equation} \label{eq:jellium}
\begin{aligned} 
 - \Delta H + \na P & = \sum_{z \in \Lambda} S \na \delta_{-z} \\
 \div H & = 0. 
 \end{aligned}
 \end{equation}
Given a solution $H = H(y)$, $P = P(y)$, we introduce for any $\eta > 0$
\begin{equation} \label{def:H_H_eta}
H^\eta := H + \sum_{z \in \Lambda} (G_S^\eta - G_S)(\cdot +z)
\end{equation}
which satisfies by \eqref{eq_on_G_S}, \eqref{eq_on_G_S_eta}: 
\begin{equation} \label{eq:jellium_eta}
\begin{aligned} 
 - \Delta H^\eta + \na P^\eta & = \sum_{z \in \Lambda} S^\eta(\cdot+z) \\
 \div H^\eta  & = 0. 
 \end{aligned}
 \end{equation}
We remark that, the set $\Lambda$ being locally finite, the sum at the right-hand side of \eqref{eq:jellium} or \eqref{eq:jellium_eta}  is well-defined as a distribution. Also, the sum at the right-hand side of  \eqref{def:H_H_eta} is well-defined pointwise, because $G_S^\eta - G_S$ is supported in $B_\eta$. 

\mspace
 As discussed at the beginning of Section \ref{sec3}, we expect the limit of $\int_{\R^3} |\na h_N^\eta|^2$ to be described in terms of quantities of the form 
 $$ \lim_{R \rightarrow +\infty} \frac{1}{R^3}\int_{K_R} |\na H^\eta(y)|^2 \, dy  $$ 
where $K_R := (-\frac{R}{2},\frac{R}{2})^3$,   for various $\Lambda$ and  solutions $H^\eta$ of  \eqref{eq:jellium_eta}. Broadly, the energy concentrated locally around a point $x$ should be understood from a blow-up of the original system \eqref{eq:hNeta}, zooming at scale $N^{-1/3}$ around $x$. Let  $x \in \mO$ (the center of the blow-up), and $\eta_N := \eta N^{-1/3}$, for a fixed $\eta > 0$. If we introduce 
\begin{equation} \label{eq:rescaling}
\begin{aligned}
& H^{\eta}_N(y) := N^{-2/3} h_N^{\eta_N}(x+N^{-1/3} y), \quad P^\eta_N(y) := N^{-1} p_N^{\eta_N}(x+N^{-1/3} y), \\
&  z_{i,N} := N^{1/3} (x-x_{i,N}) 
\end{aligned}
\end{equation}
  we find that 
\begin{equation} \label{eq:zoom_Stokes}
- \Delta  H_N^\eta + \na P_N^\eta = \sum_{i=1}^N S^{\eta}(\cdot + z_{i,N}) - N^{-1/3} S\na_x f(x+N^{-1/3}y), \quad \div H_N^\eta = 0. 
\end{equation}
System  \eqref{eq:jellium_eta} corresponds to a formal asymptotics where one replaces $\sum_{i=1}^N \delta_{z_{i,N}}$ by  $\sum_{i=1}^\infty \delta_{z_i}$, with $\Lambda = \{z_i\}$ a point distribution. Note that, under \eqref{H2}, we expect this point distribution to  be {\em well-separated}, meaning  that there is $c >0$ such that: for all $z' \neq z \in \Lambda$,  $|z'-z| \ge c$. Still, we insist that this asymptotics is purely formal and requires much more to be made rigorous. Such rigorous asymptotics will be carried in Section \ref{sec4} for various classes of point configurations. 


\mspace
We now collect several general remarks on  the  blown-up system \eqref{eq:jellium}. We start by defining a renormalized energy. For any $L > 0$, we denote $K_L := (-\frac{L}{2}, \frac{L}{2})^3$. 
\begin{definition} 
Given a point distribution $\Lambda$, we say that a solution $H$ of \eqref{eq:jellium} is admissible if for all $\eta > 0$, the field  $H^\eta$ defined by \eqref{def:H_H_eta} satisfies  $\na H^\eta \in L^2_{loc}(\R^3)$.  

\mspace
Given an admissible solution $H$ and $\eta > 0$, we say that $H^\eta$ is of finite renormalized energy if  
$$  \mathcal{W}^\eta(\na H) :=  -\lim_{R \rightarrow +\infty} \frac{1}{R^3}  \left( \int_{K_R} |\na H^\eta|^2  -\frac{1}{\eta^3}\left| \Lambda \cap K_R \right| \Bigl( \int_{B^1} |\na G_S^1|^2 + \frac{3}{10\pi} |S|^2 \Bigr) \right) $$
exists in $\R$. We say that $H$ is of finite renormalized energy if $H^\eta$ is  for all $\eta$, and 
$$\mathcal{W}(\na H) :=  \lim_{\eta \rightarrow 0}  \mathcal{W}^\eta(\na H)$$
exists in $\R$. 
\end{definition}
\begin{remark}
From formula \eqref{def:H_H_eta}, it is easily seen that $H$ is admissible if and only if there exists one $\eta > 0$ with $\na H^\eta \in L^2_{loc}(\R^3)$. 
\end{remark}

\begin{proposition}  \label{jellium_uniqueness}
If $H_1$ and $H_2$ are admissible solutions of \eqref{eq:jellium} satisfying for some $\eta > 0$:
$$  \limsup_{R \rightarrow +\infty} \frac{1}{R^3} \int_{K_R} |\na H^\eta_1|^2 < +\infty, \quad   \limsup_{R \rightarrow +\infty} \frac{1}{R^3} \int_{K_R} |\na H^\eta_2|^2 < +\infty $$
then $\na H_1$ and $\na H_2$ differ from a constant matrix.  
\end{proposition}

\mspace
{\em Proof.} We set  $H  := H_1 - H_2 =   H_1^\eta - H_2^\eta$. It is a solution of the homogeneous Stokes equation with 
$$\limsup_{R \rightarrow +\infty} \frac{1}{R^3} \int_{K_R}|\na H|^2 < +\infty.$$
By standard elliptic regularity, any solution $v$ of the Stokes equation in the unit ball: 
$$ - \Delta v + \na p  = 0, \quad \div v = 0 \: \text{ in } \: B(0,1) $$
satisfies for some absolute constant $C$, 
$$ |\na^2 v(0)| \le  C \|\na v\|_{L^2(B(0,1))} . $$
We apply this inequality to $v(x) = H(x_0 + R x)$,  $x_0$ arbitrary. After rescaling, we find that  
$$ |\na^2 H(x_0)| \le \frac{C}{R} \Big( \frac{1}{R^{3/2}} \|\na H(x_0 + \cdot)\|_{L^2(B(0,R))}\Big). $$
As $R \rightarrow +\infty$, the right hand-side goes to zero, which concludes the proof. 
\begin{proposition} \label{prop:stationarity:Weta}
Let $\Lambda$ be a well-separated point distribution, meaning there exists $c > 0$ such that for all $z' \neq z \in \Lambda$,  $|z'-z| \ge c$.  Let $0 < \alpha <  \eta < \frac{\min(c,1)}{4}$.  
 Let $H$ be an admissible solution of \eqref{eq:jellium} such that  $H^\eta$ is of finite renormalized energy. Then, $H^\alpha$ is also of finite renormalized energy, and 
$$ \mathcal{W}^\alpha(\na H) = \mathcal{W}^\eta(\na H).$$
In particular, $H$ is of finite renormalized energy as soon as $H^\eta$ is for some $\eta \in (0, \frac{c}{4})$, and $ \mathcal{W}(\na H) =  \mathcal{W}^\eta(\na H)$ for all $\eta <  \frac{\min(c,1)}{4}$.   
\end{proposition}

\mspace
{\em Proof.} Let $R > 0$.  As $\Lambda$ is well-separated, 
\begin{equation} \label{cardinal_couronne} 
\left| \Lambda \cap (K_{R+2}\setminus K_{R-2}) \right| \le C R^{2}. 
\end{equation}
From this and the fact that the limit $\mathcal{W}^\eta(\na H)$ exists (in $\R$), it follows that 
\begin{equation} \label{cardinal_couronne_2}
\lim_{R \rightarrow +\infty} \frac{1}{R^3} \int_{K_{R+2}\setminus K_{R-2}} |\na H^\eta|^2 = 0. 
\end{equation}  
Let $\Omega_R$ be an open set such that $K_{R-1} \subset \Omega_R \subset K_R$ and such that 
\begin{align} \label{condition_distance}
\textrm{dist}\bigl( \pa \Omega_R \, , \,  \cup_{z \in \Lambda} B(-z,\eta)\bigr) \ge c' > 0 
\end{align}
 where $c'$ depends on $c$ only. This implies that  $G^{\eta}(\cdot + z)$, $G^\alpha(\cdot + z)$ are smooth at $\pa \Omega_R$ for all $z \in \Lambda$,  and that $H^{\eta}$, $H^\alpha$  are smooth at $\pa \Omega_R$. 
 
\mspace 
We now proceed as in the proof of Proposition \ref{prop:h_N_eta_h_N_alpha}. We write 
$$ H^\eta = H^\alpha + \sum_{z \in \Lambda} (G_S^\eta - G_S^\alpha)(\cdot + z),  $$
\begin{align*}
 \int_{\Omega_R} |\na H^\eta|^2 & = \int_{\Omega_R} |\na H^\alpha|^2 + 2 \sum_{z \in \Lambda}  \int_{\Omega_R} \na H^\alpha : \na  (G_S^\eta - G_S^\alpha)(\cdot + z)  \\ 
& + \sum_{z,z' \in \Lambda}  \int_{\Omega_R} \na  (G_S^\eta - G_S^\alpha)(\cdot + z)  : \na  (G_S^\eta - G_S^\alpha)(\cdot + z')  
\end{align*}
After integration by parts, and manipulations similar to those used to show Proposition \ref{prop:h_N_eta_h_N_alpha}, we end up with  
\begin{equation} \label{eq:Heta_Halpha}
 \int_{\Omega_R} |\na H^\eta|^2 -  \int_{\Omega_R} |\na H^\alpha|^2  =   \sum_{z \in \Lambda} \int_{\Omega_R}  (G_S^\eta - G_S^\alpha)(\cdot + z) d S^\alpha(\cdot + z) \\ 
 \end{equation}
Let us emphasize that the contribution of the boundary terms at $\pa \Omega_R$ is zero: indeed, thanks to \eqref{condition_distance}, $(G_S^\eta - G_S^\alpha)(\cdot + z)$ is zero at $\pa \Omega_R$ for any $z \in \Lambda$. Similarly, 
\begin{align*}
 \sum_{z \in \Lambda} \int_{\Omega_R}  (G_S^\eta - G_S^\alpha)(\cdot + z) d S^\alpha(\cdot + z) 
& = \sum_{z \in \Lambda \cap \Omega_R}  \int_{\Omega_R}  (G_S^\eta - G_S^\alpha)(\cdot + z) d S^\alpha(\cdot + z) \\
 & = \sum_{z \in \Lambda \cap \Omega_R}  \int_{\R^3} (G_S^\eta - G_S^\alpha)(\cdot + z) d S^\alpha(\cdot + z)
\end{align*}
The integral in the right-hand side was computed above, see \eqref{eq:S_alpha_diff_g_S} and the lines after: 
\begin{align*}
\sum_{z \in \Lambda \cap \Omega_R}  \int_{\R^3} (G_S^\eta - G_S^\alpha)(\cdot + z) d S^\alpha(\cdot + z)  = \left|\Lambda \cap \Omega_R\right|  \left( \frac{1}{\eta^3} -  \frac{1}{\alpha^3} \right) \left(  \int_{B_1} |\na G_S^1|^2 + \frac{3}{10\pi} |S|^2 \right)
\end{align*}
Back to \eqref{eq:Heta_Halpha}, we find 
\begin{align*}
 \int_{\Omega_R} |\na H^\eta|^2 -  \int_{\Omega_R} |\na H^\alpha|^2 =  \left|\Lambda \cap \Omega_R\right|  \left( \frac{1}{\eta^3} -  \frac{1}{\alpha^3} \right) \left(  \int_{B_1} |\na G_S^1|^2 + \frac{3}{10\pi} |S|^2 \right). 
 \end{align*}
 We  deduce from this identity, \eqref{cardinal_couronne} and  \eqref{cardinal_couronne_2} that 
 \begin{align*}
\lim_{R \rightarrow +\infty} \frac{1}{R^3}  \left(  \int_{\Omega_R} |\na H^\alpha|^2 -    \frac{\left|\Lambda \cap K_R\right|}{\alpha^3} \left(  \int_{B_1} |\na G_S^1|^2 + \frac{3}{10\pi} |S|^2 \right) \right)  = \mW^\eta(\na H), 
\end{align*}
and replacing $R$ by $R+1$: 
 \begin{align*}
\lim_{R \rightarrow +\infty} \frac{1}{R^3}  \left(  \int_{\Omega_{R+1}} |\na H^\alpha|^2 -    \frac{\left|\Lambda \cap K_R\right|}{\alpha^3} \left(  \int_{B_1} |\na G_S^1|^2 + \frac{3}{10\pi} |S|^2 \right) \right)  = \mW^\eta(\na H). 
\end{align*}
As $\Omega_R \subset K_R \subset \Omega_{R+1}$, the result follows.

\subsection{Resolution of the blown-up system for stationary point processes} \label{subsec:stationary}
As pointed out several times, we follow the strategy described in \cite{MR3309890} for the treatment of minimizers and  minima of Coulomb energies. But in our effective viscosity problem,  the points $x_{i,N}$ do not minimize the analogue $\mV_N$ of the Coulomb energy $\mathcal{H}_N$. Actually, although we consider the steady Stokes equation, our point distribution may be time dependent.  More precisely, in many settings, the dynamics of the suspension evolves on a timescale associated with viscous transport (scaling like $a^2$, with $a$ the radius of the particle), which is much smaller than the convective time scale (scaling like $a$). This allows to neglect the time derivative in the Stokes equation:  system \eqref{Sto}-\eqref{Sto2} corresponds then to a snapshot of the flow at a given time $t$. Even when one is interested in the long time behaviour, the existence of an equilibrium measure for the system of particles  is a very difficult problem. To bypass this issue, a usual point of view in the physics literature is to assume that the distribution of points is given by a stationary random process (whose refined description is an issue {\it per se}).  

\mspace
We will follow this point of view here, and  introduce a class of random point processes for which we can solve  \eqref{eq:jellium}. Let $X = \R$ or $X =\T_L := \R/(L\Z)$ for some 
$L > 0$. We denote by $Point_X$ the set of point distributions in $X^3$:  an element of $Point_X$ is a locally finite subset of $X^3$, in particular a  finite subset when $X=\T_L$. We endow $Point_X$ with the smallest $\sigma$-algebra $\mathcal{P}_X$  which makes measurable all the mappings 
$$Point_X \rightarrow \N, \quad  \omega \rightarrow |A \cap \omega|, \quad A \text{ borelian bounded subset of } X. $$ 
Given a probability space $(\Omega,\mathcal{A},P)$,  a random point process $\Lambda$ with values in $X^3$ is a measurable map from $\Omega$ to $Point_X$, see \cite{MR2371524}.  By pushing forward the probability $P$ with $\Lambda$, we can always assume that the process is in canonical form, that is  $\Omega = Point_X$, $\mathcal{A} =  \mathcal{P}_X$, and $\Lambda(\omega) = \omega$. 

\mspace
We shall consider processes that,  once in canonical form, are
\begin{description}
\item[(P1)] stationary:  the probability $P$ on $\Omega$ is invariant by the shifts    
$$ \tau_y : \Omega \rightarrow \Omega, \quad \omega \rightarrow  y + \omega, \quad y \in X^3. $$
\item[(P2)] ergodic: if $A \in \mathcal{A}$ satisfies  $\tau_y(A) = A$ for all $y$, then $P(A) = 0$ or $P(A) = 1$. 
\item[(P3)] uniformly well-separated: we mean that there exists $c > 0$ such that almost surely, $|z - z'| \ge c$ for all $z \neq z'$ in $\omega$. 
\end{description}
These properties are satisfied in two important contexts:
\begin{example}[Periodic point distributions] \label{example:periodic}
Namely, for $L > 0$, $a_1, \dots, a_M$ in $K_L$, we introduce  the set $\Lambda_0 := \{a_1, \dots, a_M\} + L \Z^d$. We can of course identify $\Lambda_0$ with a point distribution in $X^3$ with $X = \T_L$. We then take  $\Omega = \T_L^3$, $P$ the normalized Lebesgue measure on $\T^3_L$, and set $\Lambda(\omega) := \Lambda_0 + \omega$. It is easily checked that this random process satisfies all assumptions. Moreover, a realization of this process is  a translate of the initial periodic point distribution $\Lambda_0$. By translation,  the almost sure results  that we will show below (well-posedness of the blown-up system, convergence of $\mW_N$)   will actually yield  results for $\Lambda_0$ itself.    
\end{example}
\begin{example}[Poisson hard core processes] \label{example:hardcore} 
These processes are obtained from Poisson point processes, by removing balls in order to guarantee the hypothesis (P3). For instance, given $c > 0$, one can remove from the Poisson process all points $z$ which are not alone  in $B(z,c)$. This leads to the so-called {\em Mat\'ern I}  hard-core process. To increase the density of points while keeping (P3), one can refine the removal process in the following way: for each point $z$ of the Poisson process,  one associates an "age" $u_z$, with $(u_z)$ a family of i.i.d. variables, uniform over $(0,1)$. Then, one retains only the points $z$ that are (strictly) the "oldest" in $B(z,c)$.  This leads to the so-called {\em Mat\'ern II} hard-core process. Obviously, these two processes satisfy (P1) by stationarity of the Poisson process, and satisfy (P2) because they have only short range of correlations. For much more on hard core processes,  we refer to \cite{Bla}.  
\end{example}

\mspace
The point is now to solve almost surely the blown-up system \eqref{eq:jellium} for point processes with properties (P1)-(P2)-(P3). We first state 
\begin{proposition} \label{prop:existence:H_eta}
Let $\Lambda = \Lambda(\omega)$  a random point process with properties (P1)-(P2)-(P3). Let $\eta > 0$. For almost every $\omega$, there exists a  solution ${\bf H}^\eta(\omega, \cdot)$ of \eqref{eq:jellium_eta} in $H^1_{loc}(X^3)$ such that 
$$ \na {\bf H}^\eta(\omega,y) = D_{\bf H}^\eta(\tau_y\omega) $$
where   $D_{\bf H}^\eta \in L^2(\Omega)$ is the unique solution of the variational formulation \eqref{variational} below.     
\end{proposition}
\begin{remark}
In the case $X = \T_L$,  point distributions and solutions $H^\eta$ over $X^3$ can be  identified  with $L\Z^3$-periodic point distributions and $L\Z^3$-periodic solutions defined on $\R^3$. This identification is implicit here and in all that follows. 
\end{remark}

\mspace
{\em Proof.}  We treat the case $X=\R$,  the case $X=\T_L$ follows the same approach. We remind that the process is in canonical form:  $\Omega = Point_\R$, $\mathcal{A} = \mathcal{P}_\R$, $\Lambda(\omega) = \omega$. The idea is to associate to \eqref{eq:jellium_eta} a probabilistic variational formulation. This approach is inspired by works of Kozlov \cite{MR1329546,MR1369834}, see also \cite{MR2410410}. Prior to the statement of this variational formulation, we introduce some vocabulary and functional spaces.  First, for any $\R^d$-valued measurable $\phi = \phi(\omega)$, we call a realization of $\phi$ an application 
$$ R_\omega[\phi](y) := \phi(\tau_y\omega), \quad \omega \in \Omega.$$
For $p \in [1,+\infty)$, $\phi \in L^p(\Omega)$, as $\tau_y$ is measure preserving, we have for all $R > 0$ that  $\mathbb{E} \int_{K_R} |R_\omega[\phi]|^p = R^3 \, \mathbb{E} |\phi|^p$. Hence,  almost surely, $R_\omega[\phi]$ is in $L^p_{loc}(\R^3)$.  Also, for $\phi \in L^\infty(\Omega)$, one finds that almost surely $R_\omega[\phi] \in L^\infty_{loc}(\R^3)$. It is a consequence of Fatou's lemma: for all $R > 0$,
\begin{align*}
\mathbb{E} \|R_\omega[\phi]\|_{L^\infty(K_R)} & = \mathbb{E} \liminf_{p \rightarrow +\infty} \|R_\omega[\phi]\|_{L^p(K_R)} \le \liminf_{p \rightarrow +\infty}  \mathbb{E} \|R_\omega[\phi]\|_{L^p(K_R)} \\
& \le   \liminf_{p \rightarrow +\infty}   \left( \mathbb{E}   \|R_\omega[\phi]\|_{L^p(K_R)}^p \right)^{1/p} =     \liminf_{p \rightarrow +\infty}   \left( \mathbb{E}   |\phi|^p \right)^{1/p}  =  \|\phi\|_{L^\infty(\Omega)}.
\end{align*}
We say that $\phi$ is smooth if, almost surely, $R_\omega[\phi]$ is. For a smooth function $\phi$, we can define its stochastic gradient $\na_\omega \phi$ by the formula 
$$ \na_\omega \phi(\omega) := \na R_\omega[\phi]\vert_{y=0},    $$
where here and below, $\na = \na_y$ refers to the usual gradient (in space).
Note that $\na_\omega \phi(\tau_y \omega) = \na R_\omega[\phi](y)$.  One can define similarly the stochastic divergence, curl, {\it etc}, and reiterate to define partial stochastic derivatives $\pa^\alpha_\omega$. 

\mspace
Starting from a function $V \in L^p(\Omega)$, $p \in [1,+\infty]$ one can build smooth functions through convolution. Namely, for $\rho \in C^\infty_c(\R^3)$, one can define 
$$ \rho \star V(\omega) := \int_{\R^3} \rho(y)  V(\tau_y \omega) dy $$
which is easily seen to be  in $L^p(\Omega)$, as 
$$ \mathbb{E} |\rho \star V(\omega)|^p \le \mathbb{E}  \left( \int_{\R^3} |\rho(y)| dy \right)^{p-1} \left( \int_{\R^3} |\rho(y)| |V(\tau_y \omega) |^p dy \right) = \left( \int_{\R^3} |\rho(y)| dy \right)^{p}  \mathbb{E} |V(\omega)|^p $$
using that $\tau_y$ is measure-preserving. Moreover, it is smooth: we leave to the reader to check  
$$ R_\omega[\rho \star V] = \check{\rho} \star R_\omega[V], \quad \na_\omega (\rho \star V) = \na \check{\rho} \star V, \quad \check{\rho}(y) := \check{\rho}(-y). $$
 We are now ready to introduce the functional spaces we need. We set 
\begin{align*}
 \mathcal{D}_\sigma & := \{ \phi : \Omega \rightarrow \R^3 \text{ smooth, }  \pa^\alpha_\omega \phi \in L^2(\Omega) \; \forall \alpha, \; \na_\omega \cdot \phi = 0 \}, \\
\mathcal{V}_\sigma  & := \text{ the closure of } \{ \na_\omega \phi, \: \phi \in \mathcal{D}_\sigma \}  \text{ in } L^2(\Omega). 
\end{align*}
We remind that $S^\eta = \div \Psi^\eta$, with $\Psi^\eta$ defined in \eqref{def:Psi_eta}. We introduce
$$ \mathsf{\Pi}^\eta(\omega) := \sum_{z \in \omega} \Psi^\eta(z) $$
Note that it is well-defined, as $\Psi^\eta$ is supported in $B_\eta$ and $\omega$ is a discrete subset. It is measurable: indeed, $\Psi^\eta$ is the pointwise limit of a sequence of simple functions of the form $\sum_i \alpha_i 1_{A_i}$, where $A_i$ are Borel subsets of $\R^3$. As
 $$ \omega \rightarrow \sum_{z \in \omega} \sum_i \alpha_i 1_{A_i}(z) = \sum_i \alpha_i  |A_i \cap \omega| $$
is measurable by definition of the $\sigma$-algebra $\mathcal{A}$, we find that $\mathsf{\Pi}^\eta$ is. Moreover, as $\Lambda$ is uniformly well-separated,  one has 
$|\mathsf{\Pi}^\eta(\omega)| \le C \|\Psi^\eta\|_{L^\infty}$
 for a constant $C$ that does not depend on $\omega$, so that $\mathsf{\Pi}^\eta$ belongs to $L^\infty(\Omega)$. 
 
 \mspace
 We now introduce the variational formulation: {\em find  $D_{\bf H}^\eta \in \mathcal{V}_\sigma$ such that  for all $D_\phi \in \mathcal{V}_\sigma$,} 
\begin{equation} \label{variational}
 \mathbb{E} \, D_{\bf H}^\eta  : D_\phi  = - \mathbb{E} \,  \mathsf{\Pi}^\eta : D_\phi .
\end{equation}  
 As $\mathcal{V}_\sigma$ is a closed subspace of $L^2(\Omega)$, existence and uniqueness of a solution comes from the Riesz theorem. 
 
 \mspace
It remains to build a solution of \eqref{eq:jellium_eta} almost surely, based on $D_{\bf H}^\eta$. Let $\phi_k = \phi_k(\omega)$ a sequence in $\mathcal{D}_\sigma$ such that $\na_\omega \phi_k$ converges to $D_{\bf H}^\eta$ in $L^2(\Omega)$. Let $\rho \in C^\infty_c(\R^3)$. It is easily seen that  $\rho \star \phi_k$ also belongs to $\mathcal{D}_\sigma$ and that $\pa^\alpha_\omega \na_\omega (\rho \star \phi_k) = \pa^\alpha_\omega  (\rho \star  \na_\omega \phi_k)$ converges to the smooth function $\pa^\alpha_\omega (\rho \star D_{\bf H}^\eta)$ in $L^2(\Omega)$, for all $\alpha$. In particular, as $\na_\omega \times  \na_\omega (\rho \star \phi_k) = 0$, we find that  $\na_\omega \times (\rho \star D_{\bf H}^\eta) = 0$.  Applying the realization operator $R_\omega$, we deduce that 
$$  \na  \times (\check{\rho} \star R_\omega[D_{\bf H}^\eta]) = \check{\rho} \star \na \times  R_\omega[D_{\bf H}^\eta] = 0. $$
We recall that $R_\omega[D_{\bf H}^\eta]$ belongs almost surely to $L^2_{loc}(\R^3)$, so that  $\na \times  R_\omega[D_{\bf H}^\eta]$  is well-defined in $H^{-1}_{loc}(\R^3)$. Taking $\rho = \rho_n$ an approximation of the identity, and sending $n$ to infinity, we end up with 
$\na  \times R_\omega[D_{\bf H}^\eta]  = 0$ in $\R^3$. As curl-free vector fields on $\R^3$ are gradients, it follows that almost surely, there exists ${\bf H}^\eta = {\bf H}^\eta(\omega, y)$ with 
$$ \na {\bf H}^\eta(\omega,y)  = R_\omega[D_{\bf H}^\eta](y) = D_{\bf H}^\eta(\tau_y(\omega)), \quad \forall y \in \R^3. $$
In the case $X = \T_L$, one can show that the mean of $R_\omega[D_{\bf H}^\eta]$ is almost surely zero, so that the same result holds. 
Besides, because the matrices $\na_\omega \phi$, $\phi \in \mathcal{D}_\sigma$, have zero trace, the same holds for $D_{\bf H}^\eta$. Hence, 
$$ \div {\bf H}^\eta(\omega,y) = \text{trace}(\na {\bf H}^\eta(\omega,y)) =  \text{trace}(D_{\bf H}^\eta)(\tau_y(\omega)) = 0.  $$
One still has to prove that the first equation of \eqref{eq:jellium_eta} is satisfied. Therefore, we use \eqref{variational} with test function $D_\phi = \na_\omega  \phi$, where the smooth function $\phi$ is  of the form 
$$ \phi = \rho \star  ( \na_\omega \times \varphi), \: \varphi  : \Omega \rightarrow \R^3 \text{ a smooth function}. $$ 
Note that for smooth functions $\varphi, \tilde \varphi$, a stochastic integration by parts formula holds:  
\begin{align*}
\mathbb{E} \,  \pa_{\omega}^i \varphi \, \tilde \varphi & =  \mathbb{E}  \int_{K_1} \pa_{i} R_\omega[\varphi]  \, R_\omega[\tilde \varphi] =  - \mathbb{E}  \int_{K_1} R_\omega[\varphi]  \,  \pa_i R_\omega[\tilde \varphi] + \mathbb{E}   \int_{\pa K_1} n_i R_\omega[\varphi]  \, R_\omega[\tilde \varphi]   \\ 
& = -  \mathbb{E}  \int_{K_1} R_\omega[\varphi]  \,  \pa_i R_\omega[\tilde \varphi] = - \mathbb{E} \,  \varphi \, \pa_{\omega,i} \tilde \varphi. 
\end{align*}
Thanks to this formula, we may write 
\begin{align*}
 \mathbb{E} \, D_{\bf H}^\eta   :   \na_\omega (\rho \star  (\na_\omega \times \varphi))  & = \mathbb{E} \, \check{\rho} \star D_{\bf H}^\eta   :   \na_\omega ( \na_\omega \times \varphi) \\
 & =  - \mathbb{E} \, \na_\omega \times  (\na_\omega \cdot  (\check{\rho} \star  D_{\bf H}^\eta))  \cdot  \varphi. 
 \end{align*}
 Similarly, we find 
 $$ - \mathbb{E} \,  \mathsf{\Pi}^\eta:  \na_\omega (\rho \star  \na_\omega \times \varphi)  = \mathbb{E}  \na_\omega \times  (\na_\omega \cdot  (\check{\rho} \star   \mathsf{\Pi}^\eta)) \cdot \varphi. $$
As this identity is valid for all smooth test fields $\varphi$, we end up with 
$$ -  \na_\omega \times  (\na_\omega \cdot  (\check{\rho} \star D_{\bf H}^\eta)) =  \na_\omega \times  (\na_\omega \cdot  (\check{\rho} \star  \mathsf{\Pi}^\eta)). $$
Proceeding as above, we find that almost surely
$$ - \na \times  \div R_\omega[D_{\bf H}^\eta]  =  \na \times  \div R_\omega[\mathsf{\Pi}^\eta] $$
which can be written 
$$ \na \times (-\Delta {\bf H}^\eta) = \na \times  \div \sum_{z \in \Omega} \Psi^\eta(\cdot +z). $$
It follows that there exists ${\bf P}^\eta = {\bf P}^\eta(\omega,y)$, such that 
$$ - \Delta {\bf H}^\eta + \na {\bf P}^\eta =   \div \sum_{z \in \omega} \Psi^\eta(\cdot +z) = \sum_{z \in \omega} S^\eta(\cdot +z) $$
which concludes the proof of the proposition. 
\begin{corollary} \label{cor:existence:H}
For random point processes with properties (P1)-(P2)-(P3),  there exists almost surely a solution $H$ of \eqref{eq:jellium} with finite renormalized energy and such that for all $\eta > 0$, the gradient field $\na H^\eta$, where $H^\eta$ is  given by \eqref{def:H_H_eta}, coincides with the gradient field $\na {\bf H}^\eta$ of Proposition \ref{prop:existence:H_eta}.  Moreover, 
$$  \mathcal{W}(\na H) =  - \lim_{\eta \rightarrow 0} \left(  \mathbb{E} \int_{K_1}|\na H^\eta|^2  - \frac{m}{\eta^3}  \Bigl( \int_{B^1} |\na G_S^1|^2 + \frac{3}{10\pi} |S|^2 \Bigr)\right) $$
where $m := \E |\Lambda \cap K_1|$ is the mean intensity of the point process, the expression at the right-hand side being actually constant for $\eta$ small enough. 
\end{corollary}

\mspace
{\em Proof}. By the definition of the mean intensity and by property (P2), which allows to apply the ergodic theorem ({\it cf.} \cite[Corollary 12.2.V]{MR2371524}), we have almost surely
\begin{equation} \label{ergodic1}
\lim_{R \rightarrow \infty} \frac{|\Lambda \cap K_R|}{R^3} = m.
\end{equation}
Let $\eta_0 < \frac{\min(c,1)}{4}$ fixed, and ${\bf H}^{\eta_0}$ given by the previous proposition. We set 
\begin{equation}
H(\omega,y) := {\bf H}^{\eta_0}(\omega,y) + \sum_{z \in \omega} (G_S- G_S^{\eta_0})(y+z). 
\end{equation}
It is clearly an admissible solution of \eqref{eq:jellium}.  By Proposition \ref{prop:stationarity:Weta}, in order to show that $H$ has almost surely finite renormalized energy, it is enough to show that  for one $\eta <  \frac{\min(c,1)}{4}$, almost surely, the function $H^\eta$ given by \eqref{def:H_H_eta}, namely 
\begin{align*} 
H^\eta(\omega,y) & := H(\omega,y) + \sum_{z \in \omega} (G^\eta_S - G_S)(y+z) \\
& = {\bf H}^{\eta_0}(\omega,y) + \sum_{z \in \omega} (G^\eta_S- G_S^{\eta_0})(y+z) 
\end{align*}
has finite renormalized energy. This holds for $\eta = \eta_0$, as $H^{\eta_0} = {\bf H}^{\eta_0}$ and the ergodic theorem applies. We then notice that   
\begin{equation} \label{station_H_eta}
\na H^\eta(\omega,y) = D_H^\eta(\tau_y(\omega)), \quad   D_H^\eta(\omega) :=   D_{\bf H}^{\eta_0}(\omega) + \sum_{z \in \omega}  \na(G^\eta_S- G_S^{\eta_0})(z). 
\end{equation}
We remark that $G^\eta_S- G_S^{\eta_0} = 0$ outside $B_{\max(\eta,\eta_0)}$, so that the sum at the r.h.s. has only a finite number of non-zero terms.  In the same way as we proved  that the function $\mathsf{\Pi}^\eta$ belongs to $L^\infty(\Omega)$, we get that $\sum_{z \in \omega}  \na (G^\eta_S- G_S^{\eta_0})(z)$ defines an element of $L^\infty(\Omega)$. Hence, by the ergodic theorem, we have  almost surely 
 $$ \lim_{R \rightarrow +\infty} \frac{1}{R^3} \int_{K_R} |\na H^\eta|^2 \: \rightarrow \mathbb{E}  \int_{K_1} |\na H^\eta|^2. $$
 Combining this with \eqref{ergodic1} and Proposition \ref{prop:stationarity:Weta}, we obtain the formula for $\mW(\na H)$. 
 
\mspace 
  The last step is to prove that for all $\eta > 0$, $\na H^\eta = \na {\bf H}^\eta$ almost surely. As a consequence of  the ergodic theorem, one has almost surely  
 $$  \limsup_{R \rightarrow +\infty} \frac{1}{R^3} \int_{K_R} |\na H^\eta|^2 < +\infty, \quad   \limsup_{R \rightarrow +\infty} \frac{1}{R^3} \int_{K_R} |\na {\bf H}^\eta|^2 < +\infty $$
 Reasoning as in the proof of Proposition \ref{jellium_uniqueness},  we find that their gradients  differ by a constant: 
$$ \na H^\eta(\omega,y) =  \na {\bf H}^\eta(\omega,y) + C(\omega). $$
Applying again the ergodic theorem, we get that almost surely
$ \E D_{H}^\eta  = \E D_{\bf H}^\eta + C(\omega)$. As $D_{\bf H}^\eta$ belongs to $\mathcal{V}_\sigma$, its expectation is easily seen to be zero. To conclude, it remains to prove that 
$\E D_{H}^\eta = \E \sum_{z \in \omega} \na (G_S^\eta - G_S^{\eta_0})(z)$ is zero. Using stationarity, we write, for all $R > 0$, 
 $$ \E \sum_{z \in \omega} \na (G_S^\eta - G_S^{\eta_0})(z) = \frac{1}{R^3} \E \sum_{z \in \omega}  \int_{K_R}  \na (G_S^\eta - G_S^{\eta_0})(z+y) dy. $$   
 We remark that for all $z$ outside a $\max(\eta, \eta_0)$-neighborhood of  $\pa K_R$,  $ \int_{K_R}  \na (G_S^\eta - G_S^{\eta_0})(z+\cdot) = \int_{\pa K_R} n \otimes  (G_S^\eta - G_S^{\eta_0})(z+\cdot) = 0$. It follows from the separation assumption and the $L^\infty$ bound on   $\na (G_S^\eta - G_S^{\eta_0})$ that 
 $$  \frac{1}{R^3} \E \sum_{z \in \omega}  \int_{K_R}  \na (G_S^\eta - G_S^{\eta_0})(z+y) dy  = O(1/R) \rightarrow 0 \quad \text{as } \: R \rightarrow +\infty. $$

\section{Convergence of $\mV_N$} \label{sec4} 
This section concludes our analysis of the quadratic correction to the effective viscosity. From Theorem \ref{thm1}, we know that this quadratic correction should be given by the limit of $\mV_N$ as $N$ goes to infinity, where $\mV_N$ was introduced in \eqref{def:VN}. We show here that the functional $\mV_N$ has indeed a limit, when the particles are given by the kind of stationary point processes seen in Section \ref{sec3}. 

\subsection{Proof of convergence} \label{subsec:conv_VN}
Let $\eps > 0$ a small parameter, and $\Lambda = \Lambda(\omega)$ a random point process with properties (P1)-(P2)-(P3): stationarity, ergodicity, and uniform separation. As seen in Examples \ref{example:periodic} and \ref{example:hardcore}, this setting covers the case of {\em periodic patterns of points} as well as {\em classical hard core processes}. We set 
$N = N(\eps)$ the cardinal of the set 
$$ \{ x \in \eps \check{\Lambda}, \: B(x,\eps) \subset \mO  \} = \{ x_{1,N}, \dots, x_{N,N} \} $$ 
where $\check{\Lambda} := -\Lambda$ and where we label the elements arbitrarily.  Note that $N$ depends on $\omega$, although it does not appear explicitly.  From the fact that $\Lambda$ is uniformly well-separated and from the ergodic theorem ({\it cf.} \cite[Corollary 12.2.V]{MR2371524}), we can deduce that almost surely,
\begin{equation} \label{lim_N_eps:random} 
\lim_{\eps \rightarrow 0} N(\eps) \eps^3 = \lim_{\eps \rightarrow 0}  |\eps \check{\Lambda}(\omega) \cap \mO| \, \eps^3  = \lim_{\eps \rightarrow 0}  \frac{ |\check{\Lambda}(\omega) \cap \eps^{-1}\mO|}{ \eps^{-3} |\mO|}  |\mO|    = m |\mO| 
\end{equation}
so that we shall note indifferently $\lim_{\eps \rightarrow 0}$ or  $\lim_{N \rightarrow +\infty}$. Note that, strictly speaking, $N = N(\eps)$ does not necessarily cover all integer values when $\eps \rightarrow 0$, but this is no difficulty.   

\mspace
More generally, for all $\varphi$ smooth and compactly supported in $\R^3$, ergodicity implies
$$ \lim_{N \rightarrow +\infty}  \frac{1}{N} \sum_{i=1}^N \varphi(x_{i}) =  \lim_{N \rightarrow +\infty}  \frac{1}{N} \sum_{x_i \in  \mO}  \varphi(x_{i})  =   \lim_{N \rightarrow +\infty}  \frac{1}{\eps^3 N} m  \int_{\mO} \varphi(x) dx = \frac{1}{|\mO|}  \int_{\mO} \varphi(x) dx $$
which shows that \eqref{H1} is satisfied with $f = \frac{1}{|\mO|} 1_{\mO}$. The hypothesis \eqref{H2} is also trivially satisfied, as well as \eqref{hyp:in_O}.  Our main theorem is 
\begin{theorem} \label{thm2} Almost surely, 
$$\lim_{N \rightarrow +\infty} \mV_N   =   \frac{25}{2 m^2} \mathcal{W}(\na H)        $$
with $m$  the mean intensity of the process, and  $H$  the solution of \eqref{eq:jellium} given in Corollary \ref{cor:existence:H}. 
\end{theorem}
\noindent
The rest of the paragraph is dedicated to the proof of this theorem. 

\mspace
Let $\eta$ satisfying $\eta  <  \frac{\min(c,1)}{4}$  and  $\eta < \frac{c}{2} (m|\mO|)^{-1/3}$. By  \eqref{lim_N_eps:random}, it follows that almost surely, for $\eps$ small enough, $\eps \eta < \frac{c}{2} N^{-1/3}$. By Corollary \ref{cor:asymptote:WN},
\begin{equation} \label{conv_VN}
 \lim_{N  \rightarrow +\infty} \quad   \mV_N  + \frac{25 |\mO|}{2 N^2} \Bigl(  \int_{\R^3} | \na h_{N}^{\eta \eps} |^2 - \frac{N}{(\eta \eps)^3} \big( \int_{B^1} |\na G_S^1|^2 + \frac{3}{10\pi} |S|^2 \big) \Bigr) = 0.  
 \end{equation}
We denote $h_\eps^\eta := h_N^{\eta \eps}$, see \eqref{def:hNeta}-\eqref{eq:hNeta}. Let $H$ be the solution  of the blown-up system \eqref{eq:jellium} provided by Corollary \ref{cor:existence:H}, $H^\eta$ given in \eqref{def:H_H_eta}, and $P^\eta$ as in \eqref{eq:jellium_eta}.  We define  new fields $\bar{h}^\eta_\eps, \bar{p}^\eta_\eps$  by the following conditions: $\bar{h}^\eta_\eps \in \dot{H}^1(\R^3)$, 
\begin{align*}
& \bar{h}^\eta_\eps(\omega,x)  =  \frac{1}{\eps^2} H^\eta\big(\frac{x}{\eps}\big) - \dashint_{\mO}  \frac{1}{\eps^2} H^\eta\big(\frac{\cdot}{\eps}\big), \: x \in \mO\\
&  \overline{p}_\eps^\eta(\omega,x)  =  \frac{1}{\eps^3} P^\eta\big(\frac{x}{\eps}\big) - \dashint_{\mO}  \frac{1}{\eps^3} P^\eta\big(\frac{\cdot}{\eps}\big), \: x \in \mO \\
& - \Delta \bar{h}^\eta_\eps + \na \bar{p}^\eta_\eps = 0, \quad \div \bar{h}^\eta_\eps = 0 \: \text{ in } \: \ext \mO . 
\end{align*}
We omit to indicate the dependence in $\omega$ to lighten notations. We claim: 
\begin{proposition} \label{prop:barh}
\begin{align*}
& \lim_{\eps \rightarrow 0} -\frac{1}{N^2} \bigl( \int_{\R^3} |\na \bar{h}^\eta_\eps|^2 - \frac{N}{(\eta \eps)^3}   ( \int_{B^1} |\na G_S^1|^2 + \frac{3}{10\pi} |S|^2 ) \bigr) = \frac{1}{m^2|\mathcal{O}|} \mathcal{W}^\eta(\na H) 
\end{align*}
\end{proposition}
\begin{proposition} \label{prop:barh_h}
$$ \lim_{\eps \rightarrow 0} \eps^6 \int_{\R^3} |\na (h^\eta_\eps - \bar{h}^\eta_\eps)|^2 = 0. $$ 
\end{proposition}
\noindent
Note that,  by Proposition \ref{prop:stationarity:Weta} and our choice of $\eta$,   $\mathcal{W}^\eta(\na H)$ =  $\mathcal{W}(\na H)$. Theorem \ref{thm2} follows directly from this fact,  \eqref{conv_VN} and the propositions. 

\mspace
{\em Proof of Proposition \ref{prop:barh}}. We know from Corollary \ref{cor:existence:H} that 
 $$ \mathcal{W}^\eta(\na H) = - \Big( \mathbb{E} \int_{K_1}|\na H^\eta|^2  - \frac{m}{\eta^3}  \bigl( \int_{B^1} |\na G_S^1|^2 + \frac{3}{10\pi} |S|^2 \bigr)\Big) 
  $$
From this and relation \eqref{lim_N_eps:random}, we see that the proposition amounts to the statement 
$$ \lim_{\eps \rightarrow 0} \frac{\eps^6}{|\mO|}\int_{\R^3} |\na \bar{h}^\eta_\eps|^2 =  \E \int_{K_1} |\na H^\eta|^2.  $$
A simple application of the ergodic theorem shows that almost surely 
$$ \frac{\eps^6}{|\mO|}\int_{\mO} |\na \bar{h}^\eta_\eps|^2 = \frac{1}{|\mO|}\ \int_{\mO} |\na_y H^\eta\big(\frac{x}{\eps}\big)|^2 dy \rightarrow  \E \int_{K_1} |\na H^\eta|^2. $$
It remains to show that 
\begin{equation} \label{residual_h_bar} 
 \lim_{\eps \rightarrow 0} \:  \eps^6 \int_{\ext \mO} |\na \bar{h}^\eta_\eps|^2 = 0. 
 \end{equation}
It will be deduced from the well-known fact that the Stokes solution $\bar{h}^\eta_\eps$ minimizes 
$$\displaystyle \int_{\ext \mO} |\na \bar{h}|^2$$
 among divergence-free fields $\bar{h}$ in $\ext \mO$ satisfying the Dirichlet condition $\bar{h}\vert_{\pa \mO} = \bar{h}^\eta_\eps\vert_{\pa \mO}$.  
 
 \mspace
 First, we prove that the $H^{1/2}(\partial \mO)$-norm of $\eps^3 \bar{h}^\eta_\eps$ goes to zero.  In this perspective, we introduce for all $\delta > 0$  a function $\chi_\delta$ with $\chi_\delta = 1$ in  a $\frac{\delta}{2}$-neighborhood of $\pa \mO$, $\chi_\delta = 0$ outside a $\delta$-neighborhood of $\pa \mO$.  We write 
\begin{align*}
 \|\eps^3 \bar{h}^\eta_\eps\|_{H^{1/2}(\pa \mO)} & = \|\eps^3 \bar{h}^\eta_\eps \chi_\delta\|_{H^{1/2}(\pa \mO)} \\ 
 & \le C \left( \|\eps^3 \bar{h}^\eta_\eps \chi_\delta\|_{L^2(\mO)} + \|\eps^3 \na \bar{h}^\eta_\eps \chi_\delta\|_{L^2(\mO)}  + \|\eps^3 \bar{h}^\eta_\eps \na \chi_\delta\|_{L^2(\mO)} \right). 
 \end{align*}
By the ergodic theorem and Corollary \ref{cor:existence:H},  $\eps^3 \na \bar{h}^\eta_\eps = \na_y H^\eta(\frac{\cdot}{\eps})$ converges almost surely  weakly in $L^2(\mO)$ to 
$\E D_{\mathbf{H}^{\eta}} = 0$. Let $\varphi \in L^2(\mO)$. By standard results on the divergence operator, {\it cf} \cite{MR1284205}, there exists $v \in H^1_0(\mO)$ with $\div v = \varphi - \dashint_{\mO} \varphi$, $\|v\|_{H^1(\mO)} \le C_\mO \|\varphi\|_{L^2(\mO)}$.  As by definition  $\bar{h}^\eta_\eps$ has zero mean over $\mO$, it follows that 
$$ \int_\mO  \eps^3 \bar{h}^\eta_\eps  \varphi =  \int_\mO   \eps^3 \bar{h}^\eta_\eps \,  (\varphi - \dashint_{\mO} \varphi) = -  \int_\mO  \eps^3 \na  \bar{h}^\eta_\eps \, v \rightarrow 0 \: \text{ as } \eps \rightarrow 0.$$ 
Hence,  $\eps^3 \bar{h}^\eta_\eps$ converges weakly to zero in $H^1(\mO)$ and therefore strongly in $L^2(\mO)$. It follows that for any given $\delta$, 
$$ \|\eps^3 \bar{h}^\eta_\eps \chi_\delta\|_{L^2(\mO)} \rightarrow 0, \quad  \|\eps^3 \bar{h}^\eta_\eps \na \chi_\delta\|_{L^2(\mO)}  \rightarrow 0 \quad \text{as } \eps \rightarrow 0. $$
To conclude, it is enough to show that  $\limsup_{\eps \rightarrow 0} \|\eps^3 \na \bar{h}^\eta_\eps \chi_\delta\|_{L^2(\mO)}$ goes to zero as $\delta \rightarrow 0$. This comes from
\begin{equation} \label{eq:delta_bound}
 \|\eps^3 \na \bar{h}^\eta_\eps \chi_\delta\|_{L^2(\mO)}^2 = \int_\mO |\na H^\eta(\cdot/\eps)|^2 \, \chi_\delta^2 \: \xrightarrow[\eps \rightarrow 0]{} \:  \E |D_{\bf H}^\eta|^2  \int_\mO \chi_\delta^2 \le C \delta.  
\end{equation}
Finally, $\|\eps^3 \bar{h}^\eta_\eps\|_{H^{1/2}(\pa \mO)} = o(1)$. To conclude that \eqref{residual_h_bar} holds, we notice that 
$$\int_{\pa \mO}  \bar{h}^\eta_\eps \cdot n =  \int_{\mO} \div  \bar{h}^\eta_\eps = 0.$$
By classical results on the right inverse of the divergence operator, see \cite{MR1284205}, one can find for $R$ such that $\mO \Subset B(0,R)$  a solution $\bar{h}$ of the equation 
$$ \div \bar{h} = 0 \quad \text{ in } \, \ext \mO \cap B(0,R), \quad \bar{h}\vert_{\pa \mO} =  \bar{h}^\eta_\eps\vert_{\pa \mO}, \quad \bar{h}\vert_{\pa B(0,R)} = 0, $$
and such that 
$$\| \bar{h} \|_{H^1(\ext \mO \cap B(0,R))} \le C \|  \bar{h}^\eta_\eps \|_{H^{1/2}(\pa \mO)} = o(\eps^{-3}).$$
Extending $\bar{h}$ by zero outside $B(0,R)$, we find 
\begin{equation} \label{bound_barh}
\int_{\ext \mO} |\na \bar{h}^\eta_\eps|^2 \le  \int_{\ext \mO} |\na \bar{h}|^2 = o(\eps^{-6}). 
\end{equation}
This concludes the proof of the proposition. 

\mspace
{\em Proof of Proposition \ref{prop:barh_h}}. Let $h := h_\eps^\eta - \bar{h}_\eps^\eta$.  It satisfies an equation of the form
$$ -\Delta h + \na p = R_1 + R_2 + R_3, \quad \div h = 0 \text{ in } \:  \R^3, $$
where the various source terms will now be defined. First,  
$$ R_1 :=  \sigma(\bar{h}_\eps^\eta, \bar{p}^\eta_\eps)n\vert_{\pa (\ext \mO)} \, s_{\pa}. $$
Here,  the value of the stress is  taken from $\ext \mO$,  $n$ refers to the normal vector pointing outward $\mO$ and $s_{\pa}$ refers to the surface measure on 
$\pa \mO$. We remind that $\bar{h}^\eta_\eps \in \dot{H}^1(\R^3)$ does not jump at the boundary, but its derivatives do, so that one must specify from which side the stress is considered. Then, 
$$ R_2 :=   -\sigma(\bar{h}_\eps^\eta, \bar{p}^\eta_\eps)n\vert_{\pa \mO}  \, s_{\pa} = - \frac{1}{\eps^3} \sigma\Big(H^\eta, P^\eta - \dashint_\mO  P^\eta(\omega,\cdot/\eps)\Big)\bigl(\frac{\cdot}{\eps}\bigr)\vert_{\pa \mO}n \, s_{\pa} $$
with the value of the stress taken from $\mO$, and $n$ as before.  Noticing that $ S \na f = - \frac{1}{|\mO|} Sn \, s_\pa$, we finally  set
$$ R_3 := -{\bf 1}_{\mO} \, \sum_{i \in I^\eta_\eps} S^{\eta \eps}(x-x_i) + \frac{N}{|\mO|} Sn \, s_{\pa}. $$
where 
$$I^\eta_\eps = \{ i, B(x_i, \eps) \not \subset \mO, \: B(x_i, \eta \eps) \cap \mO \neq \emptyset\}.$$ 
Note that the term $R_3$ is supported in pieces of spheres. From \eqref{def:Seta}, we know that for all $\eta > 0$
\begin{equation} \label{zeroint_Seta}
 \int_{\R^3} S^\eta =  \int_{\R^3} \left( -\Delta G_S^\eta + \na p_S^\eta\right)  = 0.
 \end{equation} 
This allows to show that the integral of $R_2+R_3$ is zero. Indeed,
\begin{align*}
 \int_{\R_3} R_2  & = \frac{1}{\eps^4} \int_{\mO} (-\Delta H^\eta + \na P^\eta)(\cdot/\eps)   =  \int_{\mO}  \sum_{i, B(x_i, \eta \eps) \cap \mO \neq \emptyset}  \hspace{-0.6cm}S^{\eta \eps}(\cdot-x_i) \: = \sum_{i \in I_\eps^\eta}    \int_{\mO}  S^{\eta \eps}(\cdot-x_i) 
\end{align*}
so that 
\begin{equation} \label{zeroint_R2_R3}
\int_{\R^3} (R_2 + R_3)  =  \frac{N}{|\mO|} \int_{\pa \mO} S n \, ds_\pa = 0. 
\end{equation}
The point is now to prove that $\eps^3 \|\na h\|_{L^2(\R^3)} \rightarrow 0$ as $\eps \rightarrow 0$. From a simple energy estimate, and taking \eqref{zeroint_R2_R3} into account, we find 
\begin{equation} \label{estimate_h}
 \|\na h\|_{L^2(\R^3)}^2  =  \langle  R_1 , h \rangle  + \langle R_2, h - \dashint_{\mO} h \rangle +  \langle R_3, h - \dashint_{\mO} h \rangle. 
\end{equation}
As $(\bar{h}_\eps^\eta, \bar{p}_\eps^\eta)$ is a solution of a homogeneous Stokes equation in $\ext \mO$, we get from an integration by parts: 
\begin{equation} \label{estimate_R1} 
 \langle  R_1 , h \rangle =  \int_{\ext \mO} \na \bar{h}_\eps^\eta \cdot \na h \le \nu(\eps)\eps^{-3} \| \na h \|_{L^2(\R^3)}, \quad \nu(\eps) \rightarrow 0 \quad \text{as } \eps \rightarrow 0, 
 \end{equation}
using  the Cauchy-Schwarz inequality and the bound \eqref{bound_barh}.

\mspace
We now wish to show that 
\begin{equation} \label{estimate_R2R3} 
 \langle  (R_2 + R_3) , h  - \dashint_{\mO} h \rangle \le \nu(\eps) \eps^{-3} \| \na h \|_{L^2(\R^3)},  
\end{equation}
for some $\nu(\eps)$ going to zero with $\eps$.  More precisely, we will prove that for any divergence-free $\varphi \in \dot{H}^1(\R^3)$, 
\begin{equation} \label{estimate_R2R3_varphi} 
 \langle  (R_2 + R_3) , \varphi \rangle \le \nu(\eps)  \eps^{-3}  (\| \na \varphi \|_{L^2(\R^3)} + \| \varphi \|_{H^1(\mO)}),  \quad \nu(\eps) \rightarrow 0 \text{ as } \eps \rightarrow 0.
\end{equation}
which implies \eqref{estimate_R2R3} by Poincar\'e inequality. We first notice that    
\begin{equation} \label{R2_varphi}
\langle  R_2 , \varphi \rangle = \frac{1}{\eps^3}   \,  \langle \, n \cdot F_2^\eps , \varphi \,  \rangle_{\langle H^{-1/2}(\pa \mO), H^{1/2}(\pa \mO)\rangle}
\end{equation}
where 
\begin{equation} \label{def_F2eps}
F_2^\eps := \eps^3 \left( 2  D(\bar{h}^\eta_\eps) -  \overline{p}_\eps^\eta \Id \right) =  2  D(H^\eta)(\omega, \cdot/\eps) +  \left(P^\eta(\omega,\cdot/\eps) - \dashint_\mO  P^\eta(\omega,\cdot/\eps)\right) \Id. 
\end{equation}
 Then, we use the relation $S^\eta = \div \Psi^\eta$, {\it cf.}  Lemma \ref{lem:smoothing} and integrate by parts to get 
\begin{align*}
\langle R_3, \varphi \rangle & =  \frac{1}{\eps^3}  \sum_{i \in I^\eps_\eta} \Biggl(  -\int_{\pa \mO} n \cdot \Psi^\eta \left(\frac{x-x_i}{\eps}\right) \cdot \varphi(x) ds_\pa(x) \\ 
& +  \int_{\mO}  \Psi^\eta \left(\frac{x-x_i}{\eps}\right) : \na \varphi(x) dx \Biggr)  +  \frac{N}{|\mO|}\int_{\pa \mO} Sn(x) \cdot \varphi(x) \, ds_\pa(x)  
\end{align*}
For a fixed $\eta$, there is a constant $C$ (depending on $\eta$) such that 
\begin{align*}
& \sum_{i \in I^\eps_\eta}  \int_{\mO} \Psi^\eta \left(\frac{x-x_i}{\eps}\right) : \na \varphi(x) dx \\
&  \le  C \sum_{i \in I^\eps_\eta}    \int_{B(x_i, \eta \eps) \cap \mO} |\na \varphi|(x) dx  \le C | \cup_{i \in I_\eps^\eta} B(x_i, \eta \eps)|^{1/2} \|\na \varphi\|_{L^2(\R^3)} \le C \eps^{1/2} \|\na \varphi\|_{L^2(\R^3)}  
\end{align*}
For the last inequality, we have used that all $x_i$'s with $i \in I_\eps^\eta$ belong to an $\eps$-neighborhood of $\pa \mO$, so that $|I_\eps^\eta| =  O(\eps^{-2})$.   Hence, 
\begin{equation} \label{ineq_R3}
\begin{aligned}
\langle R_3, \varphi \rangle & \le  \frac{1}{\eps^3} \sum_{i \in I^\eps_\eta}  - \int_{\pa \mO} n \cdot \Psi^\eta \left(\frac{x-x_i}{\eps}\right) \cdot \varphi(x) ds_\pa(x) \\
& +  \frac{N}{|\mO|}\int_{\pa \mO} Sn(x) \cdot \varphi(x) \, ds_\pa(x) + \nu(\eps) \eps^{-5/2} \|\na \varphi\|_{L^2(\mO)}. 
\end{aligned}
\end{equation} 
Let
\begin{equation}  \label{def_F3eps}
F_3(\omega) :=   -\sum_{z \in \Lambda} \Psi^\eta(z) + m S, \quad F_3^\eps(x) := F_3(\tau_{x/\eps}(\omega))
\end{equation}
 We claim  that $\E \int_{K_1} F_3 = 0$. Indeed, by stationarity, for all $R > 0$ 
\begin{align*}
\E  \sum_{z \in \Lambda}  \Psi^\eta(z) & = \frac{1}{R^3} \E  \sum_{z \in \Lambda} \int_{K_R} \Psi^\eta(y+z) dy   \\
& = \frac{1}{R^3} \E  \hspace{-0.3cm} \sum_{\substack{z \in \Lambda, \\ K_R \supset B(-z,\eta)}} \hspace{-0.3cm} \int_{K_R} \Psi^\eta(y+z) dy + \frac{1}{R^3} \E \hspace{-0.3cm}  \sum_{\substack{z \in \Lambda, \\ \pa K_R \cap B(-z,\eta) \neq \emptyset}} \hspace{-0.3cm} \int_{K_R} \Psi^\eta(y+z) dy \\
& + \frac{1}{R^3} \E \hspace{-0.3cm} \sum_{\substack{z \in \Lambda, \\ K_R \cap B(-z,\eta)  = \emptyset}} \hspace{-0.3cm}  \int_{K_R} \Psi^\eta(y+z) dy  \\
& = \frac{1}{R^3} \E  \hspace{-0.3cm}  \sum_{\substack{z \in \Lambda, \\ K_R \supset B(-z,\eta)}} \int_{K_R}  \Psi^\eta(y+z) dy + \frac{1}{R^3} \E     \hspace{-0.3cm}  \sum_{\substack{z \in \Lambda, \\ \pa K_R \cap B(-z,\eta)  \neq \emptyset}}  \hspace{-0.3cm}   \int_{K_R} \Psi^\eta(y+z) dy  \\
& = \frac{1}{R^3} \E \big|\{z, K_R \supset B(0,\eta) - z\}\big| \int_{B(0,\eta)} \Psi^\eta(y) dy \:  +  \:  O\big(\frac{1}{R}\big), \quad R \gg 1.
\end{align*}
We have used crucially the fact that $\Psi^\eta$ is supported in $B(0,\eta)$. The $O(\frac{1}{R})$-term is associated to the points  $z \in \Lambda$ which lie in a $\delta$-neighborhood of $\pa K_R$: see the end of the proof of Corollary \ref{cor:existence:H} for similar reasoning. By sending $R$ to infinity, we find that almost surely 
$$ \E F_3 = - m \int_{B(0,\eta)} \Psi^\eta \: + m S. $$
The last step is to compute  $\int_{B(0,\eta)} \Psi^\eta$, which is independent of $\eta$ by homogeneity. It is in particular equal to $\lim_{\eta \rightarrow 0} \langle \Psi^\eta, 1 \rangle$, a limit that was already computed in the proof of Lemma \ref{lem:smoothing}, {\it cf.} \eqref{lim_Psi1}-\eqref{lim_Psi2}. We get $\int_{B(0,\eta)} \Psi^\eta = S$, which shows that 
$\E F_3 = 0$.

\mspace    
By the definition of $F^\eps_3$, we can write 
\begin{align*}
& \frac{1}{\eps^3}\sum_{i \in I^\eps_\eta}  \int_{\pa \mO \cap B(x_i, \eta \eps)} \hspace{-0.2cm} - n \cdot \Psi^\eta \left(\frac{x-x_i}{\eps}\right) \cdot \varphi(x) ds_\pa(x)  \: + \:  \frac{N}{|\mO|}\int_{\pa \mO} Sn(x) \cdot \varphi(x) \, ds_\pa(x)  \\ 
 &= \frac{1}{\eps^3}\int_{\pa \mO} n(x) \cdot F^\eps_3(x) \cdot \varphi(x) ds_\pa(x) + \left( \frac{N}{|\mO|} - \frac{m}{\eps^3} \right) \int_{\pa \mO} Sn \cdot \varphi \\
 & \le  \frac{1}{\eps^3}\int_{\pa \mO} n(x) \cdot F^\eps_3(x) \cdot \varphi(x) ds_\pa(x)  + \nu(\eps) \eps^{-3}\|\varphi\|_{H^1(\mO)}, \quad \nu(\eps) \xrightarrow[\eps \rightarrow 0]{} 0
\end{align*} 
where the last inequality follows from \eqref{lim_N_eps:random}. Plugging this inequality in \eqref{ineq_R3}, and combining with \eqref{R2_varphi}, we see that to derive \eqref{estimate_R2R3_varphi}, it remains to show that  almost surely,  for all divergence-free fields $\varphi \in H^1(\mO)$, 
\begin{equation} \label{eq:vanish:F_eps} 
 |  \langle  n \cdot F^\eps , \varphi   \rangle_{\langle H^{-1/2}(\pa \mO), H^{1/2}(\pa \mO)\rangle} | \le   \nu(\eps) \|\varphi \|_{H^1(\mO)}, \quad \nu(\eps) \rightarrow 0 \text{ as } \eps \rightarrow 0
 \end{equation}
where $F^\eps :=  F_2^\eps +F_3^\eps$.  Notice that $\div(F_2^\eps + F_3^\eps) = 0$.   We introduce again the functions $\chi_\delta$, $\delta > 0$, seen above. We get 
\begin{align*}
& \langle  n \cdot F^\eps , \varphi   \rangle_{\langle H^{-1/2}(\pa \mO), H^{1/2}(\pa \mO)\rangle}  =   \langle  n \cdot  \chi_\delta F^\eps , \varphi   \rangle_{\langle H^{-1/2}(\pa \mO), H^{1/2}(\pa \mO)\rangle} \\
= &   \int_\mO  (\na \chi_\delta \cdot F^\eps)  \cdot  \varphi  -   \int_\mO  \chi_\delta F^\eps \cdot \na \varphi    
\end{align*}
For the last term, we take into account that $\varphi$ is divergence-free, so that the pressure disappears: we find 
$$  |\int_\mO  \chi_\delta F^\eps \cdot \na \varphi|  \le  \left( \|2 \chi_\delta D(H)(\cdot/\eps) \|_{L^2(\mO)}  + \|\chi_\delta F_3(\cdot/\eps) \|_{L^2(\mO)} \right)  \, \|\varphi\|_{H^1(\mO)}.$$
As seen in \eqref{eq:delta_bound}, we have 
$$   \lim_{\eps \rightarrow 0} \|2 \chi_\delta D(H)(\cdot/\eps) \|_{L^2(\mO)}^2 \le C \delta  $$
and similarly, 
$$    \lim_{\eps \rightarrow 0}  \|\chi_\delta F_3^\eps \|_{L^2(\mO)}^2   \le C \delta.  $$
For the first term, we write 
\begin{align*} \int_\mO  (\na \chi_\delta \cdot F^\eps)  \cdot  \varphi  
& =   2 \int_\mO   \na \chi_\delta \cdot \eps^3 D(\bar{h}^\eta_\eps) \cdot  \varphi   -  \int_\mO  (\na \chi_\delta \, \eps^3 \overline{p}_\eps^\eta)  \cdot  \varphi  + \int_\mO  (\na \chi_\delta \cdot F^\eps_3)  \cdot  \varphi  .
\end{align*} 
We know that  $\eps^3 D(\bar{h}^\eta_\eps)$ goes weakly to zero in $L^2(\mO)$, so that it converges strongly to zero in $H^{-1}(\mO)$. As $\na \chi_\delta \otimes \varphi$ belongs to $H^1_0(\mO)$, we find that for a fixed $\delta$
$$ |2 \int_\mO   \na \chi_\delta \cdot \eps^3 D(\bar{h}^\eta_\eps) \cdot  \varphi| \le C \|\eps^3 D(\bar{h}^\eta_\eps) \|_{H^{-1}(\mO)} \|\na  \chi_\delta \varphi\|_{H^1(\mO)} \le \nu(\eps) \|\varphi\|_{H^1(\mO)}.   $$  
 Similarly, as $\E F_3 = 0$, $F_3^\eps$ converges weakly to zero in $L^2(\mO)$ and we get 
 $$|\int_\mO  (\na \chi_\delta \cdot F^\eps_3)  \cdot  \varphi| \le   \nu(\eps) \|\varphi\|_{H^1(\mO)}. $$
 The last step is to prove that $\eps^3 \overline{p}^\eta_\eps$ converges weakly  to zero in $L^2(\mO)$, which will yield 
 $$  |\int_\mO  (\na \chi_\delta \, \eps^3 \overline{p}_\eps^\eta)  \cdot  \varphi| \le  \nu(\eps) \|\varphi\|_{H^1(\mO)}. $$
 As above, for $\phi \in L^2(\mO)$, we introduce $v \in H^1_0(\mO)$ such that $\div v = \phi - \dashint_\mO \phi$, $\|v\|_{H_1(\mO)} \le C_\mO  \|\phi\|_{L^2(\mO)}$. Then, using the equation satisfied by $\overline{p}_\eps^\eta$ in $\mO$: 
$$ - \Delta \eps^3 \bar{h}^\eta_\eps + \na \eps^3  \overline{p}^\eta_\eps  =  \div F_3^\eps $$
 we find  after integration by parts
\begin{align*}  
\int_{\mO} \eps^3 \overline{p}_\eps^\eta \phi  =  \int_{\mO} \eps^3 \overline{p}_\eps^\eta (\phi - \dashint_\mO \phi)  
=  \int_{\mO} \eps^3 \na \bar{h}^\eta_\eps :  \na v  +  \int_{\mO} F_3^\eps : \na v  \xrightarrow[\eps \rightarrow 0]{} 0. 
\end{align*} 
This concludes the proof of \eqref{eq:vanish:F_eps}, of  Proposition \ref{prop:barh_h} and of the theorem. 

\subsection{Formula for periodic point distributions}
Theorem \ref{thm2} gives the limit of $\mV_N$ for  properly rescaled stationary and ergodic point processes, under uniform separation of the points. Such setting includes periodic point distributions, as well as Poisson hard core processes. We focus here on the periodic case, for which further explicit formula can be given.  For   $L > 0$, we consider distinct points $a_1, \dots, a_M$ in $K_L$, and set $\Lambda_0 := \{a_1, \dots, a_M\} + L \Z^d$, which can be seen as a subset of $\T_L^3$.  In  Example \ref{example:periodic}, we explained how to build a  process on $\T_L^3$ out of $\Lambda_0$,  with  $\Lambda(\omega) = \Lambda_0 + \omega$, $\omega \in \T_L^3$. By a simple translation, the results above, that are valid for $\Lambda_0 + \omega$ for a.e. $\omega$, are still valid for $\omega=0$.  Thus, for $\Lambda = \Lambda_0$, we deduce from Proposition \ref{prop:existence:H_eta} the existence of an $L\Z^3$-periodic solution ${\bf H^\eta}$ of \eqref{eq:jellium_eta} with $\na {\bf H^\eta} \in L^2_{loc}$. If we further assume that  ${\bf H^\eta}$ is mean-free, it is clearly unique. Then, following Corollary \ref{cor:existence:H} and Theorem \ref{thm2},  there exists an $L\Z^3$-periodic solution $H$ of \eqref{eq:jellium},  such that 
\begin{equation}
\begin{aligned} \label{lim_VN_periodic}
 & \lim_{N \rightarrow +\infty}  \mV_N = \frac{25L^6}{2M^2} \mathcal{W}(\na H), \quad \mW(\na H) = \lim_{\eta \rightarrow 0} \mW^\eta(\na H), \\
&  \mathcal{W}^\eta(\na H) = -  \left(  \dashint_{K_L}|\na H^\eta|^2  - \frac{M}{L^3\eta^3}  \Bigl( \int_{B^1} |\na G_S^1|^2 + \frac{3}{10\pi} |S|^2 \Bigr)\right) 
\end{aligned}
\end{equation}
 where $H^\eta$  is associated to $H$ by \eqref{def:H_H_eta}.    We have used that in the periodic case, the intensity of the process is $m=\frac{M}{L^3}$, while the expectation is simply the average over $K_L$. 
 
 \mspace
 To make things more explicit,  we introduce the periodic Green function $G_{S,L} : \R^3 \rightarrow \R^3$ satisfying: 
\begin{equation} \label{eq:per_Green}
 -\Delta G_{S,L} + \na p_{S,L} =  S\na \delta_0 \,, \;\; \div G_{S,L} = 0 \: \text{in }  K_L,  \;\; G_{S,L} \: \text{ $L\Z^3$-periodic}, \;\; \int_{K_L} G_{S,L} = 0. 
 \end{equation}
The Green function $G_{S,L}$ is easily expressed in Fourier series. If we write 
$$ G_{S,L}(y) =  \sum_{k \in \Z^3_*} e^{\frac{2i\pi k}{L} \cdot y} \widehat{G}_{S,L}(k)$$
a straightforward calculation shows that for all $k \in \Z^3_*$
$$ \widehat{G}_{S,L}(k) = \frac{i}{2\pi L^2 |k|} \left( S \frac{k}{|k|} - \frac{Sk \cdot k}{|k|^2}  \frac{k}{|k|} \right) = \frac{i}{2\pi L^2 |k|^2} \pi^\perp_kSk $$
where $\pi^\perp_k$ denotes the projection orthogonally to the line $\R k$. Note that the Fourier series for $G_{S,L}$ converges for instance in the quadratic sense. 
\begin{proposition} \label{prop: periodic}
$$\lim_{N \rightarrow +\infty}  \mV_N \: = \: \frac{25 L^3}{2M^2} \, \left(  \sum_{i\neq j \in \{1,\dots,M\}} \!\!\!\!\! S \na \cdot G_{S,L}(a_i - a_j)  + M \lim_{y \rightarrow 0} S\na \cdot (G_{S,L}(y) - G_{S}(y)) \right). $$  
\end{proposition}

\mspace
\noindent
{\em Proof}. Clearly, the  $L\Z^3$-periodic field defined on $K_L$ by  $\tilde H(y) :=  \sum_{i=1}^M G_{S,L}(y+a_i)$  is a solution of \eqref{eq:jellium}, and by Proposition \ref{jellium_uniqueness}  $\na \tilde H$ and  $\na H$ differ from a constant matrix. As $\na (\tilde H - H) = \na (\tilde{H}^\eta - H^\eta)$ is the gradient of a periodic function, we have eventually $\na \tilde H = \na H$. Up to adding a constant field to $H$, we can  assume that 
$$ H(y) = \sum_{i=1}^M G_{S,L}(y+a_i).$$
Then, if $\eta$ is small enough so that $B(a_i, \eta) \subset K_L$ for all $i$,  $H^\eta$ is the $L$-periodic field given on $K_L$  by  
$$ H^\eta(y) = \sum_{i=1}^M \left(  G_{S,L}(y+a_i) +  (G_S^\eta - G_S)(y+a_i) \right). $$ 
We integrate by parts to find 
\begin{align*}
 \frac{1}{L^3} \int_{K_L} |\na H^\eta|^2 & = \frac{1}{L^3} \int_{K_L}  \sum_{i=1}^M  H^\eta dS^\eta(\cdot +a_i) \\ 
 & =  \frac{1}{L^3}  \int_{K_L}  \sum_{i,j } G_{S,L}(\cdot +a_j)  dS^\eta(\cdot  +a_i) + \frac{1}{L^3}  \int_{K_L}  \sum_{i,j } (G_S^\eta - G_S)(\cdot +a_j) dS^\eta(\cdot +a_i) \\
 & =  \frac{1}{L^3}   \sum_{i \neq j}   \int_{K_L} G_{S,L}(\cdot +a_j)  dS^\eta(\cdot +a_i) +   \frac{1}{L^3} \sum_i  \int_{K_L}  G_{S,L}(\cdot +a_i)  dS^\eta(\cdot +a_i) 
 \end{align*}
where we have used that the last term of the second line vanishes identically. We then write $G_{S,L} = G_S + \phi_{S,L}$ with $\phi_{S,L}$ smooth near $0$ to obtain 
\begin{align*}
 \frac{1}{L^3} \int_{K_L} |\na H^\eta|^2 & = \sum_{i \neq j}  \frac{1}{L^3} \int_{K_L} G_{S,L}(\cdot +a_j)  dS^\eta(\cdot +a_i)  \\ 
& +  \frac{1}{L^3} \sum_i  \int_{K_L}  \phi_{S,L}(\cdot +a_i)  dS^\eta(\cdot +a_i) +  \frac{M}{L^3} \int_{\R^3} G_S dS^\eta .
\end{align*}
Combining with Lemma  \ref{lem:int_GA} and \eqref{lim_VN_periodic}, we get
$$ \lim_{N \rightarrow \infty} \mV_N = -\frac{25L^6}{2M^2} \lim_{\eta \rightarrow 0} \Big(  \sum_{i \neq j} \frac{1}{L^3}   \int_{K_L} G_{S,L}(\cdot +a_j)  dS^\eta(\cdot +a_i)  + \frac{1}{L^3} \sum_i  \int_{K_L}  \phi_{S,L}(\cdot +a_i)  dS^\eta(\cdot +a_i) \Big). $$
We conclude by the last point of Lemma \ref{lem:smoothing} that 
$$ \lim_{N \rightarrow \infty} \mV_N = \frac{25 L^3}{2M^2}  \Big( \sum_{i \neq j}  S \na \cdot G_{S,L}(a_i - a_j) \: + \: M  S\na \cdot \phi_{S,L}(0) \Big).  $$

 \begin{proposition}[{\bf Simple cubic lattice}] \label{prop:array}
 In  the special case where $L=M=1$,  we find 
 $$ \lim_{N \rightarrow \infty} \mV_N  = \alpha \sum_i S_{ii}^2 + \beta \sum_{i \neq j} S_{ij}^2  $$
 with $\alpha = \frac{5}{2} (1 - 60 a)$, $\beta =  \frac{5}{2} (1 + 40 a)$, and $a \approx -0,04655$ is defined in \eqref{formula_S2}. 
 \end{proposition}
  
  \mspace
 {\em Proof.} When $M=L=1$, the formula from the last proposition simplifies into
 $\lim_N \mV_N = \frac{25}{2}    S\na \cdot \phi_{S,1}(0)$, with $\phi_{S,1} = G_{S,1} - G_S$. The periodic Green function $G_{S,1}$ was computed using Fourier series in the last paragraph. We found
\begin{align*}
G_{S,1}(y) & =  \sum_{k \in \Z^3_*}   \frac{i}{2\pi |k|}   \left( S \frac{k}{|k|} - \frac{Sk \cdot k}{|k|^2}  \frac{k}{|k|} \right) e^{2i\pi k \cdot y} \\
& = S \na  \left(  \sum_{k \in \Z^3_*}    \frac{1}{4\pi^2 |k|^2}    e^{2i\pi k \cdot y}  \right) \: + \: S : (\na \otimes \na) \na \left( \sum_{k \in \Z^3_*}  \frac{1}{16\pi^4 |k|^4}   e^{2i\pi k \cdot y}  \right)
\end{align*}
We use formulas from \cite{Has}, see also \cite[equations (64)-(65)]{ZuAdBr}: 
$$ \sum \frac{1}{4\pi^2 |k|^2} e^{2i\pi k \cdot y} = \frac{1}{4\pi} \left( \frac{1}{|y|} - c_1 + \frac{2\pi}{3} |y|^2 + O(|y|^4) \right) $$
and 
\begin{equation} \label{formula_S2}
 \sum \frac{1}{16\pi^4 |k|^4} e^{2i\pi k \cdot y} = -\frac{1}{4\pi} \left( \frac{|y|}{2} - c_2 - \frac{c_1}{6}|y|^2 + \frac{\pi}{30} |y|^4 +  a P(y) +  O(|y|^6) \right)
\end{equation}
where $c_1$ and $c_2$ are constants, and 
$$P(y) =  \frac{4\pi}{3} \Big( \frac{5}{8} (y_1^4 + y_2^4 + y_3^4) - \frac{15}{4} (y_1^2 y_2^2 + y_1^2 y_3^2 + y_2^2 y_3^2) + \frac{3}{8} |y|^4 \Big).$$ 
Note that the formula  \eqref{formula_S2} defines implicitly $a$. A numerical computation was carried in \cite{ZuAdBr}, see also \cite{NuKe}, giving $a \approx  -0,04655$. 

\mspace
Inserting in the expression for $G_{S,1}$, we find after a tedious calculation that 
$$ S \na \cdot G_{S,1}(y) =  S \na \cdot  \big( - \frac{3}{8\pi}  \frac{(Sy \cdot y) y}{|y|^5} \big) + \frac{1}{5} |S|^2 - 12 a \sum_i  S_{ii}^2  + 8 a \sum_{i \neq j} |S_{ij}|^2  + O(|y|). $$
Note that to carry out this calculation, we used the fact that $S$ is trace-free, which leads to the identity 
$$ 0 = ( \sum_i  S_{ii})^2  = \sum_i  S_{ii}^2 +  \sum_{i \neq j} S_{ii} S_{jj}. $$
  Moreover, we know from \eqref{def:G_S_p_S} that 
 $$ G_S(y) = -\frac{3}{8\pi}   \frac{(Sy \cdot y) y}{|y|^5}.$$
 We end up with 
$$  S\na \cdot \phi_{S,1}(0) = \frac{1}{5}|S|^2 - 12 a \sum_i  S_{ii}^2  + 8 a  \sum_{i \neq j} |S_{ij}|^2 $$
and the right formula for $\lim_N \mV_N$. 


\subsection{Formula in the stationary case with the $2$-point correlation function}
We consider here  the case of random point processes in $\R^3$ ($X=\R$), such that (P1)-(P2)-(P3) hold. We further assume that the mean density is $m=1$. We assume moreover that this point process admits a $2$-point correlation function, that is a function $\rho_2=\rho_2(x,y) \in L^1_{loc}(\R^3 \times \R^3)$ such that for all bounded set $K$ and all smooth $F$ in a neighborhood of $K$: 
$$ \E \sum_{z \neq z' \in K} F(z,z') = \int_{K \times K}  F(x,y) \rho_2(x,y) dx dy. $$
 As the process is stationary, one can write $\rho_2(x,y) = \rho(x-y)$.  Our goal is to prove the following formula: 
 \begin{proposition} \label{prop:explicit_energy} Almost surely,
\begin{align*} 
\lim_N \mV_N & = \frac{25}{2}  \lim_{L \rightarrow +\infty}   \frac{1}{L^3} \,  \sum_{\substack{z\neq z' \in \Lambda \cap K_{L-1}}} \!\!\! S \na \cdot G_{S,L}(z - z')  \\
&= \frac{25}{2} \lim_{L \rightarrow +\infty} \frac{1}{L^3} \int_{K_{L-1} \times K_{L-1}} S\na \cdot G_{S,L}(z-z') \rho(z-z') dz dz'. 
\end{align*}
where $G_{S,L}$ refers to the $L \Z^3$-periodic Green function introduced in \eqref{eq:per_Green}. 
\end{proposition}
\begin{remark}
We remind that the periodic Green function $G_{S,L}$ has singularities at each point of $L \Z^d$. But as the sum is restricted to points $z,z'$ in $\Lambda \cap K_{L-1}$, $z-z'$ is always away from this set of singularities. In the same way, the integral over $K_{L-1} \times K_{L-1}$ in the second equality is well-defined. Under further assumption on the two-point correlation function $\rho$, one could make sense of the integral over $K_L \times K_L$ and replace the former by the latter. 
\end{remark}
\mspace
{\em Proof.} 
Let $\eta$ small enough so that Proposition \ref{prop:stationarity:Weta} holds. We have 
$$  \mathcal{W}(\na H) = \mathcal{W}^\eta(\na H) = - \left(  \mathbb{E} \int_{K_1}|\na H^\eta|^2  - \frac{1}{\eta^3}  \Bigl( \int_{B^1} |\na G_S^1|^2 + \frac{3}{10\pi} |S|^2 \Bigr)\right) .$$
Let $H_L = \sum_{i=1}^M G_{S,L}(\cdot +a_i)$, where $\{ a_1, \dots, a_M\}  = \Lambda \cap K_{L-1}$. Note that $H_L$ is associated to the point process $\Lambda_L$ obtained by $L \Z^d$-periodization of $\Lambda \cap K_{L-1}$.  We shall prove below that,
\begin{equation} \label{periodic_cell_approx}
  \mathbb{E} \int_{K_1}|\na H^\eta|^2 =   \lim_{L \rightarrow +\infty}  \frac{1}{L^3} \int_{K_L} |\na H^\eta_L|^2 \text{ almost surely} 
\end{equation}
As $\frac{M}{L^3} = \frac{|\Lambda \cap K_L|}{L^3}  \rightarrow 1$ as $L \rightarrow +\infty$, it follows from  \eqref{periodic_cell_approx} that 
\begin{align}  
\label{eq: W_periodic_eta} \mathcal{W}(\na H) & = \lim_{L \rightarrow +\infty} -  \left( \frac{1}{L^3} \int_{K_L} |\na H^\eta_L|^2  - \frac{M}{(\eta L)^3}  \Bigl( \int_{B^1} |\na G_S^1|^2 + \frac{3}{10\pi} |S|^2 \Bigr)\right) \\
\nonumber &  = \lim_{L \rightarrow +\infty} \mW^\eta(\na H_L) = \lim_{L \rightarrow +\infty}  \mW(\na H_L),
\end{align}
where the last equality comes from  Proposition  \ref{prop:stationarity:Weta}. One can apply such proposition because the $L\Z^d$-periodized network $\Lambda_L$ has a  minimal distance between points that is independent of $L$. This is the reason why we used $K_{L-1}$ instead of $K_L$ in the definition of $\Lambda_L$.   Eventually, by Proposition  \ref{prop: periodic},
\begin{align*}
 \lim_{L \rightarrow +\infty} \mW(\na H_L) & =  \lim_{L \rightarrow +\infty} \Big(  \frac{1}{L^3} \sum_{i\neq j \in \{1,\dots,M\}} \!\!\!\!\! S \na \cdot G_{S,L}(a_i - a_j)  + \frac{M}{L^3} \lim_{y \rightarrow 0} S\na \cdot (G_{S,L}(y) - G_{S}(y)) \Big)
 \end{align*}
Using that 
\begin{equation*} 
G_{S,L}(y) = \frac{1}{L^2} G_{S,1}\left(\frac{\cdot}{L}\right), \quad G_{S}(y) = \frac{1}{L^2} G_{S}\left(\frac{\cdot}{L}\right) 
\end{equation*}
we get that 
\begin{align*}
 \frac{M}{L^3} \Big|\lim_{y \rightarrow 0} S\na \cdot (G_{S,L} - G_{S})(y)\Bigr|  & \le C \Big|\lim_{y \rightarrow 0} S\na \cdot (G_{S,L} - G_{S})(y)\Big| \\
 & \le \frac{C'}{L^3}  \Big| \lim_{y \rightarrow 0}  S \na \cdot (G_{S,1} - G_{S})(y/L)\Big| = O(L^{-3}). 
 \end{align*}
 We obtain
 \begin{equation} \label{explicit_formula_1}
\begin{aligned} 
\mathcal{W}(\na H) & =  \lim_{L \rightarrow +\infty}   \frac{1}{L^3} \,  \sum_{i\neq j \in \{1,\dots,M\}} \!\!\!\!\! S \na \cdot G_{S,L}(a_i - a_j) 
\end{aligned}
\end{equation}
This is the first formula of the proposition. To prove the second one, one can go back to formula \eqref{eq: W_periodic_eta} and take the expectation of both sides. The left-hand side, which is deterministic, is of course unchanged. As regards the r.h.s., one can swap the limit in $L$ and the expectation by invoking the dominated convergence theorem. Indeed, both terms $\frac{1}{L^3} \int_{K_L} |\na H^\eta_L|^2$ and  $\frac{M}{(\eta L)^3}  \Bigl( \int_{B^1} |\na G_S^1|^2 + \frac{3}{10\pi} |S|^2 \Bigr)$ are bounded uniformly in $n$ and in the random parameter $\omega$ (but not uniformly on $\eta$): the first term is bounded through a simple energy estimate, while the second one is bounded thanks to the almost sure separation assumption.  


\mspace
The final step is to prove  \eqref{periodic_cell_approx} almost surely. 
We set $\eps := \frac{1}{L}$, and introduce for all $x \in K_1$, 
$$ h^\eta_{\eps}(x) = \frac{1}{\eps^2} H^\eta_L(\frac{x}{\eps}), \quad      p_{\eps}^\eta(x) = \frac{1}{\eps^3} p^\eta_L(\frac{x}{\eps}) $$
and similarly, for all $x \in K_1$,  
\begin{align*}
 \bar{h}^\eta_\eps(x) & =  \frac{1}{\eps^2} H^\eta(\frac{x}{\eps}) - \dashint_{K_1}  \frac{1}{\eps^2} H^\eta(\frac{\cdot}{\eps}),\\
  \overline{p}^\eta_\eps(x) & =  \frac{1}{\eps^3} P^\eta(\frac{x}{\eps}) - \dashint_{K_1}  \frac{1}{\eps^3} P^\eta(\frac{\cdot}{\eps}),  
\end{align*}
where $(H^\eta, P^\eta)$ refers to the field built in  Proposition \ref{prop:existence:H_eta}. Clearly, 
$$ \eps^6 \int_{K_1} |\na h^\eta_{\eps}|^2 = \frac{1}{L^3} \int_{K_L} |\na H_{L}^\eta|^2  $$
while, by the ergodic theorem, one has almost surely:
$$\eps^6 \int_{K_1} |\na \bar{h}^\eta_{\eps}|^2 = \frac{1}{L^3} \int_{K_L} |\na H^\eta|^2 \xrightarrow[\eps \rightarrow 0]{} \E \int_{K_1} |\na H^\eta|^2. $$
It remains to show that 
$$ \eps^6 \int_{K_1} |\na (\bar{h}^\eta_{\eps} - h^\eta_{\eps})|^2 \rightarrow 0 \quad \text{as } \eps \rightarrow 0 $$ 
We notice that the difference $h_\eps = \bar{h}^\eta_\eps - h^\eta_{\eps}$ satisfies the  Stokes equation 
$$ - \Delta h_\eps + \na p_\eps = \frac{1}{\eps^3}\div (R_\eps - R_{\eps,L}), \quad \div h_{\eps} = 0 \quad \text{ in } K_1, $$
where 
$$R_\eps :=  \sum_{z \in \Lambda} \Psi^\eta(x/\eps + z), \quad  \quad R_{\eps,L} :=  \sum_{z \in \Lambda_L} \Psi^\eta(x/\eps + z)$$
and 	where we recall that $\Lambda_L$ is obtained by $L\Z^3$-periodization of $\Lambda \cap K_{L-1}$. Testing against $\eps^6 h_\eps$, we find 
\begin{equation} \label{estimate:heps}
\begin{aligned} 
\eps^6 \int_{K_1} |\na h_\eps|^2 & =  - \int_{K_1} (R_\eps - R_{\eps,L}) \eps^3 \na h_\eps    
   + \int_{\pa K_1}  F_\eps n \cdot \eps^3 (h_\eps - \dashint_{K_1} h_\eps)  \: - \:  \int_{\pa K_1} G_\eps n \cdot \eps^3 h_\eps
  \end{aligned}
  \end{equation}
 where 
\begin{align*}
 F_{\eps}(x)  & :=  \na H^\eta(\frac{x}{\eps}) - P^\eta\big(\frac{x}{\eps}\big) I_d + \int_{K_1} P^\eta\big(\frac{\cdot}{\eps}\big) I_d +  \tilde F(x),  \\
 G_\eps(x) & :=  \na H^\eta_L(\frac{x}{\eps}) -  P^\eta_L(\frac{x}{\eps}) I_d +  \tilde G(x) ,
 \end{align*}
 with 
 \begin{align*}
 \tilde F(x) :=  \sum_{z \in \Lambda} \Psi^\eta(x/\eps+z) - S, \quad 
 \tilde G(x) := \sum_{z \in \Lambda_L} \Psi^\eta(x/\eps+z) - S. 
 \end{align*}
  Note that both $F_\eps$ and $G_\eps$ are divergence-free. 
  
  \mspace
 To handle the first term at the right-hand side of \eqref{estimate:heps}, we notice that 
 $$ \left| \{ z \in \Lambda  \bigtriangleup \Lambda_L,  K_L  \cap B(-z,\eta) \neq \emptyset \} \right| = O(L^2) = O(\eps^{-2}), $$
 resulting in 
 $$  \int_{K_1} (R_\eps - R_{\eps,L}) \eps^3 \na h_\eps  \le C \Big( \eps \int_{\R^3} |\Psi^\eta|^2\Big)^{1/2}    \|\eps^3 \na h_\eps\|_{L^2(K_1)} \le C \eps^{1/2} \|\eps^3 \na h_\eps\|_{L^2(K_1)}. $$
As regards the second term, one proceeds exactly as  in  Paragraph \ref{subsec:conv_VN},  replacing $\mO$ by $K_1$: see the treatment of  $F^2_\eps$ and $F^3_\eps$, defined in  \eqref{def_F2eps} and \eqref{def_F3eps}. One gets in this way  that for all divergence-free $\varphi \in H^1(K_1)$,
$$  |\int_{\pa K_1}  F_\eps  n \cdot \varphi | \le \nu(\eps) \|\na \varphi\|_{L^2(K_1)}, \quad \nu(\eps) \rightarrow 0 \quad \text{as }  \eps \rightarrow 0.  $$
As regards the last term, we take into account the periodicity of $H^\eta_L$ and $\tilde G$ to write 
\begin{align*}  
\int_{\pa K_1} G_\eps n \cdot \eps^3 h_\eps  &  = \int_{\pa K_1} G_\eps n \cdot \eps^3 \bar{h}^\eta_\eps dx. 
\end{align*}
As $\int_{\pa K_1}  \bar{h}^\eta_\eps  \cdot n = 0$, we can introduce a solution $\Phi_\eps$ of 
$$ \div \Phi_\eps = 0 \quad \text{ in } \: K_1, \quad \Phi_\eps\vert_{\pa K_1} =  \eps^3 \bar{h}^\eta_\eps\vert_{\pa K_1}, \quad \|\Phi_\eps\|_{H^1(K_1)} \le C \|\eps^3 \bar{h}^\eta_\eps\vert_{\pa K_1}\|_{H^{1/2}(\pa K_1)}. $$
Proceeding as in  Paragraph \ref{subsec:conv_VN} (replacing $\mO$ by $K_1$), one can show  that $\|\eps^3 \bar{h}^\eta_\eps\|_{H^{1/2}(\pa K_1)}$ goes to zero with $\eps$, and so $\|\Phi^\eps\|_{H^1(K_1)}$ goes to zero as well. Eventually, we write 
\begin{align*} 
|\int_{\pa K_1} G_\eps n \cdot \eps^3 \bar{h}^\eta_\eps dx|  & =  |\int_{K_1} G_\eps \cdot \na \Phi_\eps|  \\
& = |\int_{K_1} \left(  2 D(H^\eta_L)(\cdot/\eps)  + \tilde G\right) \cdot \na \Phi_\eps| \\
& \le C \left( \frac{1}{L^3} \|\na H^\eta_L\|^2_{L^2(K_L)} + \|\Psi^{\eta}\|^2_{L^2} + 1  \right)^{1/2} \|\na \Phi_\eps\|_{L^2} \le C' \|\na \Phi_\eps\|_{L^2}. 
\end{align*}
Hence, we find 
$$ \eps^6 \int_{K_1} |\na h|^2 \le C \left( \eps + \nu(\eps)^2  + \|\na \Phi_\eps\|^2_{L^2} \right) \xrightarrow[\eps \rightarrow 0]{} 0 $$
which concludes the proof.

\section*{Acknowledgements}
We express our gratitude to Sylvia Serfaty for explaining to us her work on Coulomb gases and being a source of fruitful suggestions. We acknowledge  the support of the SingFlows project, grant ANR-18-CE40-0027 of the French National Research Agency (ANR). D. G.-V. acknowledges the support of the Institut Universitaire de France. 
M.H. acknowledges the support of Labex Numev Convention grants ANR-10-LABX-20.

\appendix

\section{Proof of Lemma \ref{Calderon}}
For any open set $U$, we denote $\dashint_U = \frac{1}{|U|} \int_U$. By \eqref{H2}, we have 
$$d := \frac{c}{4} N^{-1/3} \le  \min_{i \neq j} \frac{|x_i - x_j|}{4}.$$ 
We write 
$$A'_i = A'_{i,1} + A'_{i,2} + A'_{i,3}$$
 with 
\begin{align*}
A'_{i,1}  & =   \sum_{j\neq i}  \dashint_{B(x_j,d)}\Big( D(v[A_j])(x_i - x_j) -  D(v[A_j])(x_i - x') \Bigr) dx', \\
A'_{i,2} &   =  \sum_{j\neq i}  \dashint_{B(x_j,d)}  \Big( D(v[A_j])(x_i - x')  - \dashint_{B_i} D(v[A_j])(x - x') dx  \Bigr) dx', \\
A'_{i,3} & = \sum_{j\neq i}  \dashint_{B(x_j,d)}  \dashint_{B_i}  D(v[A_j])(x - x') dx  dx'. 
\end{align*}
Setting $y_i = N^{-1/3} x_i$, using that for $i \neq j$, $|y_i - y_j| \ge \frac{1}{2} (c + |y_i - y_j| ) \ge c$,   
$$ |A'_{i,1}| \le C a^3 \sum_{j\neq i} \frac{d}{|x_i - x_j|^4} |A_j|  \le  C'  \phi \sum_{j} \frac{|A_j|}{(c + |y_i - y_j|)^4} $$
From the  inequality \eqref{general_convolution}, applied with $a_{ij} = \frac{1}{(c + |y_i - y_j|)^4}$ and $b_j = A_j$, we deduce 
$$ \sum_{i} |A'_{i,1}|^q \:  \le \mathcal{C}  \phi^q \sum_{j} |A_j|^q. $$
Similarly, 
$$ |A'_{i,2}| \le C a^3 \sum_{j\neq i} \frac{a}{|x_i - x_j|^4} |A_j|  \le  C'  \phi^{\frac{4}{3}} \sum_{j} \frac{|A_j|}{c + |y_i - y_j|^4} $$
This leads to 
$$ \sum_{i} |A'_{i,2}|^q \: \le \mathcal{C}  \phi^{\frac{4q}{3}} \sum_{j} |A_j|^q $$
The last term is  the most difficult. We follow \cite{HiWu}. Let us remind that 
$$ v[A] =-\frac{5}{2} A : (x \otimes x) \frac{a^3 x}{|x|^5}. $$
Let $\chi_d(x) = \chi(x/d)$ a smooth function that is $0$ in $B(0,d)$, $1$ outside $B(0, 2d)$. Introducing the function $F_A = \sum_{j} A_j 1_{B(x_j,d)}$, using that $d \le \min_{i \neq j} \frac{|x_i - x_j|}{4}$,  we can write that  
\begin{align*}
A'_{i,3} =   \frac{1}{d^3} \int_{B_i} \int_{\R^3} \chi_d(x_i-x')  \mathbf{K}(x-x')F_A(x') dx' dx
\end{align*}
where $\mathbf{K}(x)$ is an endomorphism of the space of symmetric matrices, defined by 
$$ \mathbf{K}(x)A  =  -\frac{5}{2} \Big( \frac{4\pi}{3}\Big)^{-2} D \Big( A : (x \otimes x)  \frac{x}{|x|^5} \Big). $$
We then split $A'_{i,3} = M_i + N_i$, with 
\begin{align*}
M_i & =  \frac{1}{d^3} \int_{B_i} \int_{\R^3} \chi_d(x-x') \mathbf{K}(x-x') F_A(x')  dx' dx, \\ 
N_i & = \frac{1}{d^3} \int_{B_i}   \int_{\R^3} (\chi_d(x_i-x') - \chi_d(x-x'))   \mathbf{K}(x-x')F_A(x') dx'  dx. 
\end{align*} 
By H\"older inequality, 
$$ |M_i|^q \le  \frac{1}{d^{3q}}a^{\frac{3q}{p}}   \| \big( \chi_d \mathbf{K}\bigr) \star F_A \|_{L^q(B_i)}^q  $$
and so 
$$ \sum_i  |M_i|^q \le  \frac{1}{d^{3q}}a^{\frac{3q}{p}}   \| \big( \chi_d\mathbf{K}\bigr) \star F_A  \|_{L^q(\R^3)}^q. $$
The kernel   $\chi_d \mathbf{K}$ enters the framework of the  Calder\'on-Zygmund theorem, see for instance \cite[Chapters 4 and 5]{Met}: for all $1 < q <+\infty$, the operator  $\big( \chi_d\mathbf{K}\bigr) \, \star $  is continuous from $L^q(\R^3)$ to $L^q(\R^3)$, with 
$$ \| \big( \chi_d\mathbf{K}\bigr) \star \|_{\mathcal{L}(L^q, L^q)} \le C_q. $$
We stress that the constant $C_q$ depends only on  $q$, and not on $d$, as can be seen from the rescaling $x' := x'/d$. It follows that 
$$   \sum_i  |M_i|^q  \le \frac{C}{d^{3q}}a^{\frac{3q}{p}}  \|F_A\|_{L^q(\R^3)}^q.  $$
As the balls $B(x_j,d)$ are disjoint, $ |\sum A_j 1_{B(x_j,d)}|^q = \sum |A_j|^q 1_{B(x_j,d)}$, so that $\|F_A\|_{L^q(\R^3)}^q = \frac{4\pi}{3}  \sum |A_j|^q d^3$, and 
$$ \sum_i  |M_i|^q  \le  C'  \left(\frac{a}{d}\right)^{\frac{3q}{p}} \sum_i |A_i|^q \le \mathcal{C} \phi^{\frac{q}{p}} \sum_i |A_i|^q  $$
To bound $N_i$, we notice that for all $x \in B_i$, the support of  $x' \rightarrow \chi_d(x_i-x') - \chi_d(x-x')$ is included in 
$$ \Big( B(x_i, 2d) \cup B(x, 2d) \Bigr) \setminus \Big(  B(x,d) \cap B(x_i,d) \Bigr)   \subset B(x, 2d+a) \setminus B(x,d-a). $$
(remark that by definition of $\phi$, $a$ is less than $d$ for $\phi$ small enough). We get 
\begin{align*}
 |N_i|^q & \le  \frac{1}{d^{3q}}a^{\frac{3q}{p}}   \| \big| 1_{B(0, 2d+a) \setminus B(0,d-a)} \mathbf{K} \big| \star \big| F_A \big| \|_{L^q(B_i)}^q  
 \end{align*}
 so that 
 \begin{align*}
 \sum_i |N_i|^q &  \le   \frac{C}{d^{3q}}a^{\frac{3q}{p}}   \| \big| 1_{B(0, 2d+a) \setminus B(0,d-a)}|x|^{-3} \big| \star \big| F_A \big| \|_{L^q(\R^3)}^q \\
&  \le  \frac{C'}{d^{3q}}a^{\frac{3q}{p}} \big| \ln\big(\frac{2d+a}{d-a}\big) \big|^q  \|F_A\|_{L^q(\R^3)}^q \\
&  \le C''  \left(\frac{a}{d}\right)^{\frac{3q}{p}} \sum_i |A_i|^q \le \mathcal{C} \phi^{\frac{q}{p}} \sum_i |A_i|^q  
  \end{align*}
using that for $\phi \ll 1$, $a \ll d$ and  $\big|\ln\big(\frac{2d+a}{d-a}\big)\big|$ is bounded by an absolute constant.

{\small
 \bibliographystyle{abbrv}
 \bibliography{effective_viscosity_refs}

\begin{thebibliography}{10}

\bibitem{AlBr}
Y.~Almog and H.~Brenner.
\newblock Global homogenization of a dilute suspension of spheres.

\bibitem{MR3025042}
H.~Ammari, P.~Garapon, H.~Kang, and H.~Lee.
\newblock Effective viscosity properties of dilute suspensions of arbitrarily
  shaped particles.
\newblock {\em Asymptot. Anal.}, 80(3-4):189--211, 2012.

\bibitem{MR2410410}
A.~Basson and D.~G\'{e}rard-Varet.
\newblock Wall laws for fluid flows at a boundary with random roughness.
\newblock {\em Comm. Pure Appl. Math.}, 61(7):941--987, 2008.

\bibitem{Batch_book}
G.~Batchelor.
\newblock {\em An introduction to fluid dynamics}.
\newblock Cambridge University Press, 2002.

\bibitem{BaGr1}
G.~Batchelor and J.~Green.
\newblock The determination of the bulk stress in a suspension of spherical
  particles at order $c^2$.
\newblock {\em J. Fluid Mech.}, 56:401--427, 1972.

\bibitem{BaGr2}
G.~Batchelor and J.~Green.
\newblock The hydrodynamic interaction of two small freely moving spheres in a
  linear flow field.
\newblock {\em J. Fluid Mech.}, 56(375-400), 1972.

\bibitem{MR1369834}
A.~Y. Beliaev and S.~M. Kozlov.
\newblock Darcy equation for random porous media.
\newblock {\em Comm. Pure Appl. Math.}, 49(1):1--34, 1996.

\bibitem{Bla}
B.~Blaszczyszyn.
\newblock Lecture notes on random geometric models --- random graphs, point
  processes and stochastic geometry.
\newblock Preprint HAL cel-01654766.

\bibitem{MR3046995}
A.~Borodin and S.~Serfaty.
\newblock Renormalized energy concentration in random matrices.
\newblock {\em Comm. Math. Phys.}, 320(1):199--244, 2013.

\bibitem{MR1688875}
J.-Y. Chemin.
\newblock {\em Perfect incompressible fluids}, volume~14 of {\em Oxford Lecture
  Series in Mathematics and its Applications}.
\newblock The Clarendon Press, Oxford University Press, New York, 1998.
\newblock Translated from the 1995 French original by Isabelle Gallagher and
  Dragos Iftimie.

\bibitem{Cla}
R.~Clausius.
\newblock {\em Die mechanische Behandlung der Elektricit\"at}.
\newblock Vieweg, Braunshweig, 1879.

\bibitem{MR2371524}
D.~J. Daley and D.~Vere-Jones.
\newblock {\em An introduction to the theory of point processes. {V}ol. {II}}.
\newblock Probability and its Applications (New York). Springer, New York,
  second edition, 2008.
\newblock General theory and structure.

\bibitem{MR2398959}
L.~Desvillettes, F.~Golse, and V.~Ricci.
\newblock The mean-field limit for solid particles in a {N}avier-{S}tokes flow.
\newblock {\em J. Stat. Phys.}, 131(5):941--967, 2008.

\bibitem{MR3458165}
M.~Duerinckx and A.~Gloria.
\newblock Analyticity of homogenized coefficients under {B}ernoulli
  perturbations and the {C}lausius-{M}ossotti formulas.
\newblock {\em Arch. Ration. Mech. Anal.}, 220(1):297--361, 2016.

\bibitem{Ein}
A.~Einstein.
\newblock Eine neue {B}estimmung der {M}olek\"uldimensionen.
\newblock {\em Ann. Physik.}, 19:289--306, 1906.

\bibitem{MR1284205}
G.~P. Galdi.
\newblock {\em An introduction to the mathematical theory of the
  {N}avier-{S}tokes equations. {V}ol. {I}}, volume~38 of {\em Springer Tracts
  in Natural Philosophy}.
\newblock Springer-Verlag, New York, 1994.
\newblock Linearized steady problems.

\bibitem{MR1373741}
P.~Gamblin and X.~Saint~Raymond.
\newblock On three-dimensional vortex patches.
\newblock {\em Bull. Soc. Math. France}, 123(3):375--424, 1995.

\bibitem{GuMo}
E.~Guazelli and J.~Morris.
\newblock {\em A Physical Introduction to Suspension Dynamics}.
\newblock Cambridge University Press, 2011.

\bibitem{MR2982744}
B.~M. Haines and A.~L. Mazzucato.
\newblock A proof of {E}instein's effective viscosity for a dilute suspension
  of spheres.
\newblock {\em SIAM J. Math. Anal.}, 44(3):2120--2145, 2012.

\bibitem{Has}
H.~Hasimoto.
\newblock On the periodic fundamental solutions of the stokes equations and
  their application to viscous flow past a cubic array of spheres.
\newblock {\em Journal of Fluid Mechanics}, 5:317--328, 1959.

\bibitem{HiWu}
M.~Hillairet and D.~Wu.
\newblock Effective viscosity of a polydispersed suspension.
\newblock Preprint arXiv:1905.12306, 2019.

\bibitem{Hin}
E.~Hinch.
\newblock An averaged-equation approach to particle interactions in a fluid
  suspension.
\newblock {\em J. Fluid Mech.}, 83:695--720, 1977.

\bibitem{MR3795188}
R.~M. H\"{o}fer.
\newblock Sedimentation of inertialess particles in {S}tokes flows.
\newblock {\em Comm. Math. Phys.}, 360(1):55--101, 2018.

\bibitem{MR3744384}
R.~M. H\"{o}fer and J.~J.~L. Vel\'{a}zquez.
\newblock The method of reflections, homogenization and screening for {P}oisson
  and {S}tokes equations in perforated domains.
\newblock {\em Arch. Ration. Mech. Anal.}, 227(3):1165--1221, 2018.

\bibitem{MR2094523}
P.-E. Jabin and F.~Otto.
\newblock Identification of the dilute regime in particle sedimentation.
\newblock {\em Comm. Math. Phys.}, 250(2):415--432, 2004.

\bibitem{MR1329546}
V.~V. Jikov, S.~M. Kozlov, and O.~A. Ole\u{\i}nik.
\newblock {\em Homogenization of differential operators and integral
  functionals}.
\newblock Springer-Verlag, Berlin, 1994.
\newblock Translated from the Russian by G. A. Yosifian.

\bibitem{KeRu}
J.~Keller and L.~Rubenfeld.
\newblock Extremum principles for slow viscous flows with applications to
  suspensions.
\newblock {\em J . Fluid Mech.}, 30:97--125, 1967.

\bibitem{MR813657}
T.~L\'{e}vy and E.~S\'{a}nchez-Palencia.
\newblock Einstein-like approximation for homogenization with small
  concentration. {II}. {N}avier-{S}tokes equation.
\newblock {\em Nonlinear Anal.}, 9(11):1255--1268, 1985.

\bibitem{Max}
J.~Maxwell.
\newblock {\em A treatise on Electricity and Magnetism}, volume~1.
\newblock Clarendon Press, 1881.

\bibitem{Mec}
A.~Mecherbet.
\newblock Sedimentation of particles in {S}tokes flow.
\newblock arXiv:1806.07795, June 2018.

\bibitem{Met}
G.~M\'etivier.
\newblock Int\'egrales singuli\`eres, cours {DEA}.
\newblock {https://www.math.u-bordeaux.fr/~gmetivie/ISf.pdf}, 1981 (revised
  2005).

\bibitem{Mos}
O.~Mossotti.
\newblock Discussione analitica sul'influenza che l'azione di un mezzo
  dielettrico ha sulla distribuzione dell'elettricit{\`a} alla superficie di
  pi{\`u} corpi elettrici disseminati in esso.
\newblock {\em Mem. Mat. Fis. della Soc. Ital. di Sci. in Modena}, 24:49--74,
  1850.

\bibitem{NiSc}
B.~Niethammer and R.~Schubert.
\newblock A local version of {E}instein's formula for the effective viscosity
  of suspensions.
\newblock arXiv:1903.08554.

\bibitem{NuKe}
K.~Nunan and J.~Keller.
\newblock Effective viscosity of a periodic suspension.
\newblock {\em J. Fluid Mech.}, 142:269--287, 1984.

\bibitem{Obr}
R.~O'brien.
\newblock A method for the calculation of the effective transport properties of
  suspensions of interacting particles.
\newblock {\em J. Fluid Mech.}, 1979.

\bibitem{MR3455593}
N.~Rougerie and S.~Serfaty.
\newblock Higher-dimensional {C}oulomb gases and renormalized energy
  functionals.
\newblock {\em Comm. Pure Appl. Math.}, 69(3):519--605, 2016.

\bibitem{Sai}
N.~Saito.
\newblock Concentration dependence of the viscosity of high polymer solutions.
  i.
\newblock {\em Journal of the Physical Society of Japan}, 5(1):4--8, 1950.

\bibitem{MR813656}
E.~S\'{a}nchez-Palencia.
\newblock Einstein-like approximation for homogenization with small
  concentration. {I}. {E}lliptic problems.
\newblock {\em Nonlinear Anal.}, 9(11):1243--1254, 1985.

\bibitem{MR2945619}
E.~Sandier and S.~Serfaty.
\newblock From the {G}inzburg-{L}andau model to vortex lattice problems.
\newblock {\em Comm. Math. Phys.}, 313(3):635--743, 2012.

\bibitem{MR3353821}
E.~Sandier and S.~Serfaty.
\newblock 2{D} {C}oulomb gases and the renormalized energy.
\newblock {\em Ann. Probab.}, 43(4):2026--2083, 2015.

\bibitem{MR3309890}
S.~Serfaty.
\newblock {\em Coulomb gases and {G}inzburg-{L}andau vortices}.
\newblock Zurich Lectures in Advanced Mathematics. European Mathematical
  Society (EMS), Z\"{u}rich, 2015.

\bibitem{ZuAdBr}
M.~Zuzovsky, P.~Adler, and H.~Brenner.
\newblock Spatially periodic suspensions of convex particles in linear shear
  flows. iii. dilute arrays of spheres suspended in newtonian fluids.
\newblock {\em Physics of Fluids}, 26:1714, 1983.

\end{thebibliography}
}

 \end{document}